\documentclass[10pt]{article} %
\usepackage[accepted]{tmlr}

\usepackage{amsmath,amsthm,epsfig,amssymb,amsbsy}
\usepackage{enumerate}
\usepackage{comment}
\usepackage{algorithm}
\usepackage{algpseudocode}
\usepackage{amsfonts}       %
\usepackage{dsfont}
\usepackage{enumitem}
\usepackage{amssymb}
\usepackage{mathtools}
\usepackage{nicefrac}

\usepackage{subcaption}
\usepackage{graphicx}

\usepackage{pgfplots}
\pgfplotsset{compat=1.18}
\usepackage{subcaption}

\usepackage{pgfplotstable}
\usepackage{multirow}

\usepackage{adjustbox}

\usepackage{tikz}
\usetikzlibrary{decorations.pathmorphing,
  patterns,arrows,positioning, backgrounds, intersections}
\tikzstyle{Nodes}=[circle, draw=black, fill=red!30, line width=0.7pt, minimum size=20pt]
\tikzstyle{arrow} = [line width=1pt,->,>=stealth]
\tikzstyle{axis} = [line width=1pt,->,>=stealth]

\makeatletter
\let\save@mathaccent\mathaccent
\newcommand*\if@single[3]{%
  \setbox0\hbox{${\mathaccent"0362{#1}}^H$}%
  \setbox2\hbox{${\mathaccent"0362{\kern0pt#1}}^H$}%
  \ifdim\ht0=\ht2 #3\else #2\fi
  }
\newcommand*\rel@kern[1]{\kern#1\dimexpr\macc@kerna}
\newcommand*\widebar[1]{\@ifnextchar^{{\wide@bar{#1}{0}}}{\wide@bar{#1}{1}}}
\newcommand*\wide@bar[2]{\if@single{#1}{\wide@bar@{#1}{#2}{1}}{\wide@bar@{#1}{#2}{2}}}
\newcommand*\wide@bar@[3]{%
  \begingroup
  \def\mathaccent##1##2{%
    \let\mathaccent\save@mathaccent
    \if#32 \let\macc@nucleus\first@char \fi
    \setbox\z@\hbox{$\macc@style{\macc@nucleus}_{}$}%
    \setbox\tw@\hbox{$\macc@style{\macc@nucleus}{}_{}$}%
    \dimen@\wd\tw@
    \advance\dimen@-\wd\z@
    \divide\dimen@ 3
    \@tempdima\wd\tw@
    \advance\@tempdima-\scriptspace
    \divide\@tempdima 10
    \advance\dimen@-\@tempdima
    \ifdim\dimen@>\z@ \dimen@0pt\fi
    \rel@kern{0.6}\kern-\dimen@
    \if#31
      \overline{\rel@kern{-0.6}\kern\dimen@\macc@nucleus\rel@kern{0.4}\kern\dimen@}%
      \advance\dimen@0.4\dimexpr\macc@kerna
      \let\final@kern#2%
      \ifdim\dimen@<\z@ \let\final@kern1\fi
      \if\final@kern1 \kern-\dimen@\fi
    \else
      \overline{\rel@kern{-0.6}\kern\dimen@#1}%
    \fi
  }%
  \macc@depth\@ne
  \let\math@bgroup\@empty \let\math@egroup\macc@set@skewchar
  \mathsurround\z@ \frozen@everymath{\mathgroup\macc@group\relax}%
  \macc@set@skewchar\relax
  \let\mathaccentV\macc@nested@a
  \if#31
    \macc@nested@a\relax111{#1}%
  \else
    \def\gobble@till@marker##1\endmarker{}%
    \futurelet\first@char\gobble@till@marker#1\endmarker
    \ifcat\noexpand\first@char A\else
      \def\first@char{}%
    \fi
    \macc@nested@a\relax111{\first@char}%
  \fi
  \endgroup
}
\makeatother

\def\arxiv{false}
\def\true{true}

\usepackage{hyperref}
\usepackage{url}

\DeclareMathOperator*{\argmin}{arg\,min}

\DeclareMathOperator{\prox}{\mathbf{prox}}
\DeclareMathOperator{\proj}{\mathbf{proj}}

\DeclareMathOperator{\dom}{dom}

\DeclareMathOperator{\intr}{int}

\DeclareMathOperator{\cl}{cl}

\DeclareMathOperator{\conv}{conv}
\DeclareMathOperator{\dist}{\mathbf{dist}}

\newcommand{\bR}{\mathbb{R}}

\newcommand{\bN}{\mathbb{N}}

\newcommand{\exR}{\overline{\mathbb{R}}}

\newcommand{\cN}{\mathcal{N}}

\newcommand{\cX}{\mathcal{X}}
\newcommand{\cS}{\mathcal{S}}

\makeatletter
\newcommand{\aprox}[3][\@nil]{%
  \def\tmp{#1}%
   \ifx\tmp\@nnil
       \operatorname{prox}_{#3}^{#2}
    \else
         \operatorname{prox}_{#3}^{#1 \star #2}
    \fi}

\newcommand{\aenv}[3][\@nil]{%
  \def\tmp{#1}%
   \ifx\tmp\@nnil
       \operatorname{env}_{#3}^{#2}
    \else
         \operatorname{env}_{#3}^{#1 \star #2}
    \fi}

\newcommand{\bprox}[3][\@nil]{%
  \def\tmp{#1}%
   \ifx\tmp\@nnil
       \operatorname{bprox}_{#3}^{#2}
    \else
        \operatorname{bprox}_{#3}^{#1 #2}
    \fi}
\makeatother

\DeclareMathOperator*\minimize{minimize}

\usepackage[capitalize]{cleveref}[0.19]

\crefname{section}{section}{sections}
\crefname{subsection}{subsection}{subsections}
\Crefname{section}{Section}{Sections}
\Crefname{subsection}{Subsection}{Subsections}

\crefname{assumption}{Assumption}{Assumptions}
\crefformat{equation}{\textup{#2(#1)#3}}
\crefrangeformat{equation}{\textup{#3(#1)#4--#5(#2)#6}}
\crefmultiformat{equation}{\textup{#2(#1)#3}}{ and \textup{#2(#1)#3}}
{, \textup{#2(#1)#3}}{, and \textup{#2(#1)#3}}
\crefrangemultiformat{equation}{\textup{#3(#1)#4--#5(#2)#6}}%
{ and \textup{#3(#1)#4--#5(#2)#6}}{, \textup{#3(#1)#4--#5(#2)#6}}{, and \textup{#3(#1)#4--#5(#2)#6}}

\Crefname{figure}{Figure}{Figures}

\crefformat{equation}{\textup{#2(#1)#3}}
\crefrangeformat{equation}{\textup{#3(#1)#4--#5(#2)#6}}
\crefmultiformat{equation}{\textup{#2(#1)#3}}{ and \textup{#2(#1)#3}}
{, \textup{#2(#1)#3}}{, and \textup{#2(#1)#3}}
\crefrangemultiformat{equation}{\textup{#3(#1)#4--#5(#2)#6}}%
{ and \textup{#3(#1)#4--#5(#2)#6}}{, \textup{#3(#1)#4--#5(#2)#6}}{, and \textup{#3(#1)#4--#5(#2)#6}}

\Crefformat{equation}{#2Equation~\textup{(#1)}#3}
\Crefrangeformat{equation}{Equations~\textup{#3(#1)#4--#5(#2)#6}}
\Crefmultiformat{equation}{Equations~\textup{#2(#1)#3}}{ and \textup{#2(#1)#3}}
{, \textup{#2(#1)#3}}{, and \textup{#2(#1)#3}}
\Crefrangemultiformat{equation}{Equations~\textup{#3(#1)#4--#5(#2)#6}}%
{ and \textup{#3(#1)#4--#5(#2)#6}}{, \textup{#3(#1)#4--#5(#2)#6}}{, and \textup{#3(#1)#4--#5(#2)#6}}

\newtheorem{theorem}{Theorem}

\newtheorem{lemma}{Lemma}
\newlist{lemenum}{enumerate}{1} %
\setlist[lemenum]{label=(\roman*), ref=\thelemma(\roman*), font=\rm}
\crefalias{lemenumi}{lemma}

\newlist{propenum}{enumerate}{1} %
\setlist[propenum]{label=(\roman*), ref=\theproposition(\roman*), font=\rm}
\crefalias{propenumi}{proposition} 

\newtheorem{definition}{Definition}

\newtheorem{assumption}{Assumption}

\theoremstyle{remark}

\newlist{assumenum}{enumerate}{1} %
\setlist[assumenum]{leftmargin=2.1cm,label=(\roman*), font=\rm}
\creflabelformat{assumenumi}{#2#1#3}
\crefname{assumenumi}{assumption}{assumptions}
\Crefname{assumenumi}{Assumption}{Assumptions}

\title{The inexact power augmented Lagrangian method for constrained nonconvex optimization}

\author{\name Alexander Bodard \email alexander.bodard@kuleuven.be \\
      \addr ESAT-STADIUS \& Leuven.AI,
      KU Leuven
      \AND
      \name Konstantinos Oikonomidis \email konstantinos.oikonomidis@kuleuven.be \\
      \addr ESAT-STADIUS \& Leuven.AI,
      KU Leuven
      \AND
      \name Emanuel Laude \email elaude@proximafusion.com\\
      \addr Proxima Fusion GmbH
      \AND 
      \name Panagiotis Patrinos \email panos.patrinos@kuleuven.be\\
      \addr ESAT-STADIUS \& Leuven.AI,
      KU Leuven
  }

\def\month{10}  %
\def\year{2025} %
\def\openreview{\url{https://openreview.net/forum?id=63ANb4r7EM}} %

\begin{document}

\maketitle

\begin{abstract}
  This work introduces an unconventional inexact augmented Lagrangian method where the augmenting term is a Euclidean norm raised to a power between one and two. The proposed algorithm is applicable to a broad class of constrained nonconvex minimization problems that involve nonlinear equality constraints.
In a first part of this work, we conduct a full complexity analysis of the method under a mild regularity condition, leveraging an accelerated first-order algorithm for solving the Hölder-smooth subproblems. 
Interestingly, this worst-case result indicates that using lower powers for the augmenting term leads to faster constraint satisfaction, albeit with a slower decrease of the dual residual. 
Notably, our analysis does not assume boundedness of the iterates.
Thereafter, we present an inexact proximal point method for solving the weakly-convex and Hölder-smooth subproblems, and demonstrate that the combined scheme attains an improved rate that reduces to the best-known convergence rate whenever the augmenting term is a classical squared Euclidean norm.
Different augmenting terms, involving a lower power, further improve the primal complexity at the cost of the dual complexity.
Finally, numerical experiments validate the practical performance of unconventional augmenting terms.
\end{abstract}

\section{Introduction}

We consider \emph{nonconvex} minimization problems with possibly \emph{nonlinear} equality constraints of the form
\begin{equation} \label{eq:nonconvex_problem} \tag{P}
    \min_{x \in \bR^n} \varphi(x) := f(x) + g(x) \quad \text{s.t.} \quad A(x) = 0.
\end{equation}
Here $f: \bR^n \to \bR$ denotes a nonconvex continuously differentiable function, $g : \bR^n \to \exR$ a proper, lsc and convex function, and $A: \bR^n \to \bR^m$ a continuously differentiable mapping. \textit{Inequality} constraints $B(x) \leq 0$ can be incorporated into \eqref{eq:nonconvex_problem} by introducing slack variables $s \in \bR^m$ subject to $B(x) + s = 0$ and $s \geq 0$.

The general formulation \eqref{eq:nonconvex_problem} incorporates a wide variety of problems arising in areas such as machine learning \citep{kulis_fast_2007,ge_efficient_2016,tepper_clustering_2018} and computer science \citep{ibaraki_resource_1988,zhao_semidefinite_1998}. Typical examples include generalized eigenvalue problems, nonconvex Burer-Monteiro reformulations \citep{burer_nonlinear_2003,burer_local_2005}, and Neyman-Pearson classification \citep{neyman_ix_1997}, while recently this formulation has also been considered for neural network training \citep{sangalli_constrained_2021,evens_neural_2021,liu_inexact_2023}.

This work focuses on \emph{first-order methods} for solving \eqref{eq:nonconvex_problem}, as they scale well with the problem dimensions and require significantly less memory compared to higher-order methods. 
These advantageous properties become increasingly important, particularly in applications that involve large amounts of data.
Whereas existing works usually focus on problems with simple (`proximable') constraints, recent works also analyze first-order methods that handle more complicated constraints, see, e.g., \citep{xu_iteration_2021,li_rate-improved_2021,lin_complexity_2022}.

This paper proposes an unconventional \emph{augmented Lagrangian method} (ALM) to solve \eqref{eq:nonconvex_problem} and analyzes its complexity under a mild regularity condition (cf.\,\cref{assump:regularity}) while taking explicitly into account the inexact solutions of the subproblems.
Notably, our framework is the first to conduct such an analysis in the presence of a generic convex term \(g\) under \cref{assump:regularity} (cf.\,\cref{sec:ipalm-complexity}), without assuming compactness of \(\dom \varphi\) (cf.\,\cref{lem:bounded-level-sets}), and under relaxed smoothness assumptions.

\subsection{Background and motivation}

\paragraph{Penalty and augmented Lagrangian methods} Penalty methods address \eqref{eq:nonconvex_problem}
by solving, for some penalty function $\phi : \bR^m \to \bR$, a sequence of problems of the form
\[
    \minimize_{x \in \bR^n} \varphi(x) + \beta \phi(A(x)),
\]
in which the penalty parameter $\beta > 0$ is gradually increased.
A common penalty function is the squared Euclidean norm, i.e., $\phi(\cdot) = \frac{1}{2} \Vert \cdot \Vert^2$.
However, it is well-known that quadratic penalty methods are empirically outperformed by ALMs.
This last class of methods is based on the augmented Lagrangian (AL) function $L_\beta : \bR^n \times \bR^m \to \exR$, which, for a penalty parameter $\beta > 0$ and squared Euclidean norm $\phi$, 
reads 
\begin{align} \label{eq:power_alm}
    L_\beta(x, y) := \varphi(x) + \langle y, A(x)\rangle + \beta \phi(A(x)).%
\end{align}
The AL involves a penalty term and an additional term $\langle y, A(x) \rangle$, where $y \in \bR^m$ are called the dual variables or multipliers.
Every classical ALM iteration updates the primal-dual pair $(x^k, y^k)$ by alternatingly minimizing $L_\beta$ with respect to the primal variable, and consecutively taking a dual ascent step on $L_\beta$, i.e.,
\begin{equation} \label{eq:alm-inner-problem}
    x^{k+1} \in \argmin_{x \in \bR^n} L_{\beta}(x, y^k), \quad y^{k+1} = y^k + \sigma A(x^{k+1})
\end{equation}
for some dual step size $\sigma \geq 0$.
The primal ALM update can in general only be computed by means of an inner solver, and hence inexactly.
In recent years, various works have taken this inexactness explicitly into account when analyzing the complexity of the method, which is then referred to as \emph{inexact} ALM (iALM), see e.g., \citep{sahin_inexact_2019}.

\paragraph{An unconventional penalty term}
In this work we consider an iALM based on a generalization of the classical AL \eqref{eq:power_alm}, in which the penalty term equals
\begin{equation} \label{eq:phi}
    \phi = \tfrac{1}{\nu+1} \Vert \cdot \Vert^{\nu+1}, \qquad \nu \in (0, 1].
\end{equation}
To distinguish from the classical setup $\nu = 1$, we refer to this as the \emph{power} AL, and corresponding method.

\paragraph{Connection to high-order proximal point methods}
In the convex setting \citep{rockafellar_augmented_1976} showed that classical ALM is equivalent to the \emph{proximal point method} (PPM) applied to the negative Lagrange dual function $\varrho$, in the sense that $\{y^k\}_{k \in \bN}$ generated by 
\begin{equation} \label{eq:ppa}
y^{k+1} = \prox_{\beta \varrho }(y^k) := \argmin_{y \in \bR^m} \varrho(y) + \tfrac{1}{2\beta} \Vert y - y^k \Vert^2
\end{equation}
coincides with the sequence of multipliers defined by \eqref{eq:alm-inner-problem} with dual step size $\sigma = \beta$.
More recently, \citep{oikonomidis_global_2024} showed that in the convex regime a similar dual interpretation exists for the so-called \emph{power} ALM \cref{eq:alm-inner-problem} (with $\phi$ as in \cref{eq:phi}) if the dual update is modified to
\(
y^{k+1} = y^k + \beta \nabla \phi(A(x^{k+1}))
\), where 
\begin{equation}
    \nabla \phi(x) = \tfrac{1}{\|x\|^{1-\nu}}x \label{eq:grad_phi}
\end{equation}
for all $x \neq 0$ and $\nabla \phi(0) = 0$.
More precisely, the convex power ALM corresponds to a \emph{high-order} PPM which is obtained by replacing the quadratic penalty in \cref{eq:ppa} with a higher power $\frac{1}{\beta^p(p+1)}\|\cdot\|_2^{p+1}$, $p\geq 1$ of the Euclidean norm $\|\cdot\|_2$ such that $\frac{1}{\nu+1} + \frac{1}{p+1} = 1$. 
The proposed nonconvex power augmented Lagrangian \eqref{eq:power_alm} reduces to the convex one when $\varphi$ is convex and $A$ is affine \cite[Example 2.1]{oikonomidis_global_2024}.
Although the connection to higher-order methods does not extend to the nonconvex setup considered in this work, it motivates considering augmenting terms $\phi$ as in \eqref{eq:phi}. Moreover, the rapid advancements regarding first-order methods for H\"older-smooth objectives over the past decade \citep{nesterov_universal_2015,lan_bundle-level_2015} potentially allow for faster solutions of the inner problems, which are H\"older smooth in this setting (cf.\,\cref{th:holder-lagrangian}).

\paragraph{Connection to sharp Lagrangians}
The proposed power AL interpolates between a classical AL with \(\nu = 1\) and a \emph{sharp} AL with \(\nu = 0\) \citep[Example 11.58]{rockafellar_variational_1998}.
The latter are appealing because they support exact penalty representations under mild conditions \citep[Example 11.62]{rockafellar_variational_1998}.

\subsection{Related work}

The augmented Lagrangian method, initially introduced in 1969 by \citet{hestenes_multiplier_1969,powell_method_1969}, is a popular algorithm that allows one to cast the constrained problem \eqref{eq:nonconvex_problem} into a sequence of smooth problems. 
Over the years it has been studied extensively, see e.g., the monographs \citep{bertsekas_constrained_1982,birgin_practical_2014}.
ALM is also closely related to the popular Alternating Direction Method of Multipliers \citep{gabay_dual_1976,glowinski_augmented_1989}, which has recently been studied for nonconvex equality constrained problems \citep{cohen_dynamic_2022,el_bourkhissi_linearized_2023}.
In the convex regime, ALM was shown to be equivalent to the proximal point method (PPM) applied to the Lagrangian dual problem \citep{rockafellar_augmented_1976}. 
However, in the nonconvex regime this interpretation is lost at least globally \citep{rockafellar_augmented_2023}. Instead, we assume validity of a mild regularity condition \citep{bourkhissi_complexity_2025} (\cref{assump:regularity}) to establish global convergence of the proposed scheme.
An alternative approach that requires a (close to) feasible initial point, is explored by \citet{grapiglia_complexity_2021,grapiglia_worst-case_2023}.

Inexact ALMs explicitly take into account inexactness in the primal update.
They have been analyzed following the paper \citep{sahin_inexact_2019}, which establishes an $\widetilde O (\varepsilon^{-4})$ complexity for finding an $\varepsilon$-stationary point of \eqref{eq:nonconvex_problem} under a slightly different regularity condition -- see \cref{sec:ipalm-complexity} for details.
More recently, \citet{li_rate-improved_2021} employ a triple-loop scheme, based on an inexact proximal point method to solve the iALM subproblems, and obtain an improved complexity $\widetilde O (\varepsilon^{-3})$ under this regularity condition.
This is the best-known rate of convergence of a first-order method for solving \eqref{eq:nonconvex_problem}.
Moreover, \citet{lu_single-loop_2022} proposes a single-loop primal-dual method that attains the same $\widetilde O (\varepsilon^{-3})$ complexity under a similar condition.
We also mention that \citet{lin_complexity_2022} present a penalty method that attains the same $\widetilde O(\varepsilon^{-3})$ complexity. 

To the best of our knowledge, the `nonlinear' or higher-order proximal point method (using prox-term ${\frac{1}{p+1} \Vert \cdot \Vert^{p+1}}$, $p \geq 1$) was first studied in \citep{luque_nonlinear_1987}, which also analyzed its dual counterpart, the `nonlinear' ALM \citep{luque_nonlinear_1986}.
Besides being restricted to the convex case, these works provide only a local complexity analysis of the outer loop, neglecting the inherent difficulty of the H\"older-smooth subproblems -- a challenge addressed by this paper.
Recently, \citep{nesterov_inexact_2023} studied the joint global complexity of high-order proximal point methods in the convex case; see also \citep{ahookhosh_high-order_2021,ahookhosh_high-order_2024}. Its dual counterpart, the power ALM for convex optimization has been recently studied in \citep{oikonomidis_global_2024}.

Since the power augmented Lagrangian is the sum of a H\"older-smooth function (cf.\,\cref{th:holder-lagrangian}) and a convex function, a first-order method that tackles such problems is desired.
Besides the universal optimal methods by \citet{nesterov_universal_2015,lan_bundle-level_2015} for convex problems, we highlight the accelerated methods by \citet{ghadimi_accelerated_2016,ghadimi_generalized_2019} for the convex and nonconvex setting respectively. 
Interesting developments in this area include the linesearch-free adaptive method by \citet{li_simple_2024} for convex problems, as well as a recent family adaptive proximal-gradient methods by \citet{malitsky_adaptive_2020,latafat_adaptive_2024} which were recently shown to also converge for locally H\"older-smooth problems \citep{oikonomidis_adaptive_2024}.

\subsection{Contributions}
We propose a novel inexact augmented Lagrangian method (iALM) with unconventional powers $\nu \in (0, 1]$ of the augmenting term $\phi$ (cf.\,\eqref{eq:phi}) for solving a general class of nonconvex problems with nonlinear constraints. The case $\nu = 1$ reduces to the standard iALM \citep{sahin_inexact_2019}.
The complexity of the proposed method is analyzed, taking explicitly into account distinct primal and dual tolerances \( \varepsilon_A, \varepsilon_\varphi > 0 \). Notably,
\begin{itemize}
    \item Under a mild regularity condition (\cref{assump:regularity}) we prove convergence to first-order stationary points at a rate of $\omega^{-k}$ and show that the constraint violation decreases at a faster rate of $\omega^{-k / \nu}$, where $k$ is the number of (outer) iterations and $\omega$ determines the rate of increase of the penalty parameters.
    The joint complexity is then analyzed with the accelerated first-order method from \citep{ghadimi_generalized_2019} as an inner solver for the nonconvex and H\"older-smooth subproblems. 
    For $\nu = 1$, we subsume the complexity of \citep{sahin_inexact_2019};
    for \( \nu < 1 \), we obtain \emph{faster constraint satisfaction at the cost of a slower decrease in suboptimality}.
    Thus, the sharper penalties of the AL are reflected in the complexity analysis.
    \item Our complexity analysis improves upon existing results even for \(\nu = 1\) by relaxing the standard Lipschitz smoothness assumption on \(f\) to \emph{local H\"older smoothness}, by \emph{not assuming boundedness of the iterates}, and by including a \emph{generic convex term \(g\)} in the objective.
    \item Under slightly more restrictive assumptions on \(f\) and \(A\), we present a \emph{novel inexact proximal point method} to exploit the structure of the weakly-convex and H\"older-smooth subproblems. 
    Using this inner solver, we further strengthen our complexity analysis to match the best-known $\widetilde O (\varepsilon^{-3})$ rate of convergence for the case $\nu = 1$.
    As before, \( \nu < 1 \) yields a better primal complexity at the cost of a worse dual complexity. 
    The complexity is further improved when \( A(x) \) is a linear mapping.
\end{itemize}
Finally, numerical simulations indicate that unconventional powers \( \nu < 1 \) also perform well in practice.

\subsection{Notation, definitions, and technical assumptions}

We denote the Euclidean inner product and norm on $\bR^n$ by $\langle \cdot, \cdot \rangle$ and $\Vert \cdot \Vert$, respectively.
For matrices $\Vert \cdot \Vert$ denotes the spectral norm.
The Euclidean distance from a point $x$ to a set $\cX$ is denoted by $\dist(x, \cX) = \min_{z \in \cX} \Vert x - z \Vert$.
Given a differentiable mapping $A : \bR^{n} \to \bR^{m}$, we denote its Jacobian at $x$ by $J_{A}(x) \in \bR^{m \times n}$. 
The $\widetilde O$-notation suppresses logarithmic dependencies. 
We denote by $\exR:=\bR \cup \{+\infty\}$ the extended real line.
For an extended-valued function $g:\bR^n \to \exR$ we denote by $\dom g=\{x \in \bR^n \mid g(x) < \infty\}$ its domain and say that $g$ is proper if $\dom g \neq \emptyset$.
A function is lower semi-continuous or simply lsc if its \textit{epigraph} is closed \cite[\S 1B]{rockafellar_variational_1998}.
For a convex function \(g : \bR^n \to \exR\) we denote by \(\partial g(x)\) its (convex) subdifferential at \(x \in \bR^n\).
A function \(\varphi : \bR^n \to \exR\) is called level-coercive if \(\varphi\) is bounded below on bounded sets and satisfies \(\liminf_{\Vert x \Vert \to \infty} \nicefrac{\varphi(x)}{\Vert x \Vert} > 0\).
For a closed, convex set $\cX$ we define for any $x \in \bR^n$ the normal cone of $\cX$ at $x$ as $N_{\cX}(x)=\{v \in \bR^n \mid \langle v, y-x\rangle \leq 0, \forall y \in \cX\}$, if $x \in \cX$, and $N_{\cX}(x)=\emptyset$ otherwise.
We define by $\mathcal{C}^1(\bR^n)$ the class of continuously differentiable functions on $\bR^n$. 
We say that a mapping $F: \bR^n \to \bR^m$ is $H$-H\"older continuous of order $\nu \in (0, 1]$ on $\cX \subseteq \bR^n$ if $\Vert F(x') - F(x)\Vert \leq H\Vert x'-x\Vert^\nu$ for all $x,x' \in \cX$ where $H > 0$.
We say that a differentiable function $f : \bR^n \to \bR$ is $(H, \nu)$-H\"older smooth on $\cX$ if $\nabla f(x)$ is $H$-H\"older-continuous of order $\nu$ on $\cX$.
When omitted, we assume $\cX = \bR^n$, and sometimes we simply write $\nu$-H\"older smooth.
Finally, the function $f$ is $L_f$-Lipschitz smooth on $\cX$ if it is $(L_f, 1)$-H\"older smooth on $\cX$.

\paragraph{Optimality conditions}

If \eqref{eq:nonconvex_problem} has a local minimizer $x \in \bR^n$ satisfying the constraint qualification
\begin{equation} \label{eq:cq} \tag{CQ}
        - J_A^\top(x) y \in N_{\dom g}(x), \qquad
        A(x) = 0 \qquad \Longrightarrow \qquad y = 0,
\end{equation}
then -- this follows by applying the subdifferential chain rule \cite[Theorem 10.6]{rockafellar_variational_1998} in a similar way as in \citet[Exercise 10.7]{rockafellar_variational_1998} -- there exists a vector $y \in \bR^m$ such that
\begin{equation*}
    A(x) = 0, \quad \quad - \nabla f(x) - J_A^\top(x) y \in \partial g(x).
\end{equation*}
This set of (generalized) equations, describing the stationary points of \eqref{eq:nonconvex_problem}, is naturally extended to accommodate approximately stationary points.
In this paper we introduce two distinct tolerances for $\dist \left( - \nabla f(x) - J_A^\top(x) y, \partial g(x) \right)$ and $\Vert A(x) \Vert$, leading to the definition of $(\varepsilon_\varphi, \varepsilon_A)$-stationary points.
We refer to $\varepsilon_A$ as the primal residual tolerance and to $\varepsilon_\varphi$ as the dual residual tolerance.
\begin{definition}[$(\varepsilon_\varphi, \varepsilon_A)$-stationary points] \label{def:stationarity}
    Given $\varepsilon_\varphi \geq 0$ and $\varepsilon_A \geq 0$, a point $x \in \bR^n$ is called an $(\varepsilon_\varphi, \varepsilon_A)$-stationary point of \eqref{eq:nonconvex_problem} if there is a vector $y \in \bR^m$ such that
    \begin{align}
        \Vert A(x) \Vert \leq \varepsilon_A, \qquad \text{and} \qquad
        \dist \left( - \nabla f(x) - J_A^\top(x) y, \partial g(x) \right) \leq \varepsilon_\varphi.
    \end{align}
\end{definition}
Existing analyses for iALMs consider the case \(\varepsilon_\varphi = \varepsilon_A > 0\), as the obtained complexities depend only on \(\min \{\varepsilon_\varphi, \varepsilon_A\}\).
In contrast, our analysis reflects that `sharper' penalties with \(\nu < 1\) yield faster constraint satisfaction (cf.\,\cref{th:complexity-analysis}), which nicely agrees with the exact penalty representation of sharp ALs \cite[Example 11.62]{rockafellar_variational_1998}.
Also in practice, distinct primal and dual tolerances \( \varepsilon_A, \varepsilon_\varphi > 0 \) are relevant.
Indeed, absolute tolerances are sensitive to a rescaling of the objective and constraints, and distinct tolerances can be used to compensate this.
In fact, 
general-purpose solvers typically support distinct primal and dual tolerances for this reason.
Moreover, it is not uncommon to analyze convergence for distinct primal and dual tolerances, see e.g.\,\citep{hermans_qpalm_2022} for a proximal ALM for quadratic programming problems.
Yet, it appears unique that the primal and dual tolerance have a distinct effect on the worst-case complexity. 

\paragraph{Assumptions}

Throughout this work, we make the following assumptions on \(f\), \(g\) and \(A\).
\begin{assumption} \label{assumption:smoothness}
    For any nonempty compact set \(\cS \subseteq \dom g\), there exist positive constants \(H_f, H_A, A_{\max}, {J_A}_{\max}, \nabla f_{\max}\) such that the following statements hold:
    \begin{assumenum}
        \item There exists an \(\nu_f \in (0, 1]\) such that \( \Vert \nabla f(x') - \nabla f(x) \Vert \leq H_f \Vert x' - x \Vert^{\nu_f} \) for all \(x, x' \in \cS\);
        \item There exists an \(\nu_A \in (0, 1]\) such that \( \Vert J_A(x') - J_A(x) \Vert \leq H_A \Vert x' - x \Vert^{\nu_A} \) for all \(x, x' \in \cS\);
        \item \(\Vert \nabla f(x) \Vert \leq \nabla f_{\max} \), \(\Vert A(x) \Vert \leq A_{\max} \) and \(\Vert J_A(x) \Vert \leq {J_A}_{\max} \) for all \( x \in \cS\).
    \end{assumenum}
\end{assumption}
\begin{assumption} \label{assumption:convex-g}
    The function \(g\) is proper, lsc and convex, and \(\dom g\) has nonempty interior.
    For any nonempty compact set \(\cS \subseteq \dom g\), there exists \(G \geq 0\) such that \(\Vert g(x') - g(x) \Vert \leq G \Vert x' - x \Vert\) for all \(x, x' \in \cS\).
\end{assumption}
\begin{assumption} \label{assumption:level-coercive}
    The function \(\varphi \equiv f + g\) is level-coercive.
\end{assumption}
\cref{assumption:smoothness} is very mild: the conditions on \(f\) hold when \(f\) is differentiable with \emph{locally} H\"older continuous gradients, and likewise the conditions on \(A\) hold when \(A\) has \emph{locally} H\"older continuous Jacobians.
Slightly more restrictive smoothness assumptions, i.e., with $\nu_f = \nu_A = 1$, have been used in various related works, see e.g., \cite{sahin_inexact_2019,li_rate-improved_2021,lu_single-loop_2022}.
Also \cref{assumption:convex-g} is very mild: a proper, lsc and convex function \(g\) is Lipschitz continuous on any nonempty compact subset of \(\intr \dom g\) \cite[Theorem 24.7]{rockafellar_convex_1970}.
Therefore, any potential issue arises at relative boundary points of the effective domain, as exemplified by the function \(x \mapsto - \sqrt{x}\).
We highlight in particular that indicators of closed convex sets, and real-valued convex functions all satisfy \cref{assumption:convex-g}. 
Finally, note that \(\varphi\) is level-coercive if either \(f\) or \(g\) is level-coercive while the other is bounded below \cite[Exercise 3.29]{rockafellar_variational_1998}.
Alternatively, it suffices that \(\varphi\) is coercive, as is the case when \(\dom \varphi \equiv \dom g\) is bounded.

Notably, and contrary to related works \citep{sahin_inexact_2019,li_rate-improved_2021,lu_single-loop_2022,bourkhissi_complexity_2025}, no compactness of \(\dom \varphi\) or boundedness of the iterates is assumed, which makes the smoothness assumptions on \(f\) and \(A\) truly local.
Moreover, we incorporate a generic convex term \(g\). Also \cite{sahin_inexact_2019,bourkhissi_complexity_2025} have done this, but, respectively, under a regularity condition which may not be ideal (cf.\,\cref{sec:ipalm-complexity}) or under additional technical assumptions which may be hard to verify in practice.
When analyzing the complexity of power ALM with the proximal point method as inner solver in \cref{sec:ippm}, we restrict our assumptions to match those of \cite{li_rate-improved_2021,lu_single-loop_2022}.
\Cref{table:assumptions} summarizes our assumptions and compares them against those in similar works.
\begin{table}[h]
    \caption{
        Comparison of the assumptions in this work against those in existing works.
    }
    \label{table:assumptions}
    \centering
    
    \begin{adjustbox}{width=\textwidth}
    \setlength\extrarowheight{3pt}

    \pgfplotstableread[row sep=\\,col sep=&]{
    Property  & penalty   & smoothness & composite & boundedness & regularity & remark \\
    \shortstack{\\\citep{sahin_inexact_2019}\\ \vphantom{.}}   & \shortstack{\\$\nu = 1$\\\vphantom{.}}   & \shortstack{\\Lipschitz\\\vphantom{.}}     & \shortstack{\\generic convex\\\vphantom{.}}     & \shortstack{\\assumed\\\vphantom{.}}     & \shortstack{\\$\pm$ \cref{assump:regularity}\\\vphantom{.}} & \shortstack{\vphantom{.}\\RC not ideal for generic \(g\),\\cf.\,\cref{sec:ipalm-complexity}}\\
    \citep{li_rate-improved_2021,lu_single-loop_2022}    & $\nu = 1$      & Lipschitz     & convex indicator   & compact $\dom \varphi$      & \cref{assump:regularity} &  \\
    \citep{bourkhissi_complexity_2025}  & $\nu = 1$  & Lipschitz   & generic nonconvex   & assumed     & \cref{assump:regularity} & + technical conditions  \\
    Ours (Alg.\,\ref{alg:algorithm_nonconvex}) & $\nu \in (0, 1]$ & H\"older (local)      & generic convex   & level-coercivity     & \cref{assump:regularity} &  \\
    Ours (Alg.\,\ref{alg:algorithm_nonconvex} + \ref{alg:algorithm_ippm}) & $\nu \in (0, 1]$   & Lipschitz      & convex indicator   & compact $\dom g$     & \cref{assump:regularity} & \\
    }\datatable

    \pgfplotstabletypeset[
        col sep=&,
        row sep=\\,
        string type,
        header=true,
        columns/Property/.style={string type,column type={l|}, column name={}},
        columns/penalty/.style={string type,column type=c, column name={\shortstack{Penalty\\term}}},
        columns/smoothness/.style={string type,column type=c, column name={\shortstack{Smoothness\\of $f$ and $A$}}},
        columns/composite/.style={string type,column type=c, column name={\shortstack{Nonsmooth\\term $g$}}},
        columns/boundedness/.style={string type,column type=c, column name={\shortstack{Boundedness\\of iterates}}},
        columns/regularity/.style={string type,column type=c, column name={\shortstack{Regularity\\condition (RC)}}},
        columns/remark/.style={string type,column type={c}, column name={\shortstack{Remark\\\vphantom{.}}}},
        every head row/.style={after row=\hline},
    ]\datatable
    \end{adjustbox}

\end{table}

\section{The inexact power augmented Lagrangian method} \label{sec:ipalm}

In this section we present the \emph{inexact power augmented Lagrangian method} (power ALM), of which the pseudocode is shown in \cref{alg:algorithm_nonconvex}.
This method generalizes the inexact augmented Lagrangian method that was proposed and analyzed by \cite{sahin_inexact_2019,li_rate-improved_2021} to settings where $\nu \neq 1$.
To adequately exploit the composite structure of the augmented Lagrangian \(L_\beta(x, y)\), we define
\begin{equation} \label{eq:psi}
    \psi_\beta(x, y) := f(x) + \langle A(x), y \rangle + \frac{\beta}{2} \Vert A(x) \Vert^{\nu+1}
\end{equation}
and note that the augmented Lagrangian \(L_\beta(x, y) = \psi_\beta(x, y) + g(x)\) is the sum of a smooth and a nonsmooth term.
\begin{algorithm}[ht]
\caption{Inexact power augmented Lagrangian method}
\label{alg:algorithm_nonconvex}
\begin{algorithmic}[1]
\Require $x^1 \in \bR^n$, $y^1 \in \bR^m$, $\lambda > 0$, $\omega > 1$, $\sigma_1$, $\beta_1 > 0$, $\nu \in (0, 1]$.
\For{$k=1, 2, \dots$}
   \State \label{alg:x_update}Update the tolerance $\varepsilon_{k+1} = \lambda / \beta_k$ and obtain $x^{k+1} \in \mathcal{X}$ such that
   \begin{equation} \label{eq:alm-inner-problem-condition}
       \dist(-\nabla_x \psi_{\beta_k}(x^{k+1}, y^k), \partial g(x^{k+1})) \leq \varepsilon_{k+1}.
   \end{equation}
   \State \label{eq:dual-step-size}Update the dual step size as \(
       \sigma_{k+1} = \sigma_1 \min(1 , \tfrac{\|A(x^1)\|^\nu \log^2(2)}{\Vert A(x^{k+1}) \Vert^\nu (k+1) \log^2(k+2)}).
   \)
   \State Update the multipliers $y_{k+1} = y_k + \sigma_{k+1}\nabla \phi(A(x^{k+1}))$.\label{eq:y_update_dual}
   \State Update the penalty parameter $\beta_{k+1} = \omega \beta_k$.
\EndFor
\end{algorithmic}
\end{algorithm}
The ALM subproblem in step \ref{alg:x_update} entails \emph{inexactly} minimizing the augmented Lagrangian, i.e., 
\begin{equation} \label{eq:alm-inner-problem-step3}
    \minimize_{x \in \bR^n} L_{\beta_k}(x, y^k) \equiv \psi_{\beta_k}(x, y^k) + g(x),
\end{equation}
until the condition \eqref{eq:alm-inner-problem-condition} is satisfied.
Observe also that the dual step size update rule (step \ref{eq:dual-step-size}) is constructed in a way that ensures boundedness of the multipliers.
This is captured by the following lemma.
\begin{lemma} \label{lem:bounded-multipliers}
    The sequence $\{ y^k \}_{k \in \bN}$ generated by \cref{alg:algorithm_nonconvex} is bounded, i.e., there exists a $y_{\max} \in \bR$, such that 
    \(
       \Vert y^{k} \Vert \leq y_{\max}
    \)
    for all \(k \geq 1\).
\end{lemma}
The proof is given in \cref{sec:ipalm-proofs}.
We remark that the results in this work can be extended to incorporate unbounded sequences of multipliers, by following the proofs of \cite{li_rate-improved_2021}.
The next lemma establishes that the iterates of \cref{alg:algorithm_nonconvex} remain in a compact set as long as the inner solver in step \ref{alg:x_update} monotonically decreases its objective.
Hereby, and contrary to existing works \citep{sahin_inexact_2019,li_rate-improved_2021,bourkhissi_complexity_2025}, no compactness assumption is required on \(\dom g\).
A proof is given in \cref{sec:ipalm-proofs}.

\begin{lemma} \label{lem:bounded-level-sets}
    Suppose that \cref{assumption:level-coercive} holds, and let \( \{ x^k \}_{k \in \bN}, \{ y^k \}_{k \in \bN}, \{ \beta^k \}_{k \in \bN} \) be generated by \cref{alg:algorithm_nonconvex}.
    For any \(\beta \in \bR, y \in \bR^m, \bar x \in \dom g\), the sublevel set
    \(
        \mathcal{L}_{\beta, y}(\bar x) := \{ x \in \bR^n : L_\beta(x, y) \leq L_\beta(\bar x, y) \}
    \)
    is nonempty and compact. %
    If, additionally, \( L_{\beta_k}(x^{k+1}, y^k) \leq L_{\beta_k}(x^k, y^k) \) for \(k \geq 1\), then:
    \begin{lemenum}
        \item \label{lem:bounded-level-sets-1}the sublevel sets satisfy \(\mathcal{L}_{\beta_{k+1}, y^{k+1}}(x^{k+1}) \subseteq \mathcal{L}_{\beta_k, y^k}(x^{k})\) for \(k \geq 1\);
        \item \label{lem:bounded-level-sets-2}the iterates \( \{ x^k \}_{k \in \bN}\) remain in the (compact) initial sublevel set, i.e., 
        \(
            x^{k+1} \in \mathcal{L}_{\beta_1, y^1}(x^1)
        \)
        for \(k \geq 1\).
    \end{lemenum}
\end{lemma}

\subsection{Complexity analysis} \label{sec:ipalm-complexity}
We analyze the computational complexity of \cref{alg:algorithm_nonconvex} under a regularity condition involving the nonlinear mapping $A$ and the normal cone \(N_{\dom g}\), which also naturally arises in the constraint qualification \eqref{eq:cq}.
\begin{assumption}[regularity] \label{assump:regularity}
    For any nonempty compact set \(\cS \subseteq \dom g\), there exists an $R > 0$ such that
    \begin{align} \label{eq:nonlinear_pl} \tag{R}
        \dist(-J_A^\top(x)A(x), N_{\dom g}(x)) \geq R\|A(x)\|, \qquad \text{for all } x \in \cS.
    \end{align}
\end{assumption}
\Cref{assump:regularity} was used by \cite{bourkhissi_complexity_2025} to analyze the complexity of an ALM.
If $g$ has full domain, then \cref{assump:regularity} reduces to a Polyak-Lojasiewicz-inequality on the feasibility problem \(\minimize_x \frac{1}{2} \Vert A(x) \Vert^2\) \citep{sahin_inexact_2019},
and is a consequence of the uniform regularity condition by \citet[Definition 3]{bolte_nonconvex_2018} in the so-called \textit{information zone}.
On the other hand, if \(g = \delta_{\cX}\) is the indicator of a closed convex set \(\cX\), and if the constraint qualification \eqref{eq:cq} -- which is itself a generalization of the Mangasarian-Fromovitz condition \citep{rockafellar_lagrange_1993} -- holds at a point \(\bar x \in \bR^n\), then \eqref{eq:nonlinear_pl} holds for the singleton \(\cS = \{ \bar x \} \).
\Cref{assump:regularity} moreover assumes the existence of a \emph{uniform constant} \(R > 0\) such that \eqref{eq:nonlinear_pl} holds for any \(\cS \in \dom g\).
Existing works on iALMs, such as \citet{li_rate-improved_2021,lu_single-loop_2022}, typically only deal with the case \(g \equiv \delta_{\cX}\), and analyze the corresponding complexity under \cref{assump:regularity}.
\citet{sahin_inexact_2019} use a condition similar to \eqref{eq:nonlinear_pl} involving the subdifferential \(\partial g(x)\) instead of \(N_{\dom g}(x)\).
If \(g = \delta_{\cX}\), then this is equivalent because \(\partial g(x) = N_{\dom g}(x) = N_{\cX}(x)\) for \(x \in \cX\). 
However, if \(g\) has full domain and is strictly continuous, it has been argued by \citet{bourkhissi_complexity_2025} that \cref{assump:regularity} should not involve \(g\), which is the case if \(N_{\dom g}\) is used, but not for \(\partial g\).
We refer to \citet[\S 5]{bourkhissi_complexity_2025} for an extensive discussion and comparison to other regularity conditions.

Various problems satisfy \cref{assump:regularity}, as shown by \cite{sahin_inexact_2019} for \(g = \delta_{\cX}\), including clustering, basis pursuit and others. Moreover, \cite{li_rate-improved_2021} demonstrate that affine equality constrained problems with an additional polyhedral constraint set or a ball constraint set also satisfy this condition.
Some interesting constraint functions \( A \) do not satisfy this condition globally.
For example, \cite[Example 5]{bolte_nonconvex_2018} show that for \(g \equiv 0\) and for spherical constraints \(A\), \eqref{eq:nonlinear_pl} only holds if \(\cS\) is bounded away from the origin.

Our first result describes the number of power ALM iterations to obtain an approximate stationary point.
\begin{theorem}[Outer complexity] \label{th:convergence-outer}
    Let $\{ x^k \}_{k \in \bN}$ denote the sequence of iterates generated by \cref{alg:algorithm_nonconvex}. If \cref{assumption:smoothness,assumption:convex-g,assumption:level-coercive,assump:regularity} hold, and if there exists a nonempty compact set \(\cS \subseteq \dom g\) containing the iterates \(\{ x^k \}_{k \in \bN}\), then $x^{k+1}$ is a $(\frac{Q_f}{\beta_1 \omega^{k-1}}, \frac{Q_A}{\beta_1^\frac{1}{\nu} \omega^\frac{k-1}{\nu}})$-stationary point of \eqref{eq:nonconvex_problem} with
    \begin{align}
        &Q_f := \lambda + {J_A}_{\max} \sigma_1 \tfrac{{\nabla f}_{\max} + G + {J_A}_{\max} y_{\max} + \varepsilon_{1}}{R}, \quad Q_A := \left( \tfrac{{\nabla f}_{\max} + G + {J_A}_{\max} y_{\max} + \varepsilon_{1}}{R} \right)^{\nicefrac{1}{\nu}}.
    \end{align}
\end{theorem}
The proof is given in \cref{sec:ipalm-proofs}.
We highlight that there exists a nonempty compact set \(\cS\) containing the iterates \(\{ x^k \}_{k \in \bN}\) when the inner solver in step \ref{alg:x_update} of \cref{alg:algorithm_nonconvex} monotonically decreases its objective (cf.\,\cref{lem:bounded-level-sets}).
If \(\dom g\) is compact, then this condition is always satisfied, regardless of the inner solver.

\Cref{th:convergence-outer} states that if such a nonempty compact set \(\cS\) exists, then \cref{alg:algorithm_nonconvex} finds a first order stationary point of \eqref{eq:nonconvex_problem} at a rate of $\omega^{-k}$, where \(\omega > 1\) determines the rate of increase of the penalty parameters.
This is identical to the result of \cite{sahin_inexact_2019} for iALM.
Remarkably, in the case of \cref{alg:algorithm_nonconvex}, the constraint violation $\Vert A(x^k) \Vert$ decreases at a \emph{faster rate} of $\omega^{-k/\nu}$.
However, it is important to note that the described rates are only in terms of the number of iterations of power ALM, i.e., the number of calls to the inner solver in step \ref{alg:x_update}. To obtain a full complexity analysis of the method, we must specify the inner solver and require an estimate of its computational cost.
This is related to the H\"older smoothness of \(x \mapsto \psi_\beta(x, y)\).
\begin{lemma}[Augmented Lagrangian smoothness] \label{th:holder-lagrangian}
    Let \cref{assumption:smoothness} hold. Then, for any $y \in \bR^m$, and on any nonempty compact set \(\cS \subseteq \bR^n\), the function $\psi_\beta(\cdot, y)$ as in \eqref{eq:psi} is $(H_\beta, q)$-H\"older smooth for some $H_\beta \geq 0$ which depends on \(\cS\), with $q = \min \left\{ \nu_f, \nu_A, \nu \right\} \in (0, 1]$. In particular, we have for all \(
    x, x' \in \cS
    \)
    that
    \[
        \Vert \nabla_x \psi_\beta(x, y) - \nabla_x \psi_\beta(x', y) \Vert \leq H_\beta \Vert x - x' \Vert^q,
    \]
    where the modulus of H\"older smoothness $H_\beta$ depends on \(H_f, H_A, J_{A_{\max}}\), and thus on \(\cS\), and is given by
    \begin{align*}
        H_{\beta} = \big[ H_f + H_A \Vert y \Vert  + \beta (2^{1-\nu}{J_A}_{\max}^{1+\nu}+A_{\max}^\nu H_A )\big] \max \left\{ 1, D^{1 - q} \right\} \quad \text{ with } \quad D := \sup_{x, x' \in \cS} \Vert x - x' \Vert.
    \end{align*}
\end{lemma}
The proof is given in \cref{sec:ipalm-proofs}.
The ALM subproblems of the form \eqref{eq:alm-inner-problem-step3} in step \ref{alg:x_update} thus involve a composite objective, i.e., the sum of a nonconvex and H\"older-smooth function and a convex function. 
Minimizing such an objective typically results in a higher computational cost than its counterpart with \(\psi_\beta\) Lipschitz-smooth \citep{grimmer_optimal_2024}.
The \emph{unified problem-parameter free accelerated gradient} (UPFAG) method of \citet{ghadimi_generalized_2019} appears one of the only accelerated first-order methods for which a worst-case complexity analysis has been derived under both nonconvexity and H\"older smoothness. We proceed by deriving the total complexity of power ALM when UPFAG is used in step \ref{alg:x_update}. For ease of notation we write
\(
    L_k(x) := L_{\beta_k}(x, y^k).
\)
Remark that the value \(
    L_k^\star := \min_{x \in \bR^n}L_k(x) > - \infty
\) is finite, and hence the inner problems \eqref{eq:alm-inner-problem-step3} in \cref{alg:algorithm_nonconvex} step \ref{alg:x_update} are well-defined.
Indeed, from \cref{assumption:smoothness,assumption:convex-g,lem:bounded-level-sets}, $L_k$ is proper, lsc and level-bounded.
Lower boundedness then follows from \cite[Theorem 1.9]{rockafellar_variational_1998}.
Moreover, since the UPFAG method enforces a monotonic decrease on its objective, \cref{lem:bounded-level-sets} ensures that all UPFAG iterates remain in the compact initial sublevel set $\mathcal{L}_{\beta_1, y^1}(x^1)$.
In the remainder of this section, we therefore define the smoothness constants from \cref{assumption:smoothness,assumption:convex-g} with respect to this set, without further mention.

The following lemma upper bounds the number of UPFAG iterations to inexactly solve an inner problem.
\begin{lemma}[Inner complexity] \label{thm:ghadimi}
    Suppose that \cref{assumption:smoothness,assumption:convex-g,assumption:level-coercive} hold. Then the total number of (inner) iterations performed by the UPFAG method, with a small enough step size and initial iterate $x^k$, to obtain an {$\varepsilon_{k+1}$-stationary} point to the power ALM inner problem (cfr. \cref{alg:algorithm_nonconvex} step \ref{alg:x_update}) is bounded by
    \begin{equation}
        O \left( H_{{\beta_k}}^{1/q} \left[ \tfrac{L_k(x^k) - L_k^\star}{\varepsilon_{k+1}^2} \right]^{\frac{1+q}{2q}} \right).
    \end{equation}
\end{lemma}
This result follows by the complexity result in \cite[Corollary 5]{ghadimi_generalized_2019}.
Since we require a different termination criterion than the original work, we provide an explicit proof in \cref{sec:ipalm-proofs} for completeness.
We can now describe the complexity of power ALM in terms of UPFAG iterations. 
\begin{theorem}[Total complexity] \label{th:complexity-analysis}
    Suppose that \cref{assumption:smoothness,assumption:convex-g,assumption:level-coercive,assump:regularity} hold and the UPFAG method from \cite{ghadimi_generalized_2019} is used for solving \cref{alg:algorithm_nonconvex} step \ref{alg:x_update} in the setting of \cref{thm:ghadimi}. Given $\varepsilon_\varphi > 0$ and $\varepsilon_A > 0$, \cref{alg:algorithm_nonconvex} finds an $(\varepsilon_\varphi, \varepsilon_A)$-stationary point of \eqref{eq:nonconvex_problem} after at most $T = \max \{ T_\varphi, T_A \}$ UPFAG iterations, where for $q = \min \{\nu_f, \nu_A, \nu\}$,
    \begin{equation}
        T_\varphi = \widetilde{O} \left( \varepsilon_\varphi^{-\frac{5 + 3q}{2q}} \right), \qquad
        T_A = \widetilde{O} \left( \varepsilon_A^{-\frac{\nu}{q} - \frac{3\nu(1+q)}{2q}} \right).
    \end{equation}
\end{theorem}
The proof is given in \cref{sec:ipalm-proofs}.
In the Lipschitz smooth setting, our analysis subsumes the one from \cite{sahin_inexact_2019}: if we choose $\nu=1$ and $\varepsilon_\varphi = \varepsilon_A = \varepsilon$, we obtain $T_\varphi = T_A = \widetilde{O}\left( \varepsilon^{-4} \right)$. In the same setting, by choosing $\nu < 1$ we get a worse convergence rate for the dual residual and a better one for the primal. This result is in line with the intuition behind choosing a sharper augmenting term for the AL function, which highly penalizes the constraint violation.
Finally, we highlight that the average number of gradient evaluations per UPFAG iteration is bounded by a constant \citep{ghadimi_generalized_2019}.
It follows that \cref{th:complexity-analysis} also describes the number of first-order oracle calls for finding an $(\varepsilon_\varphi, \varepsilon_A)$-stationary point of \eqref{eq:nonconvex_problem}.

\section{An inexact proximal point inner solver with improved complexity} \label{sec:ippm}

This section presents an inexact proximal point method for solving the inner problems of \cref{alg:algorithm_nonconvex}. The proposed scheme is essentially a double-loop algorithm that uses an accelerated gradient method for computing the proximal point updates, inspired by \cite{kong_complexity_2019,li_rate-improved_2021} and adapted to the H\"older smooth setting of our paper.
However, we emphasize that this extension is not straightforward, since a H\"older-smooth subproblem cannot be made strongly convex by adding a sufficiently large quadratic term, and strong convexity of the inner-most problem is essential in obtaining an improved overall complexity.
Henceforth we restrict \cref{assumption:smoothness,assumption:convex-g,assumption:level-coercive} as follows. 
\begin{assumption} \label{assumption:problem-lipschitz}    
    The function \(g = \delta_{\cX}\) is the indicator of a non-empty, convex and compact set \(\cX \subseteq \bR^n\) with diameter \(D > 0\).
    There exist positive constants \(L_f, L_A, A_{\max}, {J_A}_{\max}, \nabla f_{\max}\) such that:
    \begin{assumenum}
        \item \( \Vert \nabla f(x') - \nabla f(x) \Vert \leq L_f \Vert x' - x \Vert \) for all \(x, x' \in \cX\);
        \item \( \Vert J_A(x') - J_A(x) \Vert \leq L_A \Vert x' - x \Vert \) for all \(x, x' \in \cX\);
        \item \(\Vert \nabla f(x) \Vert \leq \nabla f_{\max} \), \(\Vert A(x) \Vert \leq A_{\max} \) and \(\Vert J_A(x) \Vert \leq {J_A}_{\max} \) for all \( x \in \cX\).
    \end{assumenum}
\end{assumption}
Under \cref{assumption:problem-lipschitz} we have $\nu_f = \nu_A = 1$, and hence by \cref{th:holder-lagrangian} the power AL has H\"older-continuous gradients of order $q = \nu$ on $\cX$.
Moreover, in this setting the power AL function is \emph{weakly-convex}:
\begin{lemma} \label{lem:al-weak-convexity}
    Suppose that \cref{assump:regularity,assumption:problem-lipschitz} hold. Then, for any $y \in \bR^m$ the power augmented Lagrangian $L_\beta(\cdot, y)$ is $\rho$-weakly convex on $\cX$, with 
    \(
        \rho := L_f + L_A (\Vert y \Vert + \beta A_{\max}^\nu).
    \)
\end{lemma}
The proof is given in \cref{sec:ippm-proofs}.
Therefore, every subproblem \eqref{eq:alm-inner-problem-step3} in \cref{alg:algorithm_nonconvex} step \ref{alg:x_update} is of the form
\begin{equation} \label{eq:ippm-problem}
    \min_{x \in \mathcal{X}} \psi(x) := L_\beta(x, y),
\end{equation}
for some \(y \in \bR^m\), where \(\psi\) is $(H_\beta, \nu)$-H\"older-smooth and $\rho$-weakly convex on the compact set $\cX$. The weak convexity of $\psi$ motivates the use of an inexact proximal point method, described in \cref{alg:algorithm_ippm}. Note that the inexact proximal point updates entail minimizing a strongly convex and H\"older smooth function over a compact set and thus we can utilize the Fast Gradient Method (FGM) from \cite{devolder_first-order_2014} to compute them. Our approach differentiates from standard analyses in that the objective function of \eqref{eq:ippm-problem} has qualitatively distinct lower and upper bounds, obtained from the weak-convexity and H\"older-smoothness, respectively. 
To the best of our knowledge, the forthcoming analysis of the proposed inexact proximal point method is the first to exploit this, and as such enables an improved total complexity of \cref{alg:algorithm_nonconvex,alg:algorithm_ippm}.

\begin{algorithm}[ht]
\caption{Inexact proximal point method for \eqref{eq:ippm-problem}}
\label{alg:algorithm_ippm}
\begin{algorithmic}[1]
\Require $x_1 \in \bR^n$, tolerance $\varepsilon > 0$.
\For{$k=1, 2, \dots$}
   \State \label{alg:inner_inner_step}Let $F(\cdot) := \psi(\cdot) + \rho \Vert \cdot - x^k \Vert^2$ \label{alg:def_F} and obtain $x^{k+1} \in \cX$ such that \[
   \quad \dist(-\nabla F(x^{k+1}), N_{\cX}(x^{k+1})) \leq \nicefrac{\varepsilon}{4}.
   \]
   \State \textbf{If} $2 \rho \Vert x^{k+1} - x^k \Vert \leq \frac{\varepsilon}{2}$ \textbf{then} return $x^{k+1}$.
\EndFor
\end{algorithmic}
\end{algorithm}

\subsection{Complexity analysis of the inexact proximal point method}

The computation of a proximal point update, defined in \cref{alg:algorithm_ippm} step \ref{alg:def_F}, involves a problem of the form
\begin{equation} \label{eq:fgm-problem}
    \min_{x \in \cX} F(x)
\end{equation}
where $F : \bR^n \to \bR$ is $(H_F, \nu)$-H\"older-smooth with $H_F := H_\beta + 2\rho \max \left\{ 1, D^{1-q}\right\}$, and $\rho$-strongly convex. We denote the minimizer of $F$ over $\cX$ by $x^\star = \argmin_{x \in \cX}F(x)$. 
The following theorem describes the number of FGM iterations needed to obtain a point satisfying the inequality in \cref{alg:algorithm_ippm} step \ref{alg:inner_inner_step}, and is based on \cite[\S6.2]{devolder_first-order_2014}.
Its proof also handles the different termination criterion that we require compared to \cite[\S6.2]{devolder_first-order_2014}, and for this reason becomes rather technical.
\begin{theorem} \label{thm:inner_inner_complexity}
    Let $F$ be as in \cref{alg:algorithm_ippm} step \ref{alg:def_F} and suppose that \cref{assumption:problem-lipschitz} hold. Then, we need at most $T$ FGM iterations to obtain a point $x^+ \in \cX$ that satisfies
    \(
        \dist(-\nabla F(x^+), N_{\cX}(x^+)) \leq \varepsilon,
    \)
    where
    \[
        T = \widetilde{O} \left( \tfrac{ H_F^{\frac{2}{1+3\nu}}}{\rho^{\frac{\nu+1}{3\nu+1}}} \widebar{H}^{\frac{1+\nu}{\nu}\frac{1-\nu}{1+3\nu}} \varepsilon^{-\frac{1+\nu}{\nu}\frac{1-\nu}{1+3\nu}} \right), \qquad \text{with} \qquad
        \widebar{H} = H_F^{1-\nu}(2 H_F(1+H_F))^{\frac{\nu}{2}} + (2 H_F(1+H_F))^{1/2}.
    \]
\end{theorem}
The proof is given in \cref{sec:ippm-proofs}. 
The next theorem, adapted from \cite[Theorem 1]{li_rate-improved_2021}, provides an upper bound on the number of iterations the inexact proximal point method presented in \cref{alg:algorithm_ippm} needs in order to terminate.
Its proof is also found in \cref{sec:ippm-proofs}.
\begin{theorem}\label{th:ippm-complexity}
    Suppose that \cref{assumption:problem-lipschitz} holds.
    Then, \cref{alg:algorithm_ippm} stops within $T$ iterations, where
    \(
        T = \left\lceil \frac{32 \rho}{\varepsilon^2} (\psi(x^1) - \psi^\star) + 1 \right\rceil,
    \)
    and the output $x^T \in \bR^n$ satisfies $\dist(- \nabla \psi(x^T), N_{\cX}(x^T)) \leq \varepsilon$.
\end{theorem}

\subsection{Joint complexity analysis of Algorithms \ref{alg:algorithm_nonconvex} and \ref{alg:algorithm_ippm}}
Having described the complexity of \cref{alg:algorithm_ippm} we now move on to the total complexity of the joint scheme, which is the main result of this section. We highlight that this result improves upon \cref{th:complexity-analysis}.
For $\nu=1$ it subsumes the complexity result in \cite[Theorem 2]{li_rate-improved_2021}, whereas for \(\nu < 1\) a better primal complexity is obtained at the cost of a worse dual complexity.
A proof is given in \cref{sec:ippm-proofs}.
\begin{theorem}[Total complexity] \label{th:triple-loop-joint-complexity}
    Suppose that \cref{assump:regularity,assumption:problem-lipschitz} hold,
    and let $\{ x^k \}_{k \in \bN}$ denote the iterates of \cref{alg:algorithm_nonconvex}.
    If \cref{alg:algorithm_ippm} is used to solve the subproblems in \cref{alg:algorithm_nonconvex} step \ref{alg:x_update}, and if the inexact proximal point updates in \cref{alg:algorithm_ippm} are computed using FGM, then an $(\varepsilon_\varphi, \varepsilon_A)$-stationary point of \eqref{eq:nonconvex_problem} is obtained after at most $T = \max \{ T_\varphi, T_A \}$ FGM iterations, where
    \begin{align*}
        T_\varphi = \widetilde O \left( \varepsilon_\varphi^{-3 - \frac{1-\nu}{1+3\nu}\left( 1 + \frac{2 (1+\nu)}{\nu} \right)} \right), \qquad T_A = \widetilde O \left( \varepsilon_{A}^{-3\nu - \frac{1-\nu}{1+3\nu}\left( 3\nu + 2 \right)} \right).
    \end{align*}
\end{theorem}
As described in \cref{sec:inner-solvers}, an FGM iteration requires a single gradient evaluation of \(f\) and two projections onto \(\cX\). 
Consequently, \cref{th:triple-loop-joint-complexity} also describes the number of first-order oracle calls needed to find an $(\varepsilon_\varphi, \varepsilon_A)$-stationary point of \eqref{eq:nonconvex_problem}.

\subsection{Improved complexity for linear constraints}
\Cref{th:triple-loop-joint-complexity} describes the joint complexity of the triple-loop version of power ALM: \cref{alg:algorithm_nonconvex} in which the primal update in step \ref{alg:x_update} is obtained through \cref{alg:algorithm_ippm}, and where in turn the inexact proximal point updates are computed using the FGM. 
If the constraint mapping $A$ is \emph{linear}, this worst-case complexity can be further improved.
It is obtained by following the exact same steps as in the proof of \cref{th:triple-loop-joint-complexity}, and by remarking that if $A$ is linear, then $L_A = 0$.
By \cref{lem:al-weak-convexity}, we have $\rho = L_f + L_A(\Vert y \Vert + \beta A_{\max}^\nu)$, and hence $\rho = L_f = O(1)$ if $A$ is linear, instead of $\rho = O(\beta)$ for nonlinear constraints.
The following result subsumes the $\widetilde O \left( \varepsilon^{-\frac{5}{2}}\right)$ complexity of \cite[Theorem 2]{li_rate-improved_2021} for $\nu = 1$. 
A proof is given in \cref{sec:ippm-proofs}.
\begin{theorem}[Total complexity] \label{th:triple-loop-joint-complexity-linear}
    Suppose that the conditions of \cref{th:triple-loop-joint-complexity} hold, and additionally assume that the mapping $A$ is \emph{linear}.
    If \cref{alg:algorithm_ippm} is used to solve the subproblems in \cref{alg:algorithm_nonconvex} step \ref{alg:x_update}, and if the inexact proximal point updates in \cref{alg:algorithm_ippm} are computed using FGM, then an $(\varepsilon_\varphi, \varepsilon_A)$-stationary point of \eqref{eq:nonconvex_problem} is obtained after at most $T = \max \{ T_\varphi, T_A \}$ FGM iterations, where
    \begin{align*}
        T_\varphi = \widetilde O \left( \varepsilon_\varphi^{-2 - 2 \frac{1-\nu}{1+3\nu}\frac{1+\nu}{\nu} - \frac{2}{1+3\nu}} \right), \qquad T_A = \widetilde O \left( \varepsilon_{A}^{-2\nu - 2 \frac{1-\nu}{1+3\nu}(1+\nu) - \frac{2\nu}{1+3\nu}} \right).
    \end{align*}
\end{theorem}

\section{Numerical results}\label{sec:numerics}

In this section, we compare the practical performance of power ALM (\cref{alg:algorithm_nonconvex}) for various choices of the parameter $\nu \in (0, 1]$ to illustrate the empirical behavior of the proposed method in practice.
Recall that the choice $\nu = 1$ reduces to the iALM from \cite{sahin_inexact_2019}, and therefore is always included as a baseline comparison.
All experiments are run in Julia on an HP EliteBook with 16 cores and 32 GB memory, and the source code is publicly available.\footnote{\url{https://github.com/alexanderbodard/tmlr_nonconvex_power_alm}}
Unless mentioned otherwise, the subproblems in \cref{alg:algorithm_nonconvex} step \ref{alg:x_update} are solved using the UPFAG method of \cite{ghadimi_generalized_2019}.
In our experience, this approach works well in practice, and does not usually underperform when compared to the presented triple-loop scheme based on \cref{alg:algorithm_ippm}, despite the improved worst-case convergence rate of the latter approach.
\ifx\arxiv\true
\else
Some additional experiments are included in \cref{sec:numerics-appendix}.
\fi

\subsection{Clustering}

First, we consider a clustering problem which, following \cite{sahin_inexact_2019}, can be reformulated in the form \eqref{eq:nonconvex_problem} through a rank-$r$ Burer-Monteiro relaxation with
\begin{equation*}
    f(x) = \sum_{i = 1}^n \sum_{j = 1}^n D_{i, j} \langle x_i, x_j \rangle, \quad A(x) = \begin{bmatrix}
        x_1^\top \sum_{i = 1}^n x_i - 1, \dots, x_n^\top \sum_{i = 1}^n x_i - 1
    \end{bmatrix}^\top.
\end{equation*}
Here $x := \begin{bmatrix}
    x_1^\top, \dots x_n^\top
\end{bmatrix}^\top \in \bR^{rn}$ with $x_i \in \bR^{r}$ for $i \in [1, n]$, and $\cX$ is the intersection of the nonnegative orthant with the Euclidean ball of radius $\sqrt{s}$.
The scalar $s$ denotes the number of clusters, and $D \in \bR^{n \times n}$ is a distance matrix generated by some given data points $\{ z^i \}_{i = 1}^n$, i.e., such that $D_{i, j} = \Vert z^i - z^j \Vert$.

We test \cref{alg:algorithm_nonconvex} on two problem instances, being the MNIST dataset \cite{deng_mnist_2012} and the Fashion-MNIST dataset \cite{xiao_fashion-mnist_2017}.
The setup is similar to that of \cite{sahin_inexact_2019}, which is in turn based on \cite{mixon_clustering_2016}.
In particular, a simple two-layer neural network was used to first extract features from the data, and then this neural network was applied to $n = 1000$ random test samples from the dataset, yielding the vectors $\{ z^i \}_{i = 1}^{n = 1000}$ that generate the distance matrix $D$.
We define $s = 10, r = 20$, tune $\sigma_1 = 10$, $\lambda = 10^{-3}$, $\beta_1 = 5$, $\omega = 1.1$, and impose a maximum of $N = 1500$ UPFAG iterations per subproblem.
Note that \cite{sahin_inexact_2019} showed that \cref{assump:regularity} is satisfied for this clustering problem.
The results are visualized in \cref{fig:clustering-fmnist,fig:clustering-mnist}.
We observe that the unconventional powers $\nu = 0.7$ and $\nu = 0.8$ perform best, requiring less than half the number of UPFAG iterations to converge to a $10^{-4}$-stationary point when compared to the iALM of \cite{sahin_inexact_2019} ($\nu = 1$).

\begin{figure*}[t!]
    \centering
    \begin{subfigure}[b]{0.32\textwidth}
        \centering
        \resizebox{\textwidth}{!}{
            \begin{tikzpicture}
\begin{axis}[
            width=3in, height=1.8in,
            at={(1.011in,0.642in)},
            legend cell align={left},
            legend pos=north east,
            scale only axis,
            minor grid style={thin,draw opacity=0.3},
            major grid style={thin,draw opacity=0.5},
            grid=both,
            xmin=0,
            xmax=125000,
            ymin={1e-6},
            ymax={100},
            ymode=log,
            xlabel = {\# UPFAG iters},
            ylabel = {$\Vert A(x) \Vert$}
]
    \addplot[color={rgb,1:red,0.0;green,0.6056;blue,0.9787}, name path={d53d4b1e-4fe0-456d-81ce-d39a537157bc}, draw opacity={1.0}, line width={1}, solid, mark={*}, mark size={3.75 pt}, mark repeat={1}, mark options={color={rgb,1:red,0.0;green,0.0;blue,0.0}, draw opacity={1.0}, fill={rgb,1:red,0.0;green,0.6056;blue,0.9787}, fill opacity={1.0}, line width={0.75}, rotate={0}, solid}]
        table[row sep={\\}]
        {
            \\
            0.0  31.5773684301463  \\
            15000.0  0.0626304004192861  \\
            30000.0  0.01341946342946356  \\
            45000.0  0.003198605649951991  \\
            60000.0  0.0010613153158189382  \\
            75000.0  0.00038952151507543157  \\
            90000.0  0.00014769343068072165  \\
            105000.0  5.662382165660462e-5  \\
            120000.0  2.179107322036059e-5  \\
        }
        ;
    \addlegendentry {$\nu =1$}
    \addplot[color={rgb,1:red,0.8889;green,0.4356;blue,0.2781}, name path={a4e6df1a-47be-483c-8a50-f3db6bf19e94}, draw opacity={1.0}, line width={1}, solid, mark={*}, mark size={3.75 pt}, mark repeat={1}, mark options={color={rgb,1:red,0.0;green,0.0;blue,0.0}, draw opacity={1.0}, fill={rgb,1:red,0.8889;green,0.4356;blue,0.2781}, fill opacity={1.0}, line width={0.75}, rotate={0}, solid}]
        table[row sep={\\}]
        {
            \\
            0.0  31.5773684301463  \\
            15000.0  0.035138596389005655  \\
            30000.0  0.0046082612142975795  \\
            45000.0  0.0018171874045153967  \\
            60000.0  0.0005531142360192736  \\
            75000.0  0.00018172298865879743  \\
            90000.0  6.181092111456561e-5  \\
            105000.0  2.129738939797757e-5  \\
        }
        ;
    \addlegendentry {$\nu =0.9$}
    \addplot[color={rgb,1:red,0.2422;green,0.6433;blue,0.3044}, name path={a662e810-cbf3-4d1a-a13b-a56716e61037}, draw opacity={1.0}, line width={1}, solid, mark={*}, mark size={3.75 pt}, mark repeat={1}, mark options={color={rgb,1:red,0.0;green,0.0;blue,0.0}, draw opacity={1.0}, fill={rgb,1:red,0.2422;green,0.6433;blue,0.3044}, fill opacity={1.0}, line width={0.75}, rotate={0}, solid}]
        table[row sep={\\}]
        {
            \\
            0.0  31.5773684301463  \\
            15000.0  0.011375645539408973  \\
            30000.0  0.000697320716201416  \\
            45000.0  0.00011870159601381927  \\
            60000.0  2.9999039262801203e-5  \\
        }
        ;
    \addlegendentry {$\nu =0.8$}
    \addplot[color={rgb,1:red,0.7644;green,0.4441;blue,0.8243}, name path={fd800ba7-bb1e-40b2-8602-5d09668c8238}, draw opacity={1.0}, line width={1}, solid, mark={*}, mark size={3.75 pt}, mark repeat={1}, mark options={color={rgb,1:red,0.0;green,0.0;blue,0.0}, draw opacity={1.0}, fill={rgb,1:red,0.7644;green,0.4441;blue,0.8243}, fill opacity={1.0}, line width={0.75}, rotate={0}, solid}]
        table[row sep={\\}]
        {
            \\
            0.0  31.5773684301463  \\
            15000.0  0.13762479330333952  \\
            30000.0  0.0041970960726778585  \\
            45000.0  0.0018367276070501881  \\
            60000.0  0.0003596710728100801  \\
            75000.0  8.555637698937271e-5  \\
            90000.0  2.1383115266280172e-5  \\
        }
        ;
    \addlegendentry {$\nu =0.7$}
\end{axis}
\end{tikzpicture}
        }
        \captionsetup{justification=centering}
        \caption{Constraint violation\label{fig:clustering-fmnist-const}}
    \end{subfigure}
    \begin{subfigure}[b]{0.32\textwidth}
        \centering
        \resizebox{\textwidth}{!}{
            \begin{tikzpicture}
\begin{axis}[
            width=3in, height=1.8in,
            at={(1.011in,0.642in)},
            legend cell align={left},
            legend pos=north east,
            scale only axis,
            minor grid style={thin,draw opacity=0.3},
            major grid style={thin,draw opacity=0.5},
            grid=both,
            xmin=0,
            xmax=125000,
            ymin={1e-6},
            ymax={100},
            ymode=log,
            xlabel = {\# UPFAG iters},
            ylabel = {$\vert f(x) - f^* \vert$}
]
    \addplot[color={rgb,1:red,0.0;green,0.6056;blue,0.9787}, name path={f6c0723b-63e8-4b9d-a8ba-de6fd3a6189f}, draw opacity={1.0}, line width={1}, solid, mark={*}, mark size={3.75 pt}, mark repeat={1}, mark options={color={rgb,1:red,0.0;green,0.0;blue,0.0}, draw opacity={1.0}, fill={rgb,1:red,0.0;green,0.6056;blue,0.9787}, fill opacity={1.0}, line width={0.75}, rotate={0}, solid}]
        table[row sep={\\}]
        {
            \\
            0.0  56.40490596715807  \\
            15000.0  0.33581365342196534  \\
            30000.0  0.08828859806897071  \\
            45000.0  0.021186248285047782  \\
            60000.0  0.0070405857734456845  \\
            75000.0  0.0025852317225840693  \\
            90000.0  0.0009803512306163498  \\
            105000.0  0.000375826477288399  \\
            120000.0  0.00014458416747231695  \\
        }
        ;
    \addlegendentry {$\nu =1$}
    \addplot[color={rgb,1:red,0.8889;green,0.4356;blue,0.2781}, name path={2472fcef-2c69-4471-bf59-cf16648bb4db}, draw opacity={1.0}, line width={1}, solid, mark={*}, mark size={3.75 pt}, mark repeat={1}, mark options={color={rgb,1:red,0.0;green,0.0;blue,0.0}, draw opacity={1.0}, fill={rgb,1:red,0.8889;green,0.4356;blue,0.2781}, fill opacity={1.0}, line width={0.75}, rotate={0}, solid}]
        table[row sep={\\}]
        {
            \\
            0.0  56.40490596715807  \\
            15000.0  0.14681424228305673  \\
            30000.0  0.012700288741328336  \\
            45000.0  0.011553699701011055  \\
            60000.0  0.0036299609087251383  \\
            75000.0  0.001194407534192976  \\
            90000.0  0.00040629399069302963  \\
            105000.0  0.00013994968549724263  \\
        }
        ;
    \addlegendentry {$\nu =0.9$}
    \addplot[color={rgb,1:red,0.2422;green,0.6433;blue,0.3044}, name path={4f9d00d1-bed2-4738-8c76-8b243a1d6509}, draw opacity={1.0}, line width={1}, solid, mark={*}, mark size={3.75 pt}, mark repeat={1}, mark options={color={rgb,1:red,0.0;green,0.0;blue,0.0}, draw opacity={1.0}, fill={rgb,1:red,0.2422;green,0.6433;blue,0.3044}, fill opacity={1.0}, line width={0.75}, rotate={0}, solid}]
        table[row sep={\\}]
        {
            \\
            0.0  56.40490596715807  \\
            15000.0  0.06715765871150126  \\
            30000.0  0.004606372134176695  \\
            45000.0  0.000788269139540887  \\
            60000.0  0.00019923498033591613  \\
        }
        ;
    \addlegendentry {$\nu =0.8$}
    \addplot[color={rgb,1:red,0.7644;green,0.4441;blue,0.8243}, name path={a51c0b9c-d7a0-4254-862e-046fd3158ddc}, draw opacity={1.0}, line width={1}, solid, mark={*}, mark size={3.75 pt}, mark repeat={1}, mark options={color={rgb,1:red,0.0;green,0.0;blue,0.0}, draw opacity={1.0}, fill={rgb,1:red,0.7644;green,0.4441;blue,0.8243}, fill opacity={1.0}, line width={0.75}, rotate={0}, solid}]
        table[row sep={\\}]
        {
            \\
            0.0  56.40490596715807  \\
            15000.0  0.8942615591811318  \\
            30000.0  0.06335620393465291  \\
            45000.0  0.010578285564342593  \\
            60000.0  0.0023781401852076556  \\
            75000.0  0.0005663590208300207  \\
            90000.0  0.0001415183674495779  \\
        }
        ;
    \addlegendentry {$\nu =0.7$}
\end{axis}
\end{tikzpicture}
        }
        \captionsetup{justification=centering}
        \caption{Suboptimality\label{fig:clustering-fmnist-subopt}}
    \end{subfigure}
    \begin{subfigure}[b]{0.33\textwidth}
        \centering
        \resizebox{\textwidth}{!}{
            \begin{tikzpicture}
\begin{axis}[
            width=3in, height=1.8in,
            at={(1.011in,0.642in)},
            legend cell align={left},
            legend pos=north east,
            scale only axis,
            minor grid style={thin,draw opacity=0.3},
            major grid style={thin,draw opacity=0.5},
            grid=both,
            xmin=0.7,
            xmax=1,
            ymin={60000},
            ymax={150000},
            xlabel = {Power $\nu$},
            ylabel = {\# UPFAG iters}
]
    \addplot[color={rgb,1:red,0.0;green,0.6056;blue,0.9787}, name path={8ce07ea7-2862-4210-be31-1ac32f57fc68}, draw opacity={1.0}, line width={1}, solid, mark={*}, mark size={3.75 pt}, mark repeat={1}, mark options={color={rgb,1:red,0.0;green,0.0;blue,0.0}, draw opacity={1.0}, fill={rgb,1:red,0.0;green,0.6056;blue,0.9787}, fill opacity={1.0}, line width={0.75}, rotate={0}, solid}]
        table[row sep={\\}]
        {
            \\
            0.7  94500.0  \\
            0.75  85500.0  \\
            0.8  69000.0  \\
            0.85  93000.0  \\
            0.9  111000.0  \\
            0.95  111000.0  \\
            1.0  126000.0  \\
        }
        ;
\end{axis}
\end{tikzpicture}
        }
        \captionsetup{justification=centering}
        \caption{\# UPFAG iters, $\varepsilon_\varphi = \varepsilon_A = 10^{-4}$.\label{fig:clustering-fmnist-powers}}
    \end{subfigure}
    \captionsetup{justification=centering}
    \caption{Comparison of power ALM with various powers $\nu$ on solving the clustering problem with Fashion-MNIST data \cite{xiao_fashion-mnist_2017}. The case $\nu = 1$ corresponds to iALM from \cite{sahin_inexact_2019}.\label{fig:clustering-fmnist}}
\end{figure*}
\begin{figure*}[t!]
    \centering
    \begin{subfigure}[b]{0.32\textwidth}
        \centering
        \resizebox{\textwidth}{!}{
            \begin{tikzpicture}
\begin{axis}[
            width=3in, height=1.8in,
            at={(1.011in,0.642in)},
            legend cell align={left},
            legend pos=north east,
            scale only axis,
            minor grid style={thin,draw opacity=0.3},
            major grid style={thin,draw opacity=0.5},
            grid=both,
            xmin=0,
            xmax=100000,
            ymin={1e-6},
            ymax={150},
            ymode=log,
            xlabel = {\# UPFAG iters},
            ylabel = {$\Vert A(x) \Vert$}
]
    \addplot[color={rgb,1:red,0.0;green,0.6056;blue,0.9787}, name path={efe19e29-517e-4e7c-957e-090bd917bb37}, draw opacity={1.0}, line width={1}, solid, mark={*}, mark size={3.75 pt}, mark repeat={1}, mark options={color={rgb,1:red,0.0;green,0.0;blue,0.0}, draw opacity={1.0}, fill={rgb,1:red,0.0;green,0.6056;blue,0.9787}, fill opacity={1.0}, line width={0.75}, rotate={0}, solid}]
        table[row sep={\\}]
        {
            \\
            0.0  31.5773684301463  \\
            15000.0  0.02535329744906394  \\
            30000.0  0.0021587325887705957  \\
            45000.0  0.0005531808297267972  \\
            60000.0  0.0001873325884047701  \\
            75000.0  6.916083240335337e-5  \\
            90000.0  2.626708270761043e-5  \\
        }
        ;
    \addlegendentry {$\nu =1$}
    \addplot[color={rgb,1:red,0.8889;green,0.4356;blue,0.2781}, name path={70857847-3082-40a6-9fb1-6b77f7dcda21}, draw opacity={1.0}, line width={1}, solid, mark={*}, mark size={3.75 pt}, mark repeat={1}, mark options={color={rgb,1:red,0.0;green,0.0;blue,0.0}, draw opacity={1.0}, fill={rgb,1:red,0.8889;green,0.4356;blue,0.2781}, fill opacity={1.0}, line width={0.75}, rotate={0}, solid}]
        table[row sep={\\}]
        {
            \\
            0.0  31.5773684301463  \\
            15000.0  0.0075650713586879964  \\
            30000.0  0.0004557477574340618  \\
            45000.0  9.75685154373063e-5  \\
            60000.0  2.8980454586351726e-5  \\
        }
        ;
    \addlegendentry {$\nu =0.9$}
    \addplot[color={rgb,1:red,0.2422;green,0.6433;blue,0.3044}, name path={f4ecb4bb-02e6-4842-b211-0584e7d2cef6}, draw opacity={1.0}, line width={1}, solid, mark={*}, mark size={3.75 pt}, mark repeat={1}, mark options={color={rgb,1:red,0.0;green,0.0;blue,0.0}, draw opacity={1.0}, fill={rgb,1:red,0.2422;green,0.6433;blue,0.3044}, fill opacity={1.0}, line width={0.75}, rotate={0}, solid}]
        table[row sep={\\}]
        {
            \\
            0.0  31.5773684301463  \\
            15000.0  0.00044769033927113727  \\
            30000.0  1.7616846131151305e-5  \\
        }
        ;
    \addlegendentry {$\nu =0.8$}
    \addplot[color={rgb,1:red,0.7644;green,0.4441;blue,0.8243}, name path={32f8528d-bd47-4c07-98ce-c5bb0f64de1f}, draw opacity={1.0}, line width={1}, solid, mark={*}, mark size={3.75 pt}, mark repeat={1}, mark options={color={rgb,1:red,0.0;green,0.0;blue,0.0}, draw opacity={1.0}, fill={rgb,1:red,0.7644;green,0.4441;blue,0.8243}, fill opacity={1.0}, line width={0.75}, rotate={0}, solid}]
        table[row sep={\\}]
        {
            \\
            0.0  31.5773684301463  \\
            15000.0  6.590321658028746e-5  \\
        }
        ;
    \addlegendentry {$\nu =0.7$}
\end{axis}
\end{tikzpicture}
        }
        \captionsetup{justification=centering}
        \caption{Constraint violation\label{fig:clustering-mnist-const}}
    \end{subfigure}
    \begin{subfigure}[b]{0.32\textwidth}
        \centering
        \resizebox{\textwidth}{!}{
            \begin{tikzpicture}
\begin{axis}[
            width=3in, height=1.8in,
            at={(1.011in,0.642in)},
            legend cell align={left},
            legend pos=north east,
            scale only axis,
            minor grid style={thin,draw opacity=0.3},
            major grid style={thin,draw opacity=0.5},
            grid=both,
            xmin=0,
            xmax=100000,
            ymin={1e-6},
            ymax={150},
            ymode=log,
            xlabel = {\# UPFAG iters},
            ylabel = {$\vert f(x) - f^* \vert$}
]
    \addplot[color={rgb,1:red,0.0;green,0.6056;blue,0.9787}, name path={5ff6dac5-dcf3-42b3-9365-b5b761ba7975}, draw opacity={1.0}, line width={1}, solid, mark={*}, mark size={3.75 pt}, mark repeat={1}, mark options={color={rgb,1:red,0.0;green,0.0;blue,0.0}, draw opacity={1.0}, fill={rgb,1:red,0.0;green,0.6056;blue,0.9787}, fill opacity={1.0}, line width={0.75}, rotate={0}, solid}]
        table[row sep={\\}]
        {
            \\
            0.0  76.27505514448157  \\
            15000.0  0.18679391691780722  \\
            30000.0  0.016424709989152575  \\
            45000.0  0.0042252724320945845  \\
            60000.0  0.0014323527618671505  \\
            75000.0  0.0005289905053160737  \\
            90000.0  0.00020093559417944107  \\
        }
        ;
    \addlegendentry {$\nu =1$}
    \addplot[color={rgb,1:red,0.8889;green,0.4356;blue,0.2781}, name path={28459e12-db3a-46e6-849e-002a0f85328c}, draw opacity={1.0}, line width={1}, solid, mark={*}, mark size={3.75 pt}, mark repeat={1}, mark options={color={rgb,1:red,0.0;green,0.0;blue,0.0}, draw opacity={1.0}, fill={rgb,1:red,0.8889;green,0.4356;blue,0.2781}, fill opacity={1.0}, line width={0.75}, rotate={0}, solid}]
        table[row sep={\\}]
        {
            \\
            0.0  76.27505514448157  \\
            15000.0  0.055230035798260246  \\
            30000.0  0.0033947205912170375  \\
            45000.0  0.0007285527735945152  \\
            60000.0  0.00021654649344782229  \\
        }
        ;
    \addlegendentry {$\nu =0.9$}
    \addplot[color={rgb,1:red,0.2422;green,0.6433;blue,0.3044}, name path={981a6784-fdfd-4ecb-9e36-3c0d26791f33}, draw opacity={1.0}, line width={1}, solid, mark={*}, mark size={3.75 pt}, mark repeat={1}, mark options={color={rgb,1:red,0.0;green,0.0;blue,0.0}, draw opacity={1.0}, fill={rgb,1:red,0.2422;green,0.6433;blue,0.3044}, fill opacity={1.0}, line width={0.75}, rotate={0}, solid}]
        table[row sep={\\}]
        {
            \\
            0.0  76.27505514448157  \\
            15000.0  0.0034972684027394507  \\
            30000.0  0.00013791963527864937  \\
        }
        ;
    \addlegendentry {$\nu =0.8$}
    \addplot[color={rgb,1:red,0.7644;green,0.4441;blue,0.8243}, name path={9d671e36-9a49-4e93-806b-79107fcdd4c1}, draw opacity={1.0}, line width={1}, solid, mark={*}, mark size={3.75 pt}, mark repeat={1}, mark options={color={rgb,1:red,0.0;green,0.0;blue,0.0}, draw opacity={1.0}, fill={rgb,1:red,0.7644;green,0.4441;blue,0.8243}, fill opacity={1.0}, line width={0.75}, rotate={0}, solid}]
        table[row sep={\\}]
        {
            \\
            0.0  76.27505514448157  \\
            15000.0  0.0004898302309612745  \\
        }
        ;
    \addlegendentry {$\nu =0.7$}
\end{axis}
\end{tikzpicture}
        }
        \captionsetup{justification=centering}
        \caption{Suboptimality\label{fig:clustering-mnist-subopt}}
    \end{subfigure}
    \begin{subfigure}[b]{0.33\textwidth}
        \centering
        \resizebox{\textwidth}{!}{
            \begin{tikzpicture}
\begin{axis}[
            width=3in, height=1.8in,
            at={(1.011in,0.642in)},
            legend cell align={left},
            legend pos=north east,
            scale only axis,
            minor grid style={thin,draw opacity=0.3},
            major grid style={thin,draw opacity=0.5},
            grid=both,
            xmin=0.65,
            xmax=1,
            ymin={15000},
            ymax={150000},
            xlabel = {Power $\nu$},
            ylabel = {\# UPFAG iters}
]
    \addplot[color={rgb,1:red,0.0;green,0.6056;blue,0.9787}, name path={7ea84435-6373-44b4-a22e-18c869d443dd}, draw opacity={1.0}, line width={1}, solid, mark={*}, mark size={3.75 pt}, mark repeat={1}, mark options={color={rgb,1:red,0.0;green,0.0;blue,0.0}, draw opacity={1.0}, fill={rgb,1:red,0.0;green,0.6056;blue,0.9787}, fill opacity={1.0}, line width={0.75}, rotate={0}, solid}]
        table[row sep={\\}]
        {
            \\
            0.65  57124.0  \\
            0.7  21000.0  \\
            0.75  24000.0  \\
            0.8  33000.0  \\
            0.85  57000.0  \\
            0.9  70500.0  \\
            0.95  85500.0  \\
            1.0  102000.0  \\
        }
        ;
\end{axis}
\end{tikzpicture}
        }
        \captionsetup{justification=centering}
        \caption{\# UPFAG iters, $\varepsilon_\varphi = \varepsilon_A = 10^{-4}$.\label{fig:clustering-mnist-powers}}
    \end{subfigure}
    \captionsetup{justification=centering}
    \caption{Comparison of power ALM with various powers $\nu$ on solving the clustering problem with MNIST data \cite{deng_mnist_2012}. The case $\nu = 1$ corresponds to iALM from \cite{sahin_inexact_2019}.\label{fig:clustering-mnist}}
\end{figure*}

\subsection{Quadratic programs} \label{sec:exp_quadratics}

Second, we consider nonconvex quadratic programs (QPs) of the form \eqref{eq:nonconvex_problem} with
\begin{equation*}
    \begin{aligned}
        f(x) = \frac{1}{2} x^\top Q x + \langle q, x \rangle, \quad \cX = \left\{ x \mid \underbar x \leq x \leq \bar x \right\}, \quad A(x) = C x - b.
    \end{aligned}
\end{equation*}
Here $Q \in \bR^{n \times n}$, $C \in \bR^{m \times n}$, $b \in \bR^m$, and $q \in \bR^n$, with $Q$ symmetric and indefinite, are randomly generated, and $\underbar x^i = -5$ and $\bar x^i = 5$ for $i \in \bN_{[1, n]}$.
We sample the entries of a diagonal matrix $\Lambda \in \bR^{n \times n}$ from a Gaussian $\cN(0, 50)$, and the entries of a matrix $\Sigma \in \bR^{n \times n}$ from the standard Gaussian. Then, we normalize $\widehat \Sigma = \Sigma / \Vert \Sigma \Vert$ and define $Q := \widehat \Sigma^\top \Lambda \widehat \Sigma$.
The entries of $q$ are sampled from a Gaussian $\cN(0, 2)$, and the entries of $C$ are sampled from the standard Gaussian.
We construct a vector $\mu \in \bR^n$ with standard Gaussian entries, and define $b := C \mu$.
A maximum of $N = 10^4$ UPFAG iterations per subproblem is imposed, and the parameters $\beta_1 = 5, \omega = 1.1$, $\lambda = 10^{-3}$, $\sigma_1 = 10$ follow the tuning from the previous subsection.
We use tolerances $\varepsilon_\varphi = \varepsilon_A = 10^{-3}$.
Remark that the AL function with penalty $\beta_k$ has $\Vert Q + \beta_k C^\top C \Vert$-Lipschitz continuous gradients for $\nu = 1$.
We highlight that the dual residual can be efficiently computed by following the procedure described by \cite[Eq. 2.7 and below]{nedelcu_computational_2014}.

We generate $100$ random QP instances of size $n = 200, m = 10$, and visualize the total number of UPFAG iterations and the objective value at the returned points by means of violin plots in \cref{fig:qp-statistics}.
We observe that $\nu = 0.8$ and $\nu = 0.9$ perform best in terms of UPFAG iterations.
Yet, we also note the more robust performance of $\nu = 0.6$.
Regarding the quality of the (local) solutions, \cref{fig:qp-statistics-objective} indicates that the use of unconventional powers $\nu$ does not yield significantly better or worse solutions than the classical choice $\nu = 1$.

\begin{figure}[h]
    \centering
    \begin{subfigure}[b]{0.40\textwidth}
        \centering
        \resizebox{\textwidth}{!}{
            \input{figures/quadratics/quadratics_gaussian_20_100_double_powers_confidence.tex}
        }
        \captionsetup{justification=centering}
        \caption{Number of UPFAG iterations.\label{fig:qp-statistics-gradients}}
    \end{subfigure}
    \begin{subfigure}[b]{0.45\textwidth}
        \centering
        \resizebox{\textwidth}{!}{
            \input{figures/quadratics/quadratics_gaussian_20_100_double_fs_confidence.tex}
        }
        \captionsetup{justification=centering}
        \caption{Objective value.\label{fig:qp-statistics-objective}}
    \end{subfigure}
    \captionsetup{justification=centering}
    \caption{Power ALM for various powers $\nu$ on solving $100$ random QPs of size $n = 100, m = 20$.\label{fig:qp-statistics}}
\end{figure}

\Cref{table:qp_gaussian_20_100_double} reports the number of UPFAG iterations performed by power ALM for the first $10$ randomly generated QPs, and also includes
the obtained primal and dual residuals.
The power $\nu = 0.8$ performs well compared to the classical setup $\nu = 1$, and requires roughly $25\%$ fewer inner iterations on average.
Note that the best primal residuals are attained by the smaller powers $\nu$, and the best dual residuals by the higher powers $\nu$.
For \(\nu = 0.6\) and \(\nu = 0.7\), the inner solver often attains the maximum number of inner iterations $N$. 
In our experience, further increasing $N$ for these powers works counterproductively and increases the total number of UPFAG iterations.

\begin{table}[h]
    \caption{
        Power ALM with UPFAG inner solver on solving random QPs of size $n = 100, m = 20$.
    }
    \label{table:qp_gaussian_20_100_double}
    \centering
    
    \begin{adjustbox}{width=\textwidth}
    \setlength\extrarowheight{3pt}
    \pgfplotstabletypeset[%
        begin table={\begin{tabular}[t]},
        every head row/.style={
            before row={%
              \hline
              \vphantom{$q = 1$}\\
              \hline
            },
            after row/.add={}{\hline},
        },
        header=true,
        col sep=&,
        row sep=\\,
        string type,
        columns/{trial}/.style ={column name={trial}, column type={|c}},
        every row no 9/.style={after row=\hline},
        every row no 10/.style={after row=\hline},
    ]{
        \\
        trial\\
        1\\
        2\\
        3\\
        4\\
        5\\
        6\\
        7\\
        8\\
        9\\
        10\\
        avg.\\
    }%
    \begin{filecontents*}{data/quadratics_gaussian_20_100_double.csv}
0.00037227868304198346	7.478891803353912e-7	75063.0	0.0005849418354300677	1.3094019512466775e-6	84362.0	0.0009758929249597895	1.8875120157495386e-6	64649.0	0.0005583997409502623	2.723630291428005e-6	48409.0	0.0007249227679283285	1.0134411847933205e-6	58903.0
0.000190972611522921	7.974188952542509e-7	75740.0	0.0006903120496581026	6.272392740184399e-6	80000.0	0.0008968653953577771	7.145439552418387e-7	60237.0	0.0005392264428549191	1.4317819783978876e-6	55371.0	0.0006919200603217987	9.055945138016989e-7	68320.0
0.0009913297361740323	6.323770853273152e-7	132470.0	0.0007955852922586492	5.908331673271621e-7	181724.0	0.0009430266009265595	8.271142920069334e-7	118910.0	0.0008334930578368042	6.025844229762794e-7	103865.0	0.0008671540409815852	6.503496655014234e-7	134378.0
0.0007461325899240209	3.4281911199221357e-6	50000.0	0.0008409715415897128	2.9630757657173744e-6	60000.0	0.0009389009602140967	1.5921610740897089e-6	69593.0	0.0009259786876941431	1.0857431718454882e-6	66592.0	0.0008788218160964451	2.462087962361703e-6	106718.0
8.611808709443407e-5	4.3727048078473375e-5	50000.0	0.0006437905857422366	4.2277290979108743e-5	50000.0	0.0001120797976764857	4.673196880396867e-6	42728.0	0.0008608969040864233	3.3404212857975593e-6	163131.0	0.000974754867264477	6.976238921458349e-7	205252.0
0.00028358152147622687	1.4778865063771156e-6	55738.0	0.0009724236263398919	1.3708034553903964e-6	57214.0	0.0004068635902275098	1.2947640267391847e-6	51055.0	0.0006871767725127374	1.4260995614995368e-6	45508.0	0.0008003346027920729	1.3681441054385577e-6	51917.0
0.0009462685527833731	9.62982927710299e-5	40000.0	0.000467948391446747	7.297074537661242e-6	50000.0	0.0002343530222465531	3.165118753300054e-5	41242.0	0.0009387765127724247	1.420795697619545e-6	37795.0	0.0005879193899365085	3.918599980835053e-6	42767.0
0.0005518090604425221	1.4537069951316203e-5	60000.0	0.0005728592556708478	1.4138179493249134e-5	70000.0	0.0007644142688267623	1.3941754036451814e-6	68184.0	0.0005803926608887695	1.054206805948099e-6	74075.0	0.000678660777320158	8.153288565091946e-7	87036.0
0.00047452419325497906	8.918024767850257e-5	50000.0	0.0003341675474360145	1.3555075960885842e-6	60000.0	0.00038103575772291984	1.8854431984788317e-6	48242.0	0.0004545453722066571	1.1566928685567523e-6	42656.0	0.0009795085730054049	7.25883835129634e-6	44826.0
0.0005592402130511983	2.8958814561265253e-6	61594.0	0.0006720392505121201	1.954432388517943e-6	72538.0	0.0005074105817483548	1.9082695572485678e-6	54563.0	0.0005307483619589323	1.9234159294070673e-6	46028.0	0.0006732290051399136	1.2433880723382833e-6	55466.0
0.0005202255248765691	2.5372230272266475e-5	65060.5	0.0006575039376084391	7.952899207449164e-6	76583.8	0.0006160842899906809	4.782836793659719e-6	61940.3	0.0006909634513762073	1.616537201347622e-6	68343.0	0.0007857225900786693	2.033339658502141e-6	85558.3
\end{filecontents*}
\pgfplotstabletypeset[%
        begin table={\begin{tabular}[t]},
        every head row/.style={
        before row={%
          \hline
          \multicolumn{3}{|c|}{$\nu = 0.6$} & \multicolumn{3}{c|}{$\nu = 0.7$} & \multicolumn{3}{c|}{$\nu = 0.8$} & \multicolumn{3}{c|}{$\nu = 0.9$} & \multicolumn{3}{c|}{$\nu = 1.0$}\\
          \hline
        },
        after row/.add={}{\hline},
        },
        header=true,
       precision=1,
       columns/0/.style ={column name={pres}, column type={|l}},
       columns/1/.style ={column name={dres}, column type={c}},
       columns/2/.style ={column name={\# iters}, column type={c|}},
       columns/3/.style ={column name={pres}, column type={c}},
       columns/4/.style ={column name={dres.}, column type={c}},
       columns/5/.style ={column name={\# iters}, column type={c|}},
       columns/6/.style ={column name={pres}, column type={c}},
       columns/7/.style ={column name={dres}, column type={c}},
       columns/8/.style ={column name={\# iters}, column type={c|}},
       columns/9/.style ={column name={pres}, column type={c}},
       columns/10/.style ={column name={dres}, column type={c}},
       columns/11/.style ={column name={\# iters}, column type={c|}},
       columns/12/.style ={column name={pres}, column type={c}},
       columns/13/.style ={column name={dres}, column type={c}},
       columns/14/.style ={column name={\# iters}, column type={c|}},
        every row no 9/.style={after row=\hline},
        every row no 10/.style={after row=\hline},
        every row 10 column 0/.style={
                postproc cell content/.style={
                @cell content/.add={$\bf}{$}
                }
        },
        every row 10 column 10/.style={
                postproc cell content/.style={
                @cell content/.add={$\bf}{$}
                }
        },
        every row 10 column 8/.style={
                postproc cell content/.style={
                @cell content/.add={$\bf}{$}
                }
        },
    ]{data/quadratics_gaussian_20_100_double.csv}
    \end{adjustbox}
    
\end{table}

We now compare these results to \Cref{table:qp_gaussian_20_100_triple}, which is similar to \cref{table:qp_gaussian_20_100_double}, but uses the triple loop version of power ALM, i.e., \cref{alg:algorithm_nonconvex} in which the primal update in step \ref{alg:x_update} is computed using \cref{alg:algorithm_ippm} (iPPM).
The inexact proximal point updates are computed using the accelerated adaptive proximal-gradient algorithm presented in \cite[Algorithm 2]{malitsky_adaptive_2020}, which, although a heuristic, significantly outperforms the other methods we tried.
We tuned $\beta_1 = 1, \omega = 1.1, \lambda = 0.01$.
Although the triple-loop version of power ALM has a better worst-case complexity, in practice we find that it is usually outperformed by its double-loop variant, as confirmed in this table.
We do, nevertheless, observe that the triple-loop method is more stable, in the sense that the number of gradient calls fluctuates less between QP realizations.
This was also observed for $\nu = 1$ in \cite{li_rate-improved_2021}.
Moreover, remark that as $\nu$ decreases, the primal residual or constraint violation also decreases, thus qualitatively confirming our theoretical results that decreasing $\nu$ yields faster constraint satisfaction.

\begin{table}[h]
    \caption{
        Power ALM with iPPM inner solver on solving random QPs of size $n = 100, m = 20$.
    }
    \label{table:qp_gaussian_20_100_triple}
    \centering
    
    \begin{adjustbox}{width=\textwidth}
    \setlength\extrarowheight{3pt}
    \pgfplotstabletypeset[%
        begin table={\begin{tabular}[t]},
        every head row/.style={
            before row={%
              \hline
              \vphantom{$q = 1$}\\
              \hline
            },
            after row/.add={}{\hline},
        },
        header=true,
        col sep=&,
        row sep=\\,
        string type,
        columns/{trial}/.style ={column name={trial}, column type={|c}},
        every row no 9/.style={after row=\hline},
        every row no 10/.style={after row=\hline},
    ]{
        \\
        trial\\
        1\\
        2\\
        3\\
        4\\
        5\\
        6\\
        7\\
        8\\
        9\\
        10\\
        avg.\\
    }%
    \begin{filecontents*}{data/quadratics_gaussian_20_100_triple.csv}
1.3958774192266347e-5	0.0009736280699351402	121671.0	3.836143278093282e-5	0.0008401750327339168	123688.0	1.1632250904063154e-5	0.0009417318664200493	138930.0	9.274078891813048e-6	0.0009977722631088124	127821.0	0.00012552260130095966	0.0009441916336383434	112122.0
1.3190761263238576e-5	0.0009369032650932832	126670.0	0.00047460681625731134	0.0007034355383409062	117531.0	0.00022396261635247266	0.0008093938727860899	110826.0	0.0003633672013608974	0.0005219879500879576	93423.0	0.0008756721066386934	0.0005351422166052075	101270.0
8.730324824355867e-6	0.000630124637373252	130836.0	2.7551194227145002e-5	0.0005395777302971338	141272.0	0.00019775314651996817	0.000882929777744552	141028.0	6.115480718199294e-5	0.0009484754281173315	155530.0	3.7528062197536256e-5	0.0009299026669133546	178944.0
1.5981959277118116e-7	0.0009060225110022594	2.949152e6	5.67068001511623e-5	0.0007244788289106296	103536.0	0.00012726319671953746	0.00047736832296246264	99070.0	0.00017839195705439656	0.0007907252642707263	87422.0	0.0007088696088739387	0.0007139401561130832	120197.0
9.50455820788226e-8	0.0009697799191236752	1.711523e6	2.9235893638405883e-7	0.0009619693944006685	194278.0	0.00015257087632182753	0.000813778413479856	216599.0	0.00016900747589465245	0.000301541468492923	157251.0	0.0005404706428804446	0.0008691337470604539	190104.0
5.0323186390072324e-8	0.000983846299343814	948562.0	8.581848302200716e-7	0.000997579676740307	155197.0	4.056726121861204e-6	0.000987725396301013	125829.0	1.8079524395180942e-5	0.0009359087816023774	124411.0	0.0007596446111235797	0.0009582411543538937	129809.0
1.02260953567529e-6	0.0008836532816851187	1.33579e6	6.550653323348798e-6	0.0009364969226142044	245898.0	2.9337309590779852e-5	0.0009254940596017428	134179.0	0.00011208214738353161	0.0008523070437636118	164119.0	0.00026780737944541845	0.00087648927433408	164780.0
2.0962483052601933e-5	0.00045753665742527994	95543.0	0.00015598521408310745	0.0005850443673882581	123413.0	0.0008735169525099308	0.000254266039946683	322421.0	0.0005188415586135917	0.00071199781638525	121251.0	8.081613120899229e-5	0.000875437134938546	122026.0
8.464255831928547e-7	0.0009666551568187813	872048.0	2.457683934780857e-6	0.0009795813408137043	301625.0	2.7292011509762208e-5	0.000955552432366758	113861.0	3.9890875435540946e-5	0.000834574720018821	125534.0	0.00016064920483607086	0.0008929508952569282	158511.0
0.00026756638166732275	0.00042392453325325877	89998.0	6.33783002657714e-5	0.0006998037606688883	87063.0	1.7115757016131616e-5	0.0008956582447319258	114669.0	9.521973485691211e-5	0.0008422352337732421	155710.0	7.843646008967722e-5	0.0009097153040561194	107181.0
3.265829484798937e-5	0.0008132074331053862	838179.3	8.267486387901642e-5	0.0007968142592908615	159350.1	0.00016645008435663348	0.0007943898426341133	151741.2	0.00015653093610685096	0.0007737525969621053	131247.2	0.00036354168085953113	0.000850514418327001	138494.4
\end{filecontents*}
\pgfplotstabletypeset[%
        begin table={\begin{tabular}[t]},
        every head row/.style={
        before row={%
          \hline
          \multicolumn{3}{|c|}{$\nu = 0.6$} & \multicolumn{3}{c|}{$\nu = 0.7$} & \multicolumn{3}{c|}{$\nu = 0.8$} & \multicolumn{3}{c|}{$\nu = 0.9$} & \multicolumn{3}{c|}{$\nu = 1.0$}\\
          \hline
        },
        after row/.add={}{\hline},
        },
        header=true,
       precision=1,
       columns/0/.style ={column name={pres}, column type={|l}},
       columns/1/.style ={column name={dres}, column type={c}},
       columns/2/.style ={column name={\# grads}, column type={c|}},
       columns/3/.style ={column name={pres}, column type={c}},
       columns/4/.style ={column name={dres.}, column type={c}},
       columns/5/.style ={column name={\# grads}, column type={c|}},
       columns/6/.style ={column name={pres}, column type={c}},
       columns/7/.style ={column name={dres}, column type={c}},
       columns/8/.style ={column name={\# grads}, column type={c|}},
       columns/9/.style ={column name={pres}, column type={c}},
       columns/10/.style ={column name={dres}, column type={c}},
       columns/11/.style ={column name={\# grads}, column type={c|}},
       columns/12/.style ={column name={pres}, column type={c}},
       columns/13/.style ={column name={dres}, column type={c}},
       columns/14/.style ={column name={\# grads}, column type={c|}},
        every row no 9/.style={after row=\hline},
        every row no 10/.style={after row=\hline},
        every row 10 column 0/.style={
                postproc cell content/.style={
                @cell content/.add={$\bf}{$}
                }
        },
        every row 10 column 10/.style={
                postproc cell content/.style={
                @cell content/.add={$\bf}{$}
                }
        },
        every row 10 column 11/.style={
                postproc cell content/.style={
                @cell content/.add={$\bf}{$}
                }
        },
    ]{data/quadratics_gaussian_20_100_triple.csv}
    \end{adjustbox}
    
\end{table}

\section{Conclusion}

In this paper we introduced and analyzed an inexact augmented Lagrangian method for nonconvex problems with nonlinear constraints, involving a potentially sharper penalty function. 
Taking into account the composite structure of the corresponding augmented Lagrangian function, we study the joint complexity of the scheme using two different subproblem oracles. 
One of these oracles is a novel proximal point scheme that exploits the specific structure of the subproblems, of which the objectives are H\"older smooth and weakly convex. 
The proposed augmented Lagrangian method generalizes existing works with conventional penalty terms that attain the best known convergence rate for nonconvex minimization with first-order methods. 
Notably, we proved that unconventional penalty terms yield faster constraint satisfaction at the cost of a slower decrease of the cost, thereby reflecting the sharper penalties in the complexity analysis.
It is noteworthy that our analysis also improves upon existing works by considering a generic convex term \(g\) in the objective, and by not assuming boundedness of the iterates. 
Numerical experiments indicate that also in practice such penalty terms perform well.

\subsubsection*{Acknowledgments}
This work was supported by: Research Foundation Flanders (FWO) research projects G081222N, G033822N, G0A0920N; Research Council KU Leuven C1 project No. C14/24/103; and EuroHPC Project: 101118139 Inno4Scale.

\bibliography{references}
\bibliographystyle{tmlr}

\clearpage
\appendix
\section {Proofs}

\subsection{Proofs of section \ref{sec:ipalm}} \label{sec:ipalm-proofs}

    \subsubsection*{Proof of \cref{lem:bounded-multipliers}}
    \begin{proof}
    By consecutively using the multiplier update (step \ref{eq:y_update_dual}), the triangle inequality, and the dual step size update rule (step \ref{eq:dual-step-size}) we obtain that
    \begin{align*}
        \Vert y^{k+1} \Vert &= \bigg\lVert y^1 + \sum_{i=2}^k \sigma_i A(x^i) \Vert A(x^i) \Vert^{q-1} \bigg\rVert \\
        &\leq \Vert y^1 \Vert + \sum_{i=2}^k \sigma_i \Vert A(x^i) \Vert^q \\
        &\leq \Vert y^1 \Vert + \sum_{i=2}^k \sigma_1 \frac{\Vert A(x^1) \Vert (\log2)^2}{i \left(\log(i+1)\right)^2} \\
        &= \Vert y^1 \Vert + c \sigma_1 \Vert A(x^1) \Vert (\log 2)^2 =: y_{\max},
    \end{align*}
    where $c=\sum_{i=2}^\infty \frac{1}{i \log^2(i+1)} < \infty$.
\end{proof}

    \subsubsection*{Proof of \cref{lem:bounded-level-sets}}
    We first establish the following lemma.
\begin{lemma} \label{lem:level-set-inclusion}
    Let \( \{ x^k \}_{k \in \bN}, \{ y^k \}_{k \in \bN}, \{ \beta^k \}_{k \in \bN} \) be generated by \cref{alg:algorithm_nonconvex}.
    Then for any \(x \in \dom g\) and \(k \geq 1\), 
    \begin{equation*}
        \mathcal{L}_{\beta_{k+1}, y^{k+1}}(x) \subseteq \mathcal{L}_{\beta_k, y^k}(x)
    \end{equation*}
\end{lemma}
\begin{proof}
    From the definition of the augmented Lagrangian, the multiplier update and the penalty parameter update with \(\omega > 1\) we have for \(k \geq 1\)
    \begin{align*}
        L_{\beta_{k+1}, y^{k+1}}(x) &= \varphi(x) + \langle A(x), y^{k+1} \rangle + \frac{\beta_{k+1}}{2} \Vert A(x) \Vert^{\nu+1}\\
        &= \varphi(x) + \langle A(x), y^{k} \rangle + \sigma_{k+1} \Vert A(x) \Vert^{1+\nu} + \omega \frac{\beta_{k}}{2} \Vert A(x) \Vert^{\nu+1}\\
        &\geq \varphi(x) + \langle A(x), y^{k} \rangle + \frac{\beta_{k}}{2} \Vert A(x) \Vert^{\nu+1}\\
        &= L_{\beta_{k}, y^{k}}(x). 
    \end{align*}
    It follows immediately that
    \begin{equation*}
        \mathcal{L}_{\beta_{k+1}, y^{k+1}}(x) \subseteq \mathcal{L}_{\beta_k, y^k}(x).
    \end{equation*}
\end{proof}

We now provide a proof for \cref{lem:bounded-level-sets}.
\begin{proof}
    The nonemptyness claim follows from the observation that \(\bar x \in \mathcal{L}_{\beta, y}(\bar x)\) always holds by definition of the sublevel set.
    For the compactness claim, recall that the augmented Lagrangian can be expressed as
    \(
        L_\beta(x, y) = \varphi(x) + \langle A(x), y \rangle + \frac{\beta}{2} \Vert A(x) \Vert^{\nu+1}.
    \)
    By \cite[Corollary 3.27]{rockafellar_variational_1998}, it suffices to establish level-coercivity of \(L_\beta(x, y)\) with respect to the primal variable.
    A sufficient condition for this is that \(\varphi\) is level-coercive, as guaranteed by \cref{assumption:level-coercive}, and that \(x \mapsto \langle A(x), y \rangle + \frac{\beta}{2} \Vert A(x) \Vert^\nu\) is bounded below \cite[Exercise 3.29 (b)]{rockafellar_variational_1998}.
    Note that \(
        \langle A(x), y \rangle + \frac{\beta}{2} \Vert A(x) \Vert^{\nu+1} \geq -\Vert A(x) \Vert \Vert y \Vert + \frac{\beta}{2} \Vert A(x) \Vert^{\nu+1}.
    \)
    If \(\Vert A(x) \Vert^{\nu} \leq \frac{2}{\beta} \Vert y \Vert\), the claimed boundedness follows from the fact that
    \begin{equation*}
        \langle A(x), y \rangle + \frac{\beta}{2} \Vert A(x) \Vert^{\nu+1} \geq -\left( \frac{2}{\beta} \Vert y \Vert \right)^{\frac{1}{\nu}} \Vert y \Vert + \frac{\beta}{2} \Vert A(x) \Vert^{\nu+1} \geq -\left( \frac{2}{\beta} \Vert y \Vert \right)^{\frac{1}{\nu}} \Vert y \Vert.
    \end{equation*}
    Else \(\Vert A(x) \Vert^{\nu} > \frac{2}{\beta} \Vert y \Vert\), and hence \[
        \langle A(x), y \rangle + \frac{\beta}{2} \Vert A(x) \Vert^{\nu+1} \geq -\Vert A(x) \Vert \Vert y \Vert + \frac{\beta}{2} \frac{2}{\beta} \Vert A(x) \Vert \Vert y \Vert = 0, 
    \]
    which also implies the claimed boundedness.

    We now show \cref{lem:bounded-level-sets-1} by induction.
    First, observe that the assumption \(L_{\beta_k}(x^{k+1}, y^k) \leq L_{\beta_k}(x^k, y^k)\) implies for \(k \geq 1\) that
    \begin{equation} \label{prf:bounded-level-sets-primal-decrease}
        \mathcal{L}_{\beta_k, y^k}(x^{k+1}) \subseteq \mathcal{L}_{\beta_k, y^k}(x^k).
    \end{equation}
    The base case \(k = 1\) follows by consecutively applying \cref{lem:level-set-inclusion} and \eqref{prf:bounded-level-sets-primal-decrease}, i.e., \[
        \mathcal{L}_{\beta_2, y^2}(x^2) \subseteq \mathcal{L}_{\beta_1, y^1}(x^2) \subseteq \mathcal{L}_{\beta_1, y^1}(x^1)
    \]
    Suppose that \cref{lem:bounded-level-sets-1} holds for some \(k \geq 1\).
    Then by \cref{lem:level-set-inclusion} we have
    \begin{align*}
        \mathcal{L}_{\beta_{k+1}, y^{k+1}}(x^{k+1}) &\subseteq \mathcal{L}_{\beta_k, y^k}(x^{k+1}) \subseteq \mathcal{L}_{\beta_k, y^k}(x^k),
    \end{align*}
    where again we applied \cref{lem:level-set-inclusion} and \eqref{prf:bounded-level-sets-primal-decrease}.
    This establishes \cref{lem:bounded-level-sets-1}.

    Finally, \cref{lem:bounded-level-sets-2} follows by consecutive application of \cref{lem:bounded-level-sets-1} and by compactness of \(\mathcal{L}_{\beta_1, y^1}(x^1)\).
\end{proof}

    \subsubsection*{Proof of \cref{th:convergence-outer}}
    \begin{proof}
    We first prove the claimed rate of $\varepsilon_{A, k+1} := \frac{Q_A}{\beta_k^\frac{1}{\nu}}$, and then that of $\varepsilon_{\varphi, k+1} := \frac{Q_f}{\beta_k}$.
    From \cref{alg:algorithm_nonconvex} step \ref{alg:x_update} we have that
    \begin{equation}
        \dist(-\nabla_x \psi_{\beta_k}(x^{k+1}, y^k), \partial g(x^{k+1})) \leq \varepsilon_{k+1}, \quad \forall k \geq 0.
    \end{equation}
    By definition of $\psi_{\beta}$ we therefore have that $\forall k \geq 0$
    \begin{align*}
        \dist(-\nabla f(x^{k+1}) - J_{A}^\top(x^{k+1}) y^k - \beta_k J_{A}^\top(x^{k+1}) \nabla \phi(A(x^{k+1})), \partial g(x^{k+1})) \leq \varepsilon_{k+1}
    \end{align*}
    which yields, after application of the triangle inequality, 
    \begin{equation} \label{prf:outer-complexity-dist-partial-g}
        \dist(- \beta_k J_{A}^\top(x^{k+1}) \nabla \phi(A(x^{k+1})), \partial g(x^{k+1})) \leq \Vert \nabla f(x^{k+1}) \Vert + \Vert J_{A}^\top(x^{k+1}) y^k \Vert + \varepsilon_{k+1}.
    \end{equation}
    By \cite[Theorem 25.6]{rockafellar_convex_1970} it follows that for all \(x \in \bR^n\), \[
        \partial g(x) = \cl (\conv S(x)) + N_{\dom g}(x),
    \]
    where \(S(x)\) is the set of all limits of sequences of the form \(\nabla g(x_1), \nabla g(x_2), \dots\) such that \(g\) is differentiable at \(x_i\) and \(x_i \to x\), and \(\cl (\conv S(x))\) is the closure of the convex hull of \(S(x)\).
    By \cref{assumption:convex-g} the function \(g\) is \(G\)-Lipschitz continuous on \(\cS\), and therefore any \(v \in \cl (\conv S(x))\) satisfies \(\Vert v \Vert \leq G\) \cite{}.

    By again applying the triangle inequality to \eqref{prf:outer-complexity-dist-partial-g}, we obtain
    \begin{align*}
        \dist(- \beta_k J_{A}^\top(x^{k+1}) \nabla \phi(A(x^{k+1})), N_{\dom g}(x^{k+1})) \leq \Vert \nabla f(x^{k+1}) \Vert + G + \Vert J_{A}^\top(x^{k+1}) y^k \Vert + \varepsilon_{k+1}.
    \end{align*}

    Since $\nabla \phi(A(x^{k+1})) = \Vert A(x^{k+1}) \Vert^{\nu-1} A(x^{k+1})$, we can further lower bound the l.h.s. of the previous inequality. We have that $\dist( \alpha x, N_{\dom g}(x^{k+1})) = \min_{z \in N_{\dom g}(x^{k+1})}\Vert \alpha x - z\Vert = \alpha \min_{z \in N_{\dom g}(x^{k+1})}\Vert x - \tfrac{z}{\alpha} \Vert$ for $\alpha > 0$ and thus since $N_{\dom g}(x^{k+1})$ is a cone, using \cref{assump:regularity} we obtain
    \begin{align*}
        \Vert \nabla f(x^{k+1}) \Vert + G + \Vert J_{A}^\top(x^{k+1}) y^k \Vert + \varepsilon_{k+1} &\geq \beta_k \Vert A(x^{k+1}) \Vert^{\nu-1} \dist(- J_{A}^\top(x^{k+1}) A(x^{k+1}), N_{\dom g}(x^{k+1}))\\
        &\geq R \beta_k \Vert A(x^{k+1}) \Vert^{\nu}.
    \end{align*}
    Therefore, we have the following inequality
    \begin{equation} \label{eq:bound_feas}
        \Vert A(x^{k+1}) \Vert^{\nu} \leq \frac{\Vert \nabla f(x^{k+1}) \Vert + G + \Vert J_{A}^\top(x^{k+1}) y^k \Vert + \varepsilon_{k+1}}{R \beta_k }.
    \end{equation}
    By \cref{lem:bounded-multipliers}, we have that $\Vert y^k \Vert \leq y_{\max}$ for all $k \geq 0$.
    Thus, the constraint violation is upper bounded by
    \begin{equation} \label{eq:feas}
        \Vert A(x^{k+1}) \Vert \leq \left( \frac{{\nabla f}_{\max} + G + {J_A}_{\max} y_{\max} + \varepsilon_{k+1}}{R \beta_k } \right)^{1/\nu},
    \end{equation}
    where we have used the Cauchy-Schwarz inequality. The claim regarding $\varepsilon_{A, k+1}$ now follows from the fact that $\varepsilon_{k+1} \leq \varepsilon_k$.
    Now for $\varepsilon_{\varphi, k+1}$, we have by the triangle inequality that
    \begin{align*}
        &\dist(-\nabla_x \psi_{\beta_{k}}(x^{k+1}, y^{k+1}), \partial g(x^{k+1}))\\
        &\leq \dist(-\nabla_x \psi_{\beta_{k}}(x^{k+1}, y^k), \partial g(x^{k+1})) + \Vert \nabla_x \psi_{\beta_k}(x^{k+1}, y^{k+1})-\nabla_x \psi_{\beta_{k}}(x^{k+1}, y^{k}) \Vert
    \end{align*}
    The first term of the r.h.s. is upper bounded by $\varepsilon_{k+1}$ due to \cref{alg:algorithm_nonconvex} step \ref{alg:x_update}, i.e.,
    \begin{equation*}
        \dist(-\nabla_x \psi_{\beta_{k}}(x^{k+1}, y^{k}), \partial g(x^{k+1})) \leq \varepsilon_{k+1}
    \end{equation*}
    and for the second term we have that
    \begin{align*}
        \Vert \nabla_x \psi_{\beta_{k}}(x^{k+1}, y^{k+1})-\nabla_x \psi_{\beta_{k}}(x^{k+1}, y^{k}) \Vert = \Vert J_A^\top(x^{k+1}) (y^{k+1} - y^{k}) \Vert &\leq {J_A}_{\max} \Vert y^{k+1} - y^{k} \Vert\\
        &= {J_A}_{\max} \sigma_{k+1} \Vert A(x^{k+1}) \Vert^\nu.
    \end{align*}
    Therefore, in combination with \eqref{eq:bound_feas}, we obtain
    \begin{align*}
        \dist(-\nabla_x \psi_{\beta_{k}}(x^{k+1}, y^{k+1}), \partial g(x^{k+1})) \leq \varepsilon_{k+1} + {J_A}_{\max} \sigma_{k+1} \frac{{\nabla f}_{\max} + G + {J_A}_{\max} y_{\max} + \varepsilon_{k+1}}{R \beta_{k} }
    \end{align*}
    Thus, since $\sigma_k \leq \sigma_1$, $\varepsilon_{k+1} \leq \varepsilon_k$ and $\varepsilon_{k+1} = \lambda / \beta_k$, we conclude that $x^{k+1}$ is $(\varepsilon_{\varphi, k+1}, \varepsilon_{A,k+1})$-stationary with multiplier $y^{k+1} + \beta_k \nabla \phi(A(x^{k+1}))$.
\end{proof}

    \subsubsection*{Proof of \cref{th:holder-lagrangian}}
    
\begin{proof}
    Since the gradient of $\psi_\beta$ with respect to its first argument, evaluated at $(x, y)$, is given by
    \begin{equation}
        \nabla_x \psi_{\beta}(x, y) = \nabla f(x) + J_{A}^\top(x) y + \beta J_{A}^\top(x) \nabla \phi(A(x)),
    \end{equation}
    it follows for $x, x' \in \cS$ that
    \begin{align} \label{prf:holder-three-terms} \nonumber
        \Vert \nabla_x \psi_{\beta}(x', y) - \nabla_x \psi_{\beta}(x, y) \Vert &\leq \Vert \nabla f(x') - \nabla f(x) \Vert + \Vert J_{A}(x') - J_{A}(x) \Vert \cdot \Vert y \Vert\\
        &+ \beta \Vert J_{A}^\top (x') \nabla \phi(A(x')) - J_{A}^\top (x) \nabla \phi(A(x)) \Vert.
    \end{align}
    Using \cref{assumption:smoothness} and \cref{lem:bounded-multipliers}, the first two terms are bounded by
    \begin{align}
        &\Vert \nabla f(x') - \nabla f(x) \Vert \leq H_f \Vert x' - x \Vert^{\nu_f},\\
        &\Vert J_{A}(x') - J_{A}(x) \Vert \cdot \Vert y \Vert \leq H_A \Vert y \Vert \cdot \Vert x' - x \Vert^{\nu_A},
    \end{align}
    As for the third term in \eqref{prf:holder-three-terms}, we have that
    \begin{align} \label{eq:hoelder_third_term} \nonumber
        &\beta \Vert J_{A}^\top (x') \nabla \phi(A(x')) - J_{A}^\top (x) \nabla \phi(A(x)) \Vert
        \\
        &\leq \beta \Vert J_{A}^\top (x') \nabla \phi(A(x')) - J_{A}^\top (x) \nabla \phi(A(x')) \Vert + \beta \Vert J_{A}^\top (x) \nabla \phi(A(x')) - J_{A}^\top (x) \nabla \phi(A(x)) \Vert \nonumber 
        \\
        &\leq \beta \Vert A(x') \Vert^\nu \ \Vert J_A(x') - J_A(x) \Vert + \beta \Vert J_A(x) \Vert \Vert \nabla \phi(A(x')) - \nabla \phi(A(x)) \Vert \nonumber 
        \\
        &\leq \beta H_A \Vert A(x') \Vert^\nu \Vert x' - x \Vert^{\nu_A} + \beta {J_A}_{\max} \Vert \nabla \phi(A(x')) - \nabla \phi(A(x)) \Vert,
    \end{align}
    where we have used consecutively the triangle inequality, the fact that $\Vert \nabla \phi(A(x')) \Vert = \Vert A(x') \Vert^\nu$, \cref{assumption:smoothness} and the boundedness of $\Vert J_A \Vert$ on $\cS$. 
    In light of \cite[Theorem 6.3]{rodomanov_smoothness_2020}, the function $\phi$ is $(2^{1-\nu}, \nu)$-H\"older smooth.
    Hence, it follows that
    \begin{equation*}
        \Vert \nabla \phi(A(x')) - \nabla \phi(A(x)) \Vert \leq 2^{1-\nu}\Vert A(x') - A(x) \Vert^\nu \leq 2^{1-\nu}{J_A}_{\max}^\nu \Vert x' - x \Vert^\nu
    \end{equation*}
    where in the last inequality we also used the Lipschitz continuity of $A$ on $\cS$. Putting this back into \eqref{eq:hoelder_third_term}, we obtain
    \begin{equation}
        \beta \Vert J_{A}^\top (x') \nabla \phi(A(x')) - J_{A}^\top (x) \nabla \phi(A(x)) \Vert \leq \beta H_A A_{\max}^\nu \|x'-x\|^{\nu_A} + 2^{1-\nu}\beta{J_A}_{\max}^{1+\nu}\Vert x' - x \Vert^\nu
    \end{equation}
    Finally, the claim follows by summing the previous inequalities and using the fact that
    \begin{alignat*}{4}
        &\Vert x' - x \Vert^{\nu_f} &&= \Vert x' - x \Vert^{\nu_f - q} \cdot \Vert x' - x \Vert^q  & &\leq D^{\nu_f - q}\Vert x' - x \Vert^q & &\leq \max \left\{ 1, D^{1 - q} \right\}\Vert x' - x \Vert^q
        \\
        &\Vert x' - x \Vert^{\nu_A} &&= \Vert x' - x \Vert^{\nu_A - q} \cdot \Vert x' - x \Vert^q & &\leq D^{\nu_A - q}\Vert x' - x \Vert^q & &\leq \max \left\{ 1, D^{1 - q} \right\}\Vert x' - x \Vert^q
        \\
        &\Vert x' - x \Vert^{\nu} &&= \Vert x' - x \Vert^{\nu - q} \cdot \Vert x' - x \Vert^q & &\leq D^{\nu - q}\Vert x' - x \Vert^q & &\leq \max \left\{ 1, D^{1 - q} \right\}\Vert x' - x \Vert^q.
    \end{alignat*}
\end{proof}

    \subsubsection*{Proof of \cref{thm:ghadimi}}
    \begin{proof}
We first remark that UPFAG monotonically decreases the objective, and hence by \cref{lem:bounded-level-sets} the iterates remain in a compact set.
Let $\{\bar{x}^{\text{ag}}_i\}_{i \in \bN}$ and $\{x^{\text{ag}}_i\}_{i \in \bN}$ be the sequences of iterates generated by the UPFAG method described in \cite[Algorithm 2]{ghadimi_generalized_2019} when applied to $\psi_k$. Let $\{\gamma_i\}_{i \in \bN}$ be the sequence of stepsizes defined in \cite[Equation (3.6)]{ghadimi_generalized_2019}. Then in light of \cite[Equation (3.33)]{ghadimi_generalized_2019}, after discarding some constants, we need $O \left( H_{{\beta_k}}^{1/q} \left[ \frac{L_k(x^k) - L_k^\star}{\varepsilon_{k+1}^2} \right]^{\frac{1+q}{2q}} \right)$ inner iterations to obtain a point $\Vert \bar{x}^{\text{ag}}_{i}-x^{\text{ag}}_{i-1} \Vert / \gamma_i \leq \tfrac{\varepsilon_{k+1}}{2}$. Now, note that from \cite[Equation (3.8)]{ghadimi_generalized_2019} we have that 
\begin{equation*}
    \bar{x}^{\text{ag}}_i = \argmin_{u \in \bR^n} \langle u, \nabla \psi_k(x^{\text{ag}}_{i-1}) \rangle + g(u) + \tfrac{1}{2\gamma_i}\|u-x^{\text{ag}}_{i-1}\|^2
\end{equation*}
and from the optimality conditions for this minimization problem \cite[Theorem 6.12]{rockafellar_variational_1998}:
\begin{equation*}
    -\nabla \psi_k(x^{\text{ag}}_{i-1}) - \tfrac{1}{\gamma_i}(\bar{x}^{\text{ag}}_i - x^{\text{ag}}_{i-1})\in \partial g(\bar{x}^{\text{ag}}_i) 
\end{equation*}
By adding and subtracting $\nabla \psi_k(\bar{x}^{\text{ag}}_i)$ we further have:
\begin{equation*}
    -\nabla \psi_k(\bar{x}^{\text{ag}}_i) + \nabla \psi_k(\bar{x}^{\text{ag}}_i)-\nabla \psi_k(x^{\text{ag}}_{i-1}) - \tfrac{1}{\gamma_i}(\bar{x}^{\text{ag}}_i - x^{\text{ag}}_{i-1})\in \partial g(\bar{x}^{\text{ag}}_i) 
\end{equation*}
Therefore, $\dist(-\nabla \psi_k(\bar{x}^{\text{ag}}_i) + \nabla \psi_k(\bar{x}^{\text{ag}}_i)-\nabla \psi_k(x^{\text{ag}}_{i-1}) - \tfrac{1}{\gamma_i}(\bar{x}^{\text{ag}}_i - x^{\text{ag}}_{i-1}), \partial g(\bar{x}^{\text{ag}}_i)) = 0$ and from the triangle inequality
\begin{align*}
    \dist(-\nabla \psi_k(\bar{x}^{\text{ag}}_i), \partial g(\bar{x}^{\text{ag}}_i)) 
    &\leq \Vert \nabla \psi_k(\bar{x}^{\text{ag}}_i)-\nabla \psi_k(x^{\text{ag}}_{i-1})\Vert + \tfrac{1}{\gamma_i}\Vert \bar{x}^{\text{ag}}_i - x^{\text{ag}}_{i-1}\Vert
    \\
    & \leq H_{\beta_k}\Vert \bar{x}^{\text{ag}}_i - x^{\text{ag}}_{i-1}\Vert^q + \tfrac{\varepsilon_{k+1}}{2} \leq H_{\beta_k} \tfrac{\varepsilon_{k+1}^q \gamma_i^q}{2^q} + \tfrac{\varepsilon_{k+1}}{2},
\end{align*}
where in the third inequality we used the bound $\Vert \bar{x}^{\text{ag}}_{i}-x^{\text{ag}}_{i-1} \Vert / \gamma_i \leq \tfrac{\varepsilon_{k+1}}{2}$ and the H\"older continuity of $\nabla \psi_k$, \cref{th:holder-lagrangian}. Therefore, by choosing $\gamma_i \leq \tfrac{\varepsilon_{k+1}^{(1-q)/q}}{(2^{1-q}H_{\beta_k})^{1/q}}$ we obtain the claimed result.
\end{proof}

    \subsubsection*{Proof of \cref{th:complexity-analysis}}
    \begin{proof}
    Remark that by \cref{lem:bounded-level-sets-2} the iterates remain in a compact set.
    Let us start by defining the \textit{first} power ALM (outer) iteration $K_A$ for which
    \begin{equation}
        \varepsilon_A \geq \varepsilon_{A, K_A+1} = \left( \frac{{\nabla f}_{\max} + G + {J_A}_{\max} y_{\max} + \varepsilon_{1}}{R \beta_{K_A} } \right)^{\frac{1}{\nu}} := \frac{Q_A}{\beta_{K_A}^{1/\nu}}.
    \end{equation}
    or, equivalently, for which $\beta_{K_A} \geq \frac{Q_A}{\epsilon_A^\nu}$. 
    Since $K_A$ is the smallest iteration index for which this holds and $\beta_k$ is increasing, it follows that $\beta_{K_A - 1} < \frac{Q_A}{\varepsilon_A^\nu}$.
    It follows from
    \begin{align} \label{prf:total-complexity-feasibility-v2}
        \Vert A(x^{K_A+1}) \Vert \leq \left( \frac{{\nabla f}_{\max} + G + {J_A}_{\max} y_{\max} + \varepsilon_{1}}{R \beta_{K_A} } \right)^{\frac{1}{\nu}}
    \end{align}
    that also $\varepsilon_A \geq \Vert A(x^{K_A+1}) \Vert$ holds.
    From the update rule for $\beta_k$ we have that
    \begin{equation}
        \omega^{K_A - 2} < \frac{1}{\varepsilon_A^\nu} \left( \frac{Q_A}{\beta_1} \right) := \frac{Q_A'}{\varepsilon_A^\nu}.
    \end{equation}
    After taking the logarithm of both sides, we obtain that
    \begin{equation} \label{eq:num_feas}
        K_A = \bigg\lceil \log_\omega \left( \frac{Q_A'}{\varepsilon_A^\nu} \right) \bigg\rceil + 2.
    \end{equation}
    Now, the number of total UPFAG iterations needed to obtain a point $x^{K_A+1}$ satisfying \eqref{prf:total-complexity-feasibility-v2} is upper bounded by the sum of the calls to UPFAG for all outer (power ALM) iterations. Note that since $\beta_k$ is increasing geometrically, the H\"older smoothness modulus defined in \cref{th:holder-lagrangian} is determined by it for large enough $k$. Therefore, in light of \cref{thm:ghadimi} we require at most $T_A$ UPFAG iterations, where
    \begin{align*}
        T_A &= \sum_{k = 1}^{K_A} O \left( H_{{\beta_k}}^{1/q} \left[ \frac{L_k(x^k) - L_k^\star}{\varepsilon_{k+1}^2} \right]^{\frac{1+q}{2q}} \right) = \sum_{k = 1}^{K_A} O \left( \beta_k^{1/q} \left[ D^{q+1} \beta_k^3 \right]^{\frac{1+q}{2q}} \right)\\
        &= O \left( K_A \beta_{K_A}^{1/q} \left[ D^{q+1} \beta_{K_A}^3 \right]^{\frac{1+q}{2q}} \right) \leq O \left( K_A \frac{{(\omega Q_A)}^{1/q}}{\varepsilon_A^{\nu/q}} \left[ D^{q+1} \left(\omega\frac{Q_A}{\varepsilon_A^\nu}\right)^3 \right]^{\frac{1+q}{2q}} \right)\\
        &= O \left( K_A \frac{{Q_A}^{1/q}}{\varepsilon_A^{\nu/q}} D^{\frac{(q+1)^2}{2q}} \frac{{Q_A}^{3\frac{1+q}{2q}}}{\varepsilon_A^{3\nu\frac{1+q}{2q}}} \right) = O \left( \frac{K_A {Q_A}^\frac{5 + 3q}{2q} D^{\frac{(1+q)^2}{2q}}}{\varepsilon_A^{\frac{\nu}{q}+ \frac{3\nu(1+q)}{2}}} \right),
    \end{align*}
    where the first equality follows by \cref{th:holder-lagrangian} and by the fact that $\varepsilon_{k+1} = \lambda / \beta_k$. The second equality follows by unrolling the sum and the first inequality by the fact that $\beta_k = \omega \beta_{k-1} \leq \omega \frac{Q_A}{\varepsilon_A^\nu}$.
    Thus, we obtain
    \begin{align*}
        T_A &= O \left( \frac{\log_\omega \left( \frac{Q_A'}{\varepsilon_A^\nu} \right) {Q_A}^\frac{5 + 3q}{2q} D^{\frac{(1+q)^2}{2q}}}{\varepsilon_A^{\frac{\nu}{q}+ \frac{3\nu(1+q)}{2}}} \right) = \tilde{O} \left( \frac{Q_A^\frac{5 + 3q}{2q} D^{\frac{(1+q)^2}{2q}}}{\varepsilon_A^{\frac{\nu}{q}+ \frac{3\nu(1+q)}{2}}} \right).
    \end{align*}
    Regarding the suboptimality tolerance, we can in a similar way define the \textit{first} power ALM (outer) iteration $K_\varphi$ for which
    \begin{align} \label{prf:total-complexity-suboptimality-v2}
        \varepsilon_\varphi \geq \varepsilon_{\varphi, K_\varphi+1} = \frac{1}{\beta_{K_\varphi}} \left( \lambda + {J_A}_{\max} \sigma_1 \frac{{\nabla f}_{\max} + G + {J_A}_{\max} y_{\max} + \varepsilon_{1}}{R} \right) := \frac{Q_f}{\beta_{K_\varphi}}
    \end{align}
    or, equivalently, $\beta_{{K_\varphi}} \geq \frac{Q_f}{\varepsilon_\varphi}$.
    Since $K_\varphi$ is the first iteration for which this holds, it follows that $\beta_{{K_\varphi} - 1} < \frac{Q_f}{\varepsilon_\varphi}$.
    From the update rule for $\beta_k$ we have that 
    \begin{equation}
        \omega^{K_\varphi - 2} < \frac{1}{\varepsilon_\varphi} \left( \frac{Q_f}{\beta_1} \right) := \frac{Q_f'}{\varepsilon_\varphi}
    \end{equation}
    After taking the logarithm on both sides, we obtain
    \begin{equation}
        K_\varphi = \bigg\lceil \log_\omega \left( \frac{Q_f'}{\varepsilon_\varphi} \right) \bigg\rceil + 2.
    \end{equation}
    Hence, we require at most $T_\varphi$ UPFAG iterations to obtain a point $x^{K_\varphi+1}$ satisfying \eqref{prf:total-complexity-suboptimality-v2}, where
    \begin{align*}
        T_\varphi &= \sum_{k = 1}^{K_\varphi} O \left( H_{{\beta_k}}^{1/q} \left[ \frac{L_k(x^k) - L_k^\star}{\varepsilon_{k+1}^2} \right]^{\frac{1+q}{2q}} \right) = \sum_{k = 1}^{K_\varphi} O \left( \beta_k^{1/q} \left[ D^{q+1} \beta_k^3 \right]^{\frac{1+q}{2q}} \right)\\
        &= O \left( K_\varphi \beta_{K_\varphi}^{1/q} \left[ D^{q+1} \beta_{K_\varphi}^3 \right]^{\frac{1+q}{2q}} \right) \leq O \left( K_\varphi \frac{Q_f^{1/q}}{\varepsilon_\varphi^{1/q}} \left[ D^{q+1} \left(\omega \frac{Q_f}{\varepsilon_\varphi}\right)^3 \right]^{\frac{1+q}{2q}} \right)\\
        &= O \left( K_\varphi \frac{Q_f^{1/q}}{\varepsilon_\varphi^{1/q}} D^{\frac{(q+1)^2}{2q}} \frac{Q_f^\frac{3(q+1)}{2q}}{\varepsilon_\varphi^\frac{3(q+1)}{2q}} \right) = O \left( \frac{K_\varphi Q_f^\frac{5 + 3q}{2q} D^{\frac{(1+q)^2}{2q}}}{\varepsilon_\varphi^{\frac{5 + 3q}{2q}}} \right),
    \end{align*}
    and the sequence of bounds follows from similar arguments as for $T_A$. Thus, we obtain
    \begin{align*}
        T_\varphi &= O \left( \frac{\log_\omega \left( \frac{Q_f'}{\varepsilon_\varphi} \right) Q_f^\frac{5 + 3q}{2q} D^{\frac{(1+q)^2}{2q}}}{\varepsilon_\varphi^{\frac{5 + 3q}{2q}}} \right) = \tilde{O} \left( \frac{Q_f^\frac{5 + 3q}{2q} D^{\frac{(1+q)^2}{2q}}}{\varepsilon_\varphi^{\frac{5 + 3q}{2q}}} \right).
    \end{align*}
\end{proof}

\subsection{Proofs of section \ref{sec:ippm}} \label{sec:ippm-proofs}

    \subsection*{Proof of \cref{lem:al-weak-convexity}}
    \begin{proof}
    The function \(f\) is $L_f$-weakly convex, since it is $L_f$-Lipschitz smooth, i.e. for any $x,x' \in \cX$ the following inequality holds:
    \begin{equation} \label{eq:weak_conv_f}
        f(x) \geq f(x') + \langle \nabla f(x'),x-x' \rangle - \tfrac{L_f}{2}\Vert x-x'\Vert^2
    \end{equation}
    Moreover, following the proof of \cite[Lemma 4.2]{drusvyatskiy_efficiency_2019}, we have for any $x, x'\in \cX$, and $y \in \bR^m$ that
    \begin{align} \label{eq:weak_conv_linear}
        \nonumber
        \langle y, A(x') \rangle & = \langle y, A(x) \rangle + \langle y, A(x') - A(x) \rangle\\
        \nonumber
        &\geq \langle y, A(x) \rangle + \langle y, J_A(x) (x'- x) \rangle - \frac{L_A \Vert y \Vert}{2} \Vert x'- x \Vert^2\\
        &= \langle y, A(x) \rangle + \langle J_A^\top (x) y, (x'- x) \rangle - \frac{L_A \Vert y \Vert}{2} \Vert x'- x \Vert^2,
    \end{align}
    where the inequality follows from the Lipschitz-continuity of $J_A$, i.e., from the fact that $\Vert A(x') - A(x) - J_A(x) (x'- x) \Vert \leq \frac{L_A}{2} \Vert x'- x \Vert^2$.
    Likewise, for any $x, x' \in \cX$ we have that
    \begin{align} \label{eq:weak_conv_penalty}
        \nonumber
        \beta \phi(A(x)) &\geq \beta \phi(A(x')) + \langle \nabla \phi(A(x')),A(x)-A(x') \rangle
        \\
        \nonumber
        & \geq \beta \phi(A(x')) + \langle \nabla \phi(A(x')), J_A(x')(x-x') \rangle - \tfrac{L_A\Vert\nabla \phi(A(x'))\Vert}{2}\Vert x-x'\Vert^2
        \\
        &\geq \beta \phi(A(x')) + \beta \langle J_A^\top(x') \nabla \phi(A(x')), x- x' \rangle - \frac{\beta L_A A_{\max}^\nu}{2} \Vert x- x' \Vert^2
    \end{align}
    by consecutively exploiting convexity of $\phi$ between points $A(x)$ and $ A(x')$, Lipschitz-smoothness of $A$ and the fact that $\Vert \nabla \phi(A(x)) \Vert = \Vert A(x) \Vert^\nu \leq A_{\max}^\nu$.
    Summing up \eqref{eq:weak_conv_f}, \eqref{eq:weak_conv_linear} and \eqref{eq:weak_conv_penalty} we obtain:
    \begin{equation*}
        L_\beta(x,y) \geq L_\beta(x',y) + \langle \nabla_x L_\beta(x',y),x-x' \rangle - \tfrac{L_f + L_A\Vert y \Vert + \beta L_A A_{\max}^\nu}{2}\Vert x-x' \Vert^2
    \end{equation*}
    These observations directly imply that for any $y \in \bR^m$ the augmented Lagrangian is $\rho$-weakly convex on $\cX$ in its first argument, with 
    \begin{equation*}
        \rho := L_f + L_A (\Vert y \Vert + \beta A_{\max}^\nu).
    \end{equation*}
\end{proof}

    \subsection*{Proof of \cref{thm:inner_inner_complexity}}
    The proof requires the following series of lemmata. Note that $F$ is defined in \cref{alg:algorithm_ippm} step \ref{alg:def_F} and is a $(H_F, \nu)$-H\"older smooth and $\rho$-strongly convex function.

\begin{lemma} \label{lem:fgm-cost-decrease-v2}
    Let $x \in \cX$ and define the point $x^+ := \argmin_{u \in \cX} \langle \nabla F(x), u - x \rangle + \frac{1}{2\gamma} \Vert u - x \Vert^{2}$ with $\gamma = \frac{1}{ H_F} \varepsilon^{\frac{1-\nu}{1+\nu}}$, $\varepsilon > 0$.
    Then, 
    \begin{equation}
        \tfrac{1}{2 \gamma} \Vert x^+ - x \Vert^{2} \leq F(x) - F(x^+) + \tfrac{H_F}{2} \varepsilon.
    \end{equation}
\end{lemma}
\begin{proof}
    The inner problem in the definition of $x^+$ is $\tfrac{1}{\gamma}$-strongly convex and as such for any $u \in \cX$ we obtain
    \begin{equation*}
        \langle \nabla F(x), x^+ - x \rangle + \tfrac{1}{2\gamma} \Vert x^+ - x \Vert^{2} \leq \langle \nabla F(x), x - u \rangle + \tfrac{1}{2\gamma} \Vert x - u \Vert^{2} - \tfrac{1}{2\gamma}\Vert x^+-u\Vert^2.
    \end{equation*}
    By choosing $u = x \in \cX$ we get
    \begin{equation} \label{eq:inner_bound}
        0 \geq \langle \nabla F(x), x^+ - x \rangle + \tfrac{1}{\gamma} \Vert x^+ - x \Vert^{2}
    \end{equation}
    In light of Young's inequality we have for any $a, b \geq 0$ that $a b \leq \frac{a^{p_1}}{p_1} + \frac{b^{p_2}}{p_2}$ with $\tfrac{1}{p_1} + \tfrac{1}{p_2} = 1$. Choosing     
    $p_1 = \frac{2}{1+\nu}, p_2 = \frac{2}{1 - \nu}$, it follows for $a = \varepsilon^{-\frac{1 - \nu}{2}}$, $b = \varepsilon^{\frac{1-\nu}{2}}$ that
    \begin{equation}
        \Vert x^+ - x \Vert^{\nu+1} \leq \tfrac{1+\nu}{2}\varepsilon^{-\frac{1-\nu}{1+\nu}} \Vert x^+ - x \Vert^{2} + \tfrac{1-\nu}{2}\varepsilon.
    \end{equation}
    Note that for $\nu=1$ the inequality becomes an equality and thus the case $\nu=1$ is also covered. Therefore, from the H\"older smoothness inequality for $F$ between points $x, x^+ \in \cX$ we get:
    \begin{align} \nonumber \label{eq:hoelder_approx}
        F(x^+) &\leq 
        F(x) + \langle \nabla F(x), x^+ - x \rangle+ \tfrac{H_F}{\nu+1}\Vert x^+-x\Vert^{\nu+1}
        \\
        &\leq F(x) + \langle \nabla F(x), x^+ - x \rangle+ \tfrac{H_F}{2}\varepsilon^{-\frac{1-\nu}{1+\nu}}\Vert x^+-x\Vert^2 + \tfrac{1-\nu}{1+\nu}\tfrac{H_F}{2}\varepsilon
    \end{align}
    Substituting \eqref{eq:hoelder_approx} in \eqref{eq:inner_bound} we obtain:
    \begin{equation*}
        0 \geq F(x^+) - F(x) + \tfrac{1}{2\gamma} \Vert x^+ - x \Vert^{2} - \tfrac{1-\nu}{1+\nu}\tfrac{H_F}{2}\varepsilon 
    \end{equation*}
    which proves the claim.
\end{proof}

\begin{lemma} \label{cor:fgm-cost-decrease-v2}
    Suppose that a point $x \in \cX$ satisfies $F(x) - F(x^\star) \leq \varepsilon$ for some $\varepsilon > 0$.
    Let $\gamma = \frac{1}{H_F} \varepsilon^{\frac{1-\nu}{1+\nu}}$ and $x^+ = \argmin_{u \in \cX} \langle \nabla F(x), u - x \rangle + \frac{1}{2\gamma} \Vert u - x \Vert^{2}$.
    Then, $\tfrac{\Vert x^+ - x \Vert^2}{\gamma^2} \leq 2 H_F(1+H_F) \varepsilon^{\frac{2\nu}{1+\nu}}$.
\end{lemma}
\begin{proof}
    From \cref{lem:fgm-cost-decrease-v2} and $F(x) - F(x^+) \leq F(x) - F(x^\star) \leq \varepsilon$ it follows that
    \begin{align}
        \frac{\Vert x^+ - x \Vert^2}{\gamma^2} &\leq \frac{2}{\gamma} \left( \varepsilon + \frac{H_F}{2} \varepsilon \right)\\
        &= 2H_F \varepsilon^{-\frac{1-\nu}{1+\nu}}  \left( \varepsilon + \frac{H_F}{2}\varepsilon \right).
    \end{align}
    Therefore, we obtain
    \[
        \frac{\Vert x^+ - x \Vert^2}{\gamma^2} \leq 2 H_F(1+H_F) \varepsilon^{\frac{2\nu}{1+\nu}}.
    \]
\end{proof}

\begin{lemma} \label{th:ippm-complexity-connection}
    Suppose that a point $x \in \cX$ satisfies $F(x) - F(x^\star) \leq \varepsilon$ for some $\varepsilon > 0$.
    Let $\gamma = \frac{1}{H_F} \varepsilon^{\frac{1-\nu}{1+\nu}}$ and $x^+ = \argmin_{u \in \cX}\langle \nabla F(x), u - x \rangle + \frac{1}{2\gamma} \Vert u - x \Vert^{2}$.
    Then, 
    \[
        \dist(-\nabla F(x^+), N_{\cX}(x^+)) \leq \Big ( H_F^{1-\nu}(2 H_F(1+H_F))^{\frac{\nu}{2}} + (2 H_F(1+H_F))^{1/2} \Big ) \varepsilon^{\frac{\nu}{1+\nu}}.
    \]
\end{lemma}
\begin{proof}
    By the optimality conditions of the inner minimization in $x^+ = \argmin_{u \in \cX}\langle \nabla F(x), u - x \rangle + \frac{1}{2\gamma} \Vert u - x \Vert^{2}$, we have that
    \begin{equation*}
        - \nabla F(x) - \tfrac{1}{\gamma} (x^+ - x) \in N_{\cX}(x^+),
    \end{equation*}
    which implies that
    \begin{equation*}
        - \nabla F(x^+) + \left( \nabla F(x^+) - \nabla F(x) - \tfrac{1}{\gamma} (x^+ - x) \right) \in N_{\cX}(x^+).
    \end{equation*}
    Observe that from the triangle inequality and the H\"older smoothness of $F$, we have
    \begin{equation}
        \Vert \nabla F(x^+) - \nabla F(x) - \tfrac{1}{\gamma}(x^+ - x) \Vert \leq H_F \Vert x^+ - x \Vert^\nu + \tfrac{1}{\gamma} \Vert x^+ - x \Vert
    \end{equation}
    and by \cref{cor:fgm-cost-decrease-v2} we can further bound 
    \begin{align*}
        H_F \Vert x^+ - x \Vert^\nu + \tfrac{1}{\gamma} \Vert x^+ - x \Vert
        &\leq 
        \gamma^\nu H_F \tfrac{\Vert x^+ - x \Vert^\nu}{\gamma^\nu} + \tfrac{1}{\gamma} \Vert x^+ - x \Vert\\
        &= 
        \gamma^\nu H_F \left[ (2 H_F(1+H_F))^{\frac{\nu}{2}}\varepsilon^{\frac{\nu^2}{1+\nu}} \right] + (2 H_F(1+H_F))^{1/2}\varepsilon^{\frac{\nu}{1+\nu}}
        \\
        &= \left[ (\frac{1}{H_F})^\nu \varepsilon^{\frac{\nu(1-\nu)}{1+\nu}} \right] H_F \left[ (2 H_F(1+H_F))^{\frac{\nu}{2}}\varepsilon^{\frac{\nu^2}{1+\nu}} \right] + (2 H_F(1+H_F))^{1/2}\varepsilon^{\frac{\nu}{1+\nu}}
        \\
        &= \varepsilon^\frac{\nu}{\nu+1}\Big ( H_F^{1-\nu}(2 H_F(1+H_F))^{\frac{\nu}{2}} + (2 H_F(1+H_F))^{1/2}
        \Big)
    \end{align*}
    The result then follows by the fact that $\dist(-\nabla F(x^+), N_{\cX}(x^+)) \leq  \dist(\nabla F(x^+) - \nabla F(x) - \frac{1}{\gamma} (x^+ - x), N_{\cX}(x^+)) \leq \Vert \nabla F(x^+) - \nabla F(x) - \frac{1}{\gamma} (x^+ - x) \Vert$.
\end{proof}
    
    \begin{proof}
    Consider a point $x \in \cX$ such that
    \[
        F(x) - F(x^\star) \leq \bar \varepsilon := \left[\widebar{H}^{-1} \varepsilon \right]^\frac{1+\nu}{\nu}. 
    \]
    Then, it follows from \cref{th:ippm-complexity-connection} that for $\gamma = \tfrac{1}{H_F} \bar \varepsilon^\frac{1-\nu}{1+\nu} = \frac{1}{H_F} \left[ \widebar H^{-1} \right]^\frac{1-\nu}{\nu} \varepsilon^\frac{1-\nu}{\nu}$ the point $x^+$ satisfies
    \begin{align*}
        \dist(-\nabla F(x^+), N_{\cX}(x^+)) &\leq \bar H \bar \varepsilon^{\frac{\nu}{1+\nu}} = \varepsilon.
    \end{align*}
    Thus, in light of \cite[\S 6.2 ]{devolder_first-order_2014}, it takes at most $T$ iterations of the FGM, and a single proximal-gradient step with stepsize $\gamma$ to find such a point $x^+$, where
    \begin{align}
        T &= \widetilde{O} \left( \frac{H_F^{\frac{2}{1+3\nu}}}{\rho^{\frac{\nu+1}{3\nu+1}}} \bar \varepsilon^{-\frac{1-\nu}{1+3\nu}} \right)
        \\
        &= \widetilde{O} \left( \frac{H_F^{\frac{2}{1+3\nu}}}{\rho^{\frac{\nu+1}{3\nu+1}}} \left ( \widebar{H}^{-\frac{1+\nu}{\nu}} \varepsilon^\frac{1+\nu}{\nu} \right)^{-\frac{1-\nu}{1+3\nu}} \right)
        \\
        &= \widetilde{O} \left( \frac{H_F^{\frac{2}{1+3\nu}}}{\rho^{\frac{\nu+1}{3\nu+1}}} \widebar{H}^{\frac{1+\nu}{\nu}\frac{1-\nu}{1+3\nu}} \varepsilon^{-\frac{1+\nu}{\nu}\frac{1-\nu}{1+3\nu}} \right).
    \end{align}
\end{proof}

    \subsection*{Proof of \cref{th:ippm-complexity}}
    
\begin{proof}
    Let $F_k(x) := \psi(x) + \rho \Vert x - x^k \Vert^2$.
    Since $\psi$ is $\rho$-weakly convex and has $(H, \nu)$-H\"older continuous gradients, it follows that $F_k$ is $\rho$-strongly convex and has $(H_F, \nu)$-H\"older continuous gradients on $\cX$, with $H_F = H + \rho \max \left\{ 1, D^{1-\nu}\right\}$.

    We can now trace the steps of \cite[Theorem 1]{li_rate-improved_2021} to conclude that \cref{alg:algorithm_ippm} terminates after at most
    \[
        T = \lceil \frac{32 \rho}{\varepsilon^2} (\psi(x^1) - \psi^\star ) + 1\rceil
    \]
    iterations and that the output $x$ satisfies $\dist(- \nabla \psi(x), N_{\cX}(x)) \leq \varepsilon$.

    Denote henceforth $F_k^\star = \min_{x \in \cX} F_k(x)$.
    By \cref{alg:algorithm_ippm} step 3 and \cite[Theorem 10.1]{rockafellar_variational_1998} we have that
    \(
        \dist(0, \partial(F_k + \delta_{\cX})(x^{k+1})) = \dist(-\nabla F_k(x^{k+1}), N_{\cX}(x^{k+1})) \leq \frac{\varepsilon}{4},
    \)
    and by $\rho$-strong convexity of $F_k$ and convexity of $\delta_{\cX}$ we have that 
    \begin{align*}
        F_k^\star \geq F_k(x^{k+1}) + \langle v^{k+1}, x^\star - x^{k+1} \rangle + \frac{\rho}{2} \Vert x^{k+1} - x^\star \Vert^2, \qquad v^{k+1} \in \partial(F_k+\delta_{\cX})(x^{k+1}).
    \end{align*}
    Combining these two inequalities yields
    \(
        F_k(x^{k+1}) - F_k^\star \leq \frac{\varepsilon^2}{32\rho}
    \)
    and hence we obtain
    \(
        \psi(x^{k+1}) + \rho \Vert x^{k+1} - x^k \Vert^2 - \psi(x^k) \leq \frac{\varepsilon^2}{32 \rho}.
    \)
    Thus, 
    \begin{align*}
        \psi(x^T) - \psi(x^1) + \rho \sum_{k = 1}^{T-1} \Vert x^{k+1} - x^k \Vert^2 &\leq (T-1) \frac{\varepsilon^2}{32\rho}\\
        (T-1) \min_{1 \leq k \leq T-1} \Vert x^{k+1} - x^k \Vert^2 &\leq \frac{1}{\rho} \left( (T-1) \frac{\varepsilon^2}{32 \rho} + \left[ \psi(x^1) - \psi(x^T) \right]\right)\\
        2 \rho \min_{1 \leq k \leq T-1} \Vert x^{k+1} - x^k \Vert &\leq 2 \sqrt{\frac{\varepsilon^2}{32} + \frac{\rho \left[\psi(x^1) - \psi^\star \right]}{T-1}}
    \end{align*}
    For $T - 1 \geq \frac{32 \rho}{\varepsilon^2} (\psi(x^1) - \psi^\star)$ this yields
    \[
        2 \rho \min_{1 \leq k \leq T} \Vert x^{k+1} - x^k \Vert \leq \frac{\varepsilon}{2},
    \]
    meaning that \cref{alg:algorithm_ippm} must have terminated.
    Denote by $x^T$ the iterate returned by \cref{alg:algorithm_ippm} upon termination.
    Since $2 \rho \Vert x^T - x^{T-1} \Vert \leq \frac{\varepsilon}{2}$ and $\dist(0, \partial (F_{T-1} + \delta_{\cX})(x^T)) \leq \frac{\varepsilon}{4}$, we conclude that
    \begin{align*}
        \dist(0, \partial(\psi+\delta_{\cX}(x^T))) &\leq \dist(0, \partial(F_{T-1}+\delta_{\cX}(x^T))) + 2 \rho \Vert x^T - x^{T-1} \Vert \leq \varepsilon. 
    \end{align*}
    This concludes the proof, since $\partial(\psi+\delta_{\cX})(\cdot) = \nabla \psi(\cdot) + N_{\cX}(\cdot)$.
\end{proof}

    \subsection*{Proof of \cref{th:triple-loop-joint-complexity}}
    \begin{proof}
    We denote the $t$-th iterate of \cref{alg:algorithm_ippm} (iPPM) within the $k$-th power ALM outer iteration by $x_k^t$.
    At $x_k^t$, Nesterov's FGM is used to minimize $F_k^t := L_{\beta_k}(\cdot, y^k) + \rho_k \Vert \cdot - x_k^t \Vert^2$, which is $\rho_k$-strongly convex and $(H_{F_k}, \nu)$-H\"older smooth with $H_{F_k}= H_k + \rho_k \max\{1, D^{1-\nu}\}$. In light of \cref{thm:inner_inner_complexity}, we require at most $T_k^{FGM}$ iterations to find an $\frac{\varepsilon_{k+1}}{4}$ stationary point of $F_k^t + \delta_{\cX}$, with
    \begin{equation}
        T_k^{FGM} = \widetilde{O} \left( \frac{H_{F_k}^{\frac{2}{1+3\nu}}}{\rho_k^{\frac{1+\nu}{1+3\nu}}} \widebar{H}_k^{\frac{1+\nu}{\nu}\frac{1-\nu}{1+3\nu}} \varepsilon_{k+1}^{-\frac{1+\nu}{\nu}\frac{1-\nu}{1+3\nu}} \right),
    \end{equation}
    where $\widebar{H}_k = H_{F_k}^{1-\nu}(2 H_{F_k}(1+H_{F_k}))^{\frac{\nu}{2}} + (2 H_{F_k}(1+H_{F_k}))^{1/2}$.
    Remark that $T_k^{FGM}$ is independent of $t$.
    Since $\rho_k = O(\beta_k), H_{F_k} = O(\beta_k), \bar H_k = O(\beta_k)$, this expression simplifies to
    \[
        T_k^{FGM} = \widetilde{O} \left( \beta_k^{\frac{1-\nu}{1+3\nu}\left( 1+ \frac{1+\nu}{\nu} \right)} \varepsilon_{k+1}^{-\frac{1+\nu}{\nu}\frac{1-\nu}{1+3\nu}} \right).
    \]
    Likewise, at $x^k$, the iPPM is used to minimize $L_{\beta_k}(\cdot, y^k)$, which is $\rho_k$-weakly convex.
    By \cref{th:ippm-complexity}, we require at most $T_k^{PPM}$ iterations of \cref{alg:algorithm_ippm} to find a point satisfying the update rule of step \ref{alg:x_update} in \cref{alg:algorithm_nonconvex}, where
    \begin{align}
        T_k^{PPM} &= \lceil \frac{32 \rho_k}{\varepsilon_{k+1}^2} (\psi_k(x^k) - \psi_k^\star) + 1 \rceil
    \end{align}
    From \eqref{eq:bound_feas} we have that $\Vert A(x^{k}) \Vert^q \leq \left( \frac{{\nabla f}_{\max} + G + {J_A}_{\max} y_{\max} + \varepsilon_{k}}{R \beta_{k-1}} \right)$ for all $k \geq 1$ and using $\beta_k = \omega \beta_{k-1}$ we obtain
    \begin{equation} \label{eq:bound_beta_A}
        \beta_k \Vert A(x^k) \Vert^\nu \leq \omega \frac{{\nabla f}_{\max} + G + {J_A}_{\max} y_{\max} + \varepsilon_{k}}{R}
    \end{equation}
    and hence we obtain the following bound for $k \geq 2$:
    \begin{align*}
        \psi_k(x^k) &= f(x^k) + \langle y^k, A(x^k) \rangle + \frac{\beta_k}{1+\nu} \Vert A(x^k) \Vert^{\nu+1}
        \\
        &\leq f_{\max} + y_{\max} \Vert A(x^k) \Vert + \frac{\omega}{1+\nu} \left( \frac{{\nabla f}_{\max} + G + {J_A}_{\max} y_{\max} + \varepsilon_{k}}{R} \right) \Vert A(x^k) \Vert
        \\
        &\leq f_{\max} + y_{\max} \left( \frac{{\nabla f}_{\max} + G + {J_A}_{\max} y_{\max} + \varepsilon_{k}}{R \beta_{k-1}} \right)^\frac{1}{\nu}
        \\
        &\quad+ \frac{\omega}{1+\nu} \left( \frac{{\nabla f}_{\max} + G + {J_A}_{\max} y_{\max} + \varepsilon_{k}}{R} \right) \left( \frac{{\nabla f}_{\max} + G + {J_A}_{\max} y_{\max} + \varepsilon_{k}}{R \beta_{k-1}} \right)^{\frac{1}{\nu}}
        \\
        &\leq f_{\max} + y_{\max} + \left( \frac{{\nabla f}_{\max} + G + {J_A}_{\max} y_{\max} + \varepsilon_{1}}{R \beta_{1}} \right)^\frac{1}{\nu}
        \\
        &\quad+ \frac{\omega}{1+\nu} \left( \frac{{\nabla f}_{\max} + G + {J_A}_{\max} y_{\max} + \varepsilon_{1}}{R} \right) \left( \frac{{\nabla f}_{\max} + G + {J_A}_{\max} y_{\max} + \varepsilon_{1}}{R \beta_{1}} \right)^{\frac{1}{\nu}} =: C.
    \end{align*}
    The first inequality follows by $f$ being  continuous on the compact set $\cX$, $f(x^k) \leq \max_{x \in \cX}f(x) := f_{\max}$, and \eqref{eq:bound_beta_A} and the second one by \eqref{eq:feas}. The last inequality follows by $\varepsilon_k < \varepsilon_1$ and $\beta_k > \beta_1$.
    Moreover, we have for all $x \in \cX$
    \begin{equation*}
        \psi_k(x) \geq f(x) + \langle y^k, A(x) \rangle \geq -f_{\max} - y_{\max} A_{\max}
    \end{equation*}
    and as such we can bound
    \[
        \psi_k(x^k) - \psi^\star \leq C + f_{\max} + y_{\max} A_{\max}.
    \]
    Using this expression and \cref{lem:al-weak-convexity}, we can rewrite the number of iPPM iterations as
    \[
        T_k^{PPM} = \widetilde O \left( \frac{\beta_k}{\varepsilon_{k+1}^2}  \right)
    \]

    We now define the \textit{first} power ALM (outer) iteration $K_A$ for which
    \begin{equation}
        \varepsilon_A \geq \varepsilon_{A, K_A+1} = \left( \frac{{\nabla f}_{\max} + G + {J_A}_{\max} y_{\max} + \varepsilon_{1}}{R \beta_{K_A} } \right)^{\frac{1}{\nu}} := \frac{Q_A}{\beta_{K_A}^{1/\nu}}
    \end{equation}
    
    and from \eqref{eq:num_feas} we have that
    \begin{equation*}
        K_A = \bigg\lceil \log_\omega \left( \frac{Q_A'}{\varepsilon_A^\nu} \right) \bigg\rceil + 2.
    \end{equation*}
    Hence, we require at most $T_A$ total FGM iterations to obtain a point $x^{K_A+1}$, where
    \begin{align*}
        T_A &= \sum_{k = 1}^{K_A} T_k^{PPM} T_k^{FGM}\\
        &= \sum_{k = 1}^{K_A} \widetilde O \left( \frac{\beta_k}{\varepsilon_{k+1}^2} \beta_k^{\frac{1-\nu}{1+3\nu}\left( 1+ \frac{1+\nu}{\nu} \right)} \varepsilon_{k+1}^{-\frac{1+\nu}{\nu}\frac{1-\nu}{1+3\nu}} \right)
        = \sum_{k = 1}^{K_A} \widetilde O \left( \beta_k^{1 + \frac{1-\nu}{1+3\nu}\left( 1+ \frac{1+\nu}{\nu} \right)} \beta_{k}^{2 +\frac{1+\nu}{\nu}\frac{1-\nu}{1+3\nu}} \right)\\
        &= \sum_{k = 1}^{K_A} \widetilde O \left( \beta_k^{3 + \frac{1-\nu}{1+3\nu}\left( 1+ 2 \frac{1+\nu}{\nu} \right)} \right)\\
        &\leq \widetilde O \left( K_A \beta_{K_A}^{3 + \frac{1-\nu}{1+3\nu}\left( 1+ 2 \frac{1+\nu}{\nu} \right)} \right)
        = \widetilde O \left( \beta_{K_A}^{3 + \frac{1-\nu}{1+3\nu}\left( 1+ 2 \frac{1+\nu}{\nu} \right)} \right)
    \end{align*}
    By substitution of \(
        \beta_{K_A} = \omega \beta_{K_{A} - 1} < \frac{Q_A}{\varepsilon_A^\nu}
    \)
    this yields
    \begin{align*}
        T_A = \widetilde O \left( \varepsilon_{A}^{-3\nu - \frac{1-\nu}{1+3\nu}\left( 3\nu + 2 \right)} \right).
    \end{align*}
    In a similar way to the proof of \cref{th:complexity-analysis} we can define the \textit{first} power ALM (outer) iteration $K_\varphi$ for which
    \begin{equation*} \label{prf:total-complexity-suboptimality}
        \varepsilon_\varphi \geq \varepsilon_{\varphi, K_\varphi+1} = \frac{1}{\beta_{K_\varphi}} \left( \lambda + {J_A}_{\max} \sigma_1 \frac{{\nabla f}_{\max} + G + {J_A}_{\max} y_{\max} + \varepsilon_{1}}{R} \right) := \frac{Q_f}{\beta_{K_\varphi}}
    \end{equation*}
    and obtain
    \begin{equation*}
        K_\varphi = \bigg\lceil \log_\omega \left( \frac{Q_f'}{\varepsilon_\varphi} \right) \bigg\rceil + 2.
    \end{equation*}
    \begin{align*}
        T_\varphi &= \widetilde O \left( \beta_{K_\varphi}^{3 + \frac{1-\nu}{1+3\nu}\left( 1+ 2 \frac{1+\nu}{\nu} \right)} \right)
    \end{align*}
    By substitution of $\beta_{K_\varphi} = \omega \beta_{K_\varphi - 1} < \frac{Q_f}{\varepsilon_\varphi}$ this yields
    \begin{align*}
        T_\varphi = \widetilde O \left( \varepsilon_\varphi^{-3 - \frac{1-\nu}{1+3\nu}\left( 1 + \frac{2 (1+\nu)}{\nu} \right)} \right).
    \end{align*}
\end{proof}

    \subsection*{Proof of \cref{th:triple-loop-joint-complexity-linear}}
    \begin{proof}
    We denote the $t$-th iterate of \cref{alg:algorithm_ippm} (iPPM) within the $k$-th power ALM outer iteration by $x_k^t$.
    At $x_k^t$, Nesterov's FGM is used to minimize $F_k^t := L_{\beta_k}(\cdot, y^k) + \rho_k \Vert \cdot - x_k^t \Vert^2$, which is $\rho_k$-strongly convex and $(H_{F_k}, \nu)$-H\"older smooth with $H_{F_k}= H_k + \rho_k \max\{1, D^{1-\nu}\}$. In light of \cref{thm:inner_inner_complexity}, we require at most $T_k^{FGM}$ iterations to find an $\frac{\varepsilon_{k+1}}{4}$ stationary point of $F_k^t + \delta_{\cX}$, with
    \begin{equation}
        T_k^{FGM} = \widetilde{O} \left( \frac{H_{F_k}^{\frac{2}{1+3\nu}}}{\rho_k^{\frac{1+\nu}{1+3\nu}}} \widebar{H}_k^{\frac{1+\nu}{\nu}\frac{1-\nu}{1+3\nu}} \varepsilon_{k+1}^{-\frac{1+\nu}{\nu}\frac{1-\nu}{1+3\nu}} \right),
    \end{equation}
    where $\widebar{H}_k = H_{F_k}^{1-\nu}(2 H_{F_k}(1+H_{F_k}))^{\frac{\nu}{2}} + (2 H_{F_k}(1+H_{F_k}))^{1/2}$.
    Remark that $T_k^{FGM}$ is independent of $t$.
    Since $\rho_k = O(1), H_{F_k} = O(\beta_k), \bar H_k = O(\beta_k)$, this expression simplifies to
    \[
        T_k^{FGM} = \widetilde{O} \left( \beta_k^{\frac{1-\nu}{1+3\nu}\frac{1+\nu}{\nu} + \frac{2}{1+3\nu}} \varepsilon_{k+1}^{-\frac{1+\nu}{\nu}\frac{1-\nu}{1+3\nu}} \right).
    \]
    Likewise, at $x^k$, the iPPM is used to minimize $L_{\beta_k}(\cdot, y^k)$, which is $\rho_k$-weakly convex.
    By \cref{th:ippm-complexity}, we require at most $T_k^{PPM}$ iterations of \cref{alg:algorithm_ippm} to find a point satisfying the update rule of step \ref{alg:x_update} in \cref{alg:algorithm_nonconvex}, where
    \begin{align}
        T_k^{PPM} &= \lceil \frac{32 \rho_k}{\varepsilon_{k+1}^2} (\psi_k(x^k) - \psi_k^\star) + 1 \rceil
    \end{align}
    Following the same steps as in the proof of \cref{th:triple-loop-joint-complexity}, we can bound for some $C > 0$,
    \[
        \psi_k(x^k) - \psi^\star \leq C + f_{\max} + y_{\max} A_{\max}.
    \]
    Using this expression and $\rho_k = O(1)$ (cf.\,\cref{lem:al-weak-convexity}), we can rewrite the number of iPPM iterations as
    \[
        T_k^{PPM} = \widetilde O \left( \varepsilon_{k+1}^{-2}  \right)
    \]

    We now define the \textit{first} power ALM (outer) iteration $K_A$ for which
    \begin{equation}
        \varepsilon_A \geq \varepsilon_{A, K_A+1} = \left( \frac{{\nabla f}_{\max} + {J_A}_{\max} y_{\max} + \varepsilon_{1}}{v \beta_{K_A} } \right)^{\frac{1}{\nu}} := \frac{Q_A}{\beta_{K_A}^{1/\nu}}
    \end{equation}
    
    and from \eqref{eq:num_feas} we have that
    \begin{equation*}
        K_A = \bigg\lceil \log_\omega \left( \frac{Q_A'}{\varepsilon_A^\nu} \right) \bigg\rceil + 2.
    \end{equation*}
    Hence, we require at most $T_A$ total FGM iterations to obtain a point $x^{K_A+1}$, where
    \begin{align*}
        T_A &= \sum_{k = 1}^{K_A} T_k^{PPM} T_k^{FGM}\\
        &= \sum_{k = 1}^{K_A} \widetilde O \left( \varepsilon_{k+1}^{-2} \beta_k^{\frac{1-\nu}{1+3\nu}\frac{1+\nu}{\nu} + \frac{2}{1+3\nu}} \varepsilon_{k+1}^{-\frac{1+\nu}{\nu}\frac{1-\nu}{1+3\nu}} \right)
        = \sum_{k = 1}^{K_A} \widetilde O \left( \beta_k^{\frac{1-\nu}{1+3\nu}\frac{1+\nu}{\nu} + \frac{2}{1+3\nu}} \beta_{k}^{2 +\frac{1+\nu}{\nu}\frac{1-\nu}{1+3\nu}} \right)\\
        &= \sum_{k = 1}^{K_A} \widetilde O \left( \beta_k^{2 + 2 \frac{1-\nu}{1+3\nu}\frac{1+\nu}{\nu} + \frac{2}{1+3\nu}} \right)\\
        &\leq \widetilde O \left( K_A \beta_{K_A}^{2 + 2 \frac{1-\nu}{1+3\nu}\frac{1+\nu}{\nu} + \frac{2}{1+3\nu}} \right)
        = \widetilde O \left( \beta_{K_A}^{2 + 2 \frac{1-\nu}{1+3\nu}\frac{1+\nu}{\nu} + \frac{2}{1+3\nu}} \right)
    \end{align*}
    By substitution of \(
        \beta_{K_A} = \omega \beta_{K_{A} - 1} < \frac{Q_A}{\varepsilon_A^\nu}
    \)
    this yields
    \begin{align*}
        T_A = \widetilde O \left( \varepsilon_{A}^{-2\nu - 2 \frac{1-\nu}{1+3\nu}(1+\nu) - \frac{2\nu}{1+3\nu}} \right).
    \end{align*}
    In a similar way to the proof of \cref{th:complexity-analysis} we can define the \textit{first} power ALM (outer) iteration $K_\varphi$ for which
    \begin{equation*}
        \varepsilon_\varphi \geq \varepsilon_{\varphi, K_\varphi+1} = \frac{1}{\beta_{K_\varphi}} \left( \lambda + {J_A}_{\max} \sigma_1 \frac{{\nabla f}_{\max} + {J_A}_{\max} y_{\max} + \varepsilon_{1}}{v} \right) := \frac{Q_f}{\beta_{K_\varphi}}
    \end{equation*}
    and obtain
    \begin{equation*}
        K_\varphi = \bigg\lceil \log_\omega \left( \frac{Q_f'}{\varepsilon_\varphi} \right) \bigg\rceil + 2.
    \end{equation*}
    \begin{align*}
        T_\varphi &= \widetilde O \left( \beta_{K_\varphi}^{2 + 2 \frac{1-\nu}{1+3\nu}\frac{1+\nu}{\nu} + \frac{2}{1+3\nu}} \right)
    \end{align*}
    By substitution of $\beta_{K_\varphi} = \omega \beta_{K_\varphi - 1} < \frac{Q_f}{\varepsilon_\varphi}$ this yields
    \begin{align*}
        T_\varphi = \widetilde O \left( \varepsilon_\varphi^{-2 - 2 \frac{1-\nu}{1+3\nu}\frac{1+\nu}{\nu} - \frac{2}{1+3\nu}} \right).
    \end{align*}
\end{proof}
\clearpage
\section{Additional experiments} \label{sec:numerics-appendix}

\subsection{Generalized eigenvalue problem}
We consider the generalized eigenvalue problem (GEVP)
\begin{equation*}
    \min_{x \in \bR^n} x^\top C x \quad \text{s.t.} \quad x^\top B x = 1,
\end{equation*}
where $B, C \in \bR^{n \times n}$ are symmetric matrices and $B$ is positive definite.
Clearly, the GEVP is of the form \eqref{eq:nonconvex_problem}, and satisfies the regularity condition \cref{assump:regularity} \cite{sahin_inexact_2019}.
We sample the entries of a matrix $\widehat C \in \bR^{n \times n}$ from a Gaussian $\cN(0, 0.1)$ and define $C := \frac{1}{2}(\widehat C + \widehat C^\top)$.
The matrix $B$ is defined as $Q^\top Q$, where $Q$ is the orthonormal matrix in the QR-decomposition of a random matrix with entries sampled uniformly from the unit interval.
We use the classical accelerated proximal gradient method (APGM) by \cite{beck_fast_2009} as an inner solver and tune its step size to $0.5 / (10 \Vert C \Vert + 5000 + 500 \beta)$.
Although APGM has no convergence guarantees that fully cover our setting, i.e., nonconvex and Hölder-smooth objectives, it appears that convergence issues can be mitigated by sufficiently decreasing the step size.
We follow the tuning of \cite{li_rate-improved_2021} $\beta_1 = 0.01, \omega = 3, \lambda = 1, \sigma_1 = 10$, and impose a maximum of $N = 10^5$ APGM iterations per subproblem.

\Cref{table:gevp_gaussian} reports the number of gradient calls that power ALM requires to attain an $(\varepsilon_\varphi, \varepsilon_A)$-stationary point with $\varepsilon_\varphi = \varepsilon_A = 10^{-3}$ for various powers $\nu \in (0, 1]$.
Also the constraint violation $\Vert A(x) \Vert$ and the suboptimality $\vert f(x) - f^* \vert$ are listed.
Every trial denotes a random realization of the GEVP with $n = 500$, and the first trial is further illustrated in \cref{fig:gevp-gaussian}.
We observe that smaller values of $\nu$ perform significantly better than larger values, with $\nu = 0.4$ requiring an order of magnitude fewer gradient evaluations than $\nu = 1$.
The figure corresponding to the first realization confirms that both constraint violation and suboptimality decrease steadily, even for a small power $\nu = 0.4$.

\begin{figure}[h]
    \centering
    \begin{subfigure}[b]{0.32\textwidth}
        \centering
        \resizebox{\textwidth}{!}{
            \begin{tikzpicture}
\begin{axis}[
            width=3in, height=1.8in,
            at={(1.011in,0.642in)},
            legend cell align={left},
            legend pos=north east,
            scale only axis,
            minor grid style={thin,draw opacity=0.3},
            major grid style={thin,draw opacity=0.5},
            grid=both,
            xmin=0,
            xmax=600000,
            ymin={1e-5},
            ymax={3e2},
            ymode=log,
            xlabel = {\# Gradients},
            ylabel = {$\Vert A(x) \Vert$}
]
    \addplot[color={rgb,1:red,0.0;green,0.6056;blue,0.9787}, name path={07b8a02d-df97-432a-b42d-a902b9f993c3}, draw opacity={1.0}, line width={1}, solid, mark={*}, mark size={3.75 pt}, mark repeat={1}, mark options={color={rgb,1:red,0.0;green,0.0;blue,0.0}, draw opacity={1.0}, fill={rgb,1:red,0.0;green,0.6056;blue,0.9787}, fill opacity={1.0}, line width={0.75}, rotate={0}, solid}]
        table[row sep={\\}]
        {
            \\
            0.0  157.7440757562481  \\
            1.0  157.64395015642333  \\
            8.0  0.9997564517595128  \\
            9.0  0.9998512144998982  \\
            10.0  0.9999079843630169  \\
            13.0  0.99997951682994  \\
            17.0  0.9999970261594365  \\
            23.0  0.999999768677971  \\
            30.0  0.9999999548390669  \\
            43.0  0.9999999930196976  \\
            68.0  0.9999999989829155  \\
            111.0  0.9999999995992577  \\
            92346.0  0.6063381644975421  \\
            105484.0  0.20179633205298497  \\
            133085.0  0.06723161609087913  \\
            185408.0  0.022406909333490765  \\
            277544.0  0.0074685793839169845  \\
            377544.0  0.0024894843815402856  \\
            477544.0  0.0008298235834601764  \\
            577544.0  0.0002766073697592075  \\
        }
        ;
    \addlegendentry {$\nu = 1$}
    \addplot[color={rgb,1:red,0.8889;green,0.4356;blue,0.2781}, name path={2be47bf0-5dcf-46be-adda-37a825e41b5a}, draw opacity={1.0}, line width={1}, solid, mark={*}, mark size={3.75 pt}, mark repeat={1}, mark options={color={rgb,1:red,0.0;green,0.0;blue,0.0}, draw opacity={1.0}, fill={rgb,1:red,0.8889;green,0.4356;blue,0.2781}, fill opacity={1.0}, line width={0.75}, rotate={0}, solid}]
        table[row sep={\\}]
        {
            \\
            0.0  157.7440757562481  \\
            1.0  157.7072352905373  \\
            13.0  0.9741288583344594  \\
            14.0  0.9780404309125846  \\
            22.0  0.998065525943819  \\
            30.0  0.9997967941418394  \\
            38.0  0.9999692327933106  \\
            48.0  0.9999954974283504  \\
            64.0  0.9999995514495431  \\
            88.0  0.9999999162380093  \\
            132.0  0.9999999698800697  \\
            31965.0  0.5533297554680371  \\
            45653.0  0.13929026713134618  \\
            68713.0  0.03521021156688142  \\
            108398.0  0.008912386078438805  \\
            177110.0  0.002256859376550957  \\
            277110.0  0.000571576750591829  \\
        }
        ;
    \addlegendentry {$\nu = 0.8$}
    \addplot[color={rgb,1:red,0.2422;green,0.6433;blue,0.3044}, name path={6899e11a-b6aa-4f06-a0ae-fe4843b9ccc4}, draw opacity={1.0}, line width={1}, solid, mark={*}, mark size={3.75 pt}, mark repeat={1}, mark options={color={rgb,1:red,0.0;green,0.0;blue,0.0}, draw opacity={1.0}, fill={rgb,1:red,0.2422;green,0.6433;blue,0.3044}, fill opacity={1.0}, line width={0.75}, rotate={0}, solid}]
        table[row sep={\\}]
        {
            \\
            0.0  157.7440757562481  \\
            1.0  157.7302375304667  \\
            19.0  0.7857097354493026  \\
            33.0  0.8494056478981256  \\
            46.0  0.979560552844509  \\
            60.0  0.9975053096257274  \\
            76.0  0.9997369999837493  \\
            95.0  0.9999655211696473  \\
            122.0  0.9999941073015883  \\
            166.0  0.9999978349015831  \\
            10431.0  0.4623467143254175  \\
            14428.0  0.07221549464744392  \\
            24037.0  0.011478352731909558  \\
            43599.0  0.00183460777927702  \\
            79695.0  0.00029374943174809154  \\
        }
        ;
    \addlegendentry {$\nu = 0.6$}
    \addplot[color={rgb,1:red,0.7644;green,0.4441;blue,0.8243}, name path={fa4df1d8-7db6-4e21-957f-a3046716f9ae}, draw opacity={1.0}, line width={1}, solid, mark={*}, mark size={3.75 pt}, mark repeat={1}, mark options={color={rgb,1:red,0.0;green,0.0;blue,0.0}, draw opacity={1.0}, fill={rgb,1:red,0.7644;green,0.4441;blue,0.8243}, fill opacity={1.0}, line width={0.75}, rotate={0}, solid}]
        table[row sep={\\}]
        {
            \\
            0.0  157.7440757562481  \\
            1.0  157.73859742830945  \\
            26.0  11.830008326351193  \\
            49.0  0.2479651964530447  \\
            69.0  0.8883631183576465  \\
            92.0  0.9878904464177041  \\
            116.0  0.9980720384882631  \\
            148.0  0.9997158061813253  \\
            194.0  0.9999132854353228  \\
            2882.0  0.4049622677684871  \\
            9268.0  0.022893691165846608  \\
            17647.0  0.001412913894379919  \\
            30808.0  8.953606965866889e-5  \\
        }
        ;
    \addlegendentry {$\nu = 0.4$}
\end{axis}
\end{tikzpicture}
        }
        \captionsetup{justification=centering}
        \caption{Constraint violation.\label{fig:gevp-gaussian-feasibility}}
    \end{subfigure}
    \begin{subfigure}[b]{0.32\textwidth}
        \centering
        \resizebox{\textwidth}{!}{
            \begin{tikzpicture}
\begin{axis}[
            width=3in, height=1.8in,
            at={(1.011in,0.642in)},
            legend cell align={left},
            legend pos=north east,
            scale only axis,
            minor grid style={thin,draw opacity=0.3},
            major grid style={thin,draw opacity=0.5},
            grid=both,
            xmin=0,
            xmax=600000,
            ymin={1e-4},
            ymax={1e1},
            ymode=log,
            xlabel = {\# Gradients},
            ylabel = {$\vert f(x) - f^* \vert$}
]
    \addplot[color={rgb,1:red,0.0;green,0.6056;blue,0.9787}, name path={b5bb362e-038d-4ddc-98ff-f69d5518891c}, draw opacity={1.0}, line width={1}, solid, mark={*}, mark size={3.75 pt}, mark repeat={1}, mark options={color={rgb,1:red,0.0;green,0.0;blue,0.0}, draw opacity={1.0}, fill={rgb,1:red,0.0;green,0.6056;blue,0.9787}, fill opacity={1.0}, line width={0.75}, rotate={0}, solid}]
        table[row sep={\\}]
        {
            \\
            0.0  3.3491248847254083  \\
            1.0  3.310461713254738  \\
            8.0  1.5627855646372701  \\
            9.0  1.5627995710996698  \\
            10.0  1.562807978453057  \\
            13.0  1.5628186072285166  \\
            17.0  1.5628212258009098  \\
            23.0  1.5628216393679528  \\
            30.0  1.562821667726611  \\
            43.0  1.5628216735812643  \\
            68.0  1.562821674502787  \\
            111.0  1.5628216745986327  \\
            92346.0  0.9476203570710628  \\
            105484.0  0.31537824787124724  \\
            133085.0  0.1050719365012911  \\
            185408.0  0.035018098849379786  \\
            277544.0  0.011672070396886891  \\
            377544.0  0.003890626950040721  \\
            477544.0  0.0012968718864834194  \\
            577544.0  0.0004322932538018964  \\
        }
        ;
    \addlegendentry {$\nu = 1$}
    \addplot[color={rgb,1:red,0.8889;green,0.4356;blue,0.2781}, name path={4bb9a71e-46d2-4872-ad25-0a0e4098c81d}, draw opacity={1.0}, line width={1}, solid, mark={*}, mark size={3.75 pt}, mark repeat={1}, mark options={color={rgb,1:red,0.0;green,0.0;blue,0.0}, draw opacity={1.0}, fill={rgb,1:red,0.8889;green,0.4356;blue,0.2781}, fill opacity={1.0}, line width={0.75}, rotate={0}, solid}]
        table[row sep={\\}]
        {
            \\
            0.0  3.3491248847254083  \\
            1.0  3.3111663572510865  \\
            13.0  1.561464733491189  \\
            14.0  1.5616643099456196  \\
            22.0  1.5627112882638816  \\
            30.0  1.56280926036818  \\
            38.0  1.5628196928531177  \\
            48.0  1.5628213689341104  \\
            64.0  1.5628216421178334  \\
            88.0  1.562821668274782  \\
            132.0  1.5628216722515473  \\
            31965.0  0.8653556050575101  \\
            45653.0  0.21777837977125358  \\
            68713.0  0.055038153399437295  \\
            108398.0  0.013929693917241526  \\
            177110.0  0.0035272041233165297  \\
            277110.0  0.0008933043686405462  \\
        }
        ;
    \addlegendentry {$\nu = 0.8$}
    \addplot[color={rgb,1:red,0.2422;green,0.6433;blue,0.3044}, name path={37ee4a2d-89d6-4983-a603-0178eef6c9a4}, draw opacity={1.0}, line width={1}, solid, mark={*}, mark size={3.75 pt}, mark repeat={1}, mark options={color={rgb,1:red,0.0;green,0.0;blue,0.0}, draw opacity={1.0}, fill={rgb,1:red,0.2422;green,0.6433;blue,0.3044}, fill opacity={1.0}, line width={0.75}, rotate={0}, solid}]
        table[row sep={\\}]
        {
            \\
            0.0  3.3491248847254083  \\
            1.0  3.311422474487814  \\
            19.0  1.527537577790941  \\
            33.0  1.558058225580901  \\
            46.0  1.561983245327504  \\
            60.0  1.5626934760412974  \\
            76.0  1.5628050577107149  \\
            95.0  1.5628191109353409  \\
            122.0  1.5628211710290336  \\
            166.0  1.562821467389301  \\
            10431.0  0.7237012239998919  \\
            14428.0  0.11328078064763103  \\
            24037.0  0.018004956394026328  \\
            43599.0  0.0028750984569576588  \\
            79695.0  0.0004600644764334305  \\
        }
        ;
    \addlegendentry {$\nu = 0.6$}
    \addplot[color={rgb,1:red,0.7644;green,0.4441;blue,0.8243}, name path={12992f5a-3249-4c73-a62b-4d592baf711a}, draw opacity={1.0}, line width={1}, solid, mark={*}, mark size={3.75 pt}, mark repeat={1}, mark options={color={rgb,1:red,0.0;green,0.0;blue,0.0}, draw opacity={1.0}, fill={rgb,1:red,0.7644;green,0.4441;blue,0.8243}, fill opacity={1.0}, line width={0.75}, rotate={0}, solid}]
        table[row sep={\\}]
        {
            \\
            0.0  3.3491248847254083  \\
            1.0  3.3115155573816866  \\
            26.0  1.2500140649239189  \\
            49.0  1.5202573565988542  \\
            69.0  1.5542519143876634  \\
            92.0  1.561573164172576  \\
            116.0  1.5625785823592055  \\
            148.0  1.5627776687975063  \\
            194.0  1.5628059050220462  \\
            2882.0  0.6659174329428743  \\
            9268.0  0.044494199348455554  \\
            17647.0  0.0030342944308767716  \\
            30808.0  0.00023690803812193018  \\
        }
        ;
    \addlegendentry {$\nu = 0.4$}
\end{axis}
\end{tikzpicture}
        }
        \captionsetup{justification=centering}
        \caption{Suboptimality.\label{fig:gevp-gaussian-objective}}
    \end{subfigure}
    \begin{subfigure}[b]{0.33\textwidth}
        \centering
        \resizebox{\textwidth}{!}{
            \begin{tikzpicture}
\begin{axis}[
            width=3in, height=1.8in,
            at={(1.011in,0.642in)},
            legend cell align={left},
            legend pos=north east,
            scale only axis,
            minor grid style={thin,draw opacity=0.3},
            major grid style={thin,draw opacity=0.5},
            grid=both,
            xmin=0.2,
            xmax=1,
            ymin={1000},
            ymax={700000},
            xlabel = {Power $\nu$},
            ylabel = {\# Gradients}
]
    \addplot[color={rgb,1:red,0.0;green,0.6056;blue,0.9787}, name path={0fe18d87-cdca-4776-830b-e8c899f3fc83}, draw opacity={1.0}, line width={1}, solid, mark={*}, mark size={3.75 pt}, mark repeat={1}, mark options={color={rgb,1:red,0.0;green,0.0;blue,0.0}, draw opacity={1.0}, fill={rgb,1:red,0.0;green,0.6056;blue,0.9787}, fill opacity={1.0}, line width={0.75}, rotate={0}, solid}]
        table[row sep={\\}]
        {
            \\
            0.1  101409.0  \\
            0.2  203098.0  \\
            0.3  17554.0  \\
            0.4  30808.0  \\
            0.5  49295.0  \\
            0.6  79695.0  \\
            0.7  152939.0  \\
            0.8  277110.0  \\
            0.9  534008.0  \\
            1.0  577544.0  \\
        }
        ;
\end{axis}
\end{tikzpicture}
        }
        \captionsetup{justification=centering}
        \caption{\# gradients, $\varepsilon_\varphi = \varepsilon_A = 10^{-3}$.\label{fig:gevp-gaussian-powers}}
    \end{subfigure}
    \captionsetup{justification=centering}
    \caption{Comparison of the proposed power ALM with various powers $\nu$ on solving a representative GEVP with $n = 500$. The case $\nu = 1$ corresponds to the iALM from \cite{sahin_inexact_2019}.\label{fig:gevp-gaussian}}
\end{figure}

\begin{table}[h]
    \caption{
        Performance of power ALM with APGM inner solver on solving random GEVPs of size $n = 500$.
    }
    \label{table:gevp_gaussian}
    \centering
    
    \begin{adjustbox}{width=\textwidth}
    \setlength\extrarowheight{3pt}
    \pgfplotstabletypeset[%
        begin table={\begin{tabular}[t]},
        every head row/.style={
            before row={%
              \hline
              \vphantom{$q = 1$}\\
              \hline
            },
            after row/.add={}{\hline},
        },
        header=true,
        col sep=&,
        row sep=\\,
        string type,
        columns/{trial}/.style ={column name={trial}, column type={|c}},
        every row no 9/.style={after row=\hline},
        every row no 10/.style={after row=\hline},
    ]{
        \\
        trial\\
        1\\
        2\\
        3\\
        4\\
        5\\
        6\\
        7\\
        8\\
        9\\
        10\\
        avg.\\
    }%
    \begin{filecontents*}{data/gevp_gaussian.csv}
0.0006123259534016601	0.0009569562914475416	203098.0	8.953606965866889e-5	0.00023690803812193018	30808.0	0.00029374943174809154	0.0004600644764334305	79695.0	0.000571576750591829	0.0008933043686405462	277110.0	0.0002766073697592075	0.0004322932538018964	577544.0
0.0001642403582569063	0.0002547530833967304	103072.0	0.00010243662941822329	0.0003983542959100639	31540.0	0.00031265616676690033	0.00048734149184004316	102014.0	0.0006195448948845828	0.0009599248066798349	281648.0	0.00029863892294657024	0.0004627029733226262	607078.0
0.00016351812308468272	0.00025764225319035994	103162.0	0.00010392585880314442	0.00022900807382342236	31498.0	0.00032168944251076415	0.0005070313261907877	105864.0	0.00016411856814269044	0.0002580613848019997	390364.0	0.00031071897679590954	0.0004885409217263437	652055.0
0.0006137203839999827	0.0009758342616070781	203314.0	9.084751028531901e-5	0.00025316684151621693	31926.0	0.0002979364994041278	0.00047482126294906557	82429.0	0.0005823256518759923	0.0009259557605576507	282520.0	0.00028160970261970064	0.00044779069021561035	637788.0
0.0001682028306907446	0.0002622623375492239	103403.0	0.00010108842715339517	0.00031581402108993295	28301.0	0.00031176037778102295	0.0004875898296388126	88925.0	0.0006180619005216048	0.0009635710250788243	283888.0	0.00029797389204888614	0.0004645101542564145	617000.0
0.00018562906779784338	0.0002870238167200778	102763.0	0.00010547400937410689	0.0007020790660567933	36529.0	0.0003244616884657825	0.0005082855196734126	162261.0	0.00016709824145566543	0.00025900844293569314	387643.0	0.0003162404783584316	0.000490319771378811	658579.0
0.00018030119431300395	0.0002775609429053194	102932.0	4.7385756811513247e-5	8.558520835810768e-5	146417.0	0.0003174949396580473	0.0005007526426332731	206851.0	0.0006323581972866865	0.000977726757939168	376557.0	0.00030440728062952616	0.00047312152982348366	656967.0
0.00018428594490493833	0.0002849779744209524	103666.0	0.00010195857310590917	0.0002774456217748966	33362.0	0.00031207298766910974	0.0004828066403961273	102324.0	0.0006181441923025499	0.000953414587177015	302904.0	0.00029800655353562355	0.00045965054913787995	650897.0
0.0001554101011870257	0.0002443527701556647	103510.0	4.707148263005667e-5	0.0006710873209918589	132516.0	0.00030895403032926705	0.0005732243230023304	227561.0	0.0006095020213795488	0.0009492202742715783	435065.0	0.00029407954862514796	0.0004615996397672184	762997.0
0.0006153785174133342	0.0009596083764380481	203791.0	0.0001003183336103497	0.0002220567847632804	25497.0	0.00030291235538071337	0.0004731885251438861	86236.0	0.0005933902926960055	0.0009253467203869103	269632.0	0.00028669614653353825	0.00044708155812123174	620934.0
0.0003043012475050122	0.0004760972107830996	133271.1	8.900426508506865e-5	0.0003391505272406503	52839.4	0.0003103687919713827	0.0004955106037901169	124416.0	0.0005176120711137156	0.000806553412846922	328733.1	0.00029649788718525414	0.0004627611041551516	644183.9
\end{filecontents*}
\pgfplotstabletypeset[%
        begin table={\begin{tabular}[t]},
        every head row/.style={
        before row={%
          \hline
          \multicolumn{3}{|c|}{$\nu = 0.2$} & \multicolumn{3}{c|}{$\nu = 0.4$} & \multicolumn{3}{c|}{$\nu = 0.6$} & \multicolumn{3}{c|}{$\nu = 0.8$} & \multicolumn{3}{c|}{$\nu = 1.0$}\\
          \hline
        },
        after row/.add={}{\hline},
        },
        header=true,
       precision=1,
       columns/0/.style ={column name={const.\,viol.}, column type={|l}},
       columns/1/.style ={column name={subopt}, column type={c}},
       columns/2/.style ={column name={\# grads}, column type={c|}},
       columns/3/.style ={column name={const.\,viol.}, column type={c}},
       columns/4/.style ={column name={subopt.}, column type={c}},
       columns/5/.style ={column name={\# grads}, column type={c|}},
       columns/6/.style ={column name={const.\,viol.}, column type={c}},
       columns/7/.style ={column name={subopt}, column type={c}},
       columns/8/.style ={column name={\# grads}, column type={c|}},
       columns/9/.style ={column name={const.\,viol.}, column type={c}},
       columns/10/.style ={column name={subopt}, column type={c}},
       columns/11/.style ={column name={\# grads}, column type={c|}},
       columns/12/.style ={column name={const.\,viol.}, column type={c}},
       columns/13/.style ={column name={subopt}, column type={c}},
       columns/14/.style ={column name={\# grads}, column type={c|}},
        every row no 9/.style={after row=\hline},
        every row no 10/.style={after row=\hline},
        every row 10 column 3/.style={
                postproc cell content/.style={
                @cell content/.add={$\bf}{$}
                }
        },
        every row 10 column 4/.style={
                postproc cell content/.style={
                @cell content/.add={$\bf}{$}
                }
        },
        every row 10 column 5/.style={
                postproc cell content/.style={
                @cell content/.add={$\bf}{$}
                }
        },
    ]{data/gevp_gaussian.csv}
    \end{adjustbox}
    
\end{table}

\clearpage
\section{Inner solvers} \label{sec:inner-solvers}
For completeness, this section summarizes the \emph{universal problem-parameter free accelerated gradient} method \cite[Algorithm 1]{ghadimi_generalized_2019} method and the \emph{fast gradient method} \cite[Algorithm 3]{devolder_first-order_2014}, which are, respectively, used as inner solvers in \cref{sec:ipalm,sec:ippm}. 

First, we consider the ALM subproblem in \cref{alg:algorithm_nonconvex} \cref{alg:inner_inner_step}, which are of the form
\begin{equation} \label{eq:upfag-problem}
    \minimize_{x \in \bR^n} \Psi(x) \equiv \psi(x) + g(x).
\end{equation}
Here, \(\psi : \bR^n \to \bR \) is H\"older-smooth and \(g : \bR^n \to \exR \) is convex.
\Cref{alg:UPFAG} specifies the UPFAG method \citep[Algorithm 1]{ghadimi_generalized_2019} for problems of the form \eqref{eq:upfag-problem}.

\begin{algorithm}[ht]
\caption{Unified problem-parameter free accelerated gradient (UPFAG) method \citep{ghadimi_generalized_2019}}
\label{alg:UPFAG}
\begin{algorithmic}[1]
\Require $x_0 \in \bR^n$, \(\gamma_1, \gamma_2, \gamma_3 \in (0, 1)\), accuracy parameter \(\delta > 0\).
\State Set \(x_0^{\text{ag}} = x_0\), and \(\Lambda_0 = 0\).
\For{$k=1, 2, \dots$}
   \State Choose initial stepsize \(\hat \lambda_k > 0\) and find the smallest integer \(\tau_{1, k} \geq 0\) such that with \[
    \eta_k = \hat \lambda_k \gamma_1^{\tau_{1, k}}, \quad \lambda_k = \frac{\eta_k + \sqrt{\eta_k^2 + 4 \eta_k \Lambda_{k-1}}}{2}, \quad \alpha_k = \frac{\lambda_k}{\Lambda_k} \quad \text{and} \quad \Lambda_k = \sum_{i = 1}^{k} \lambda_i,
   \]
   \qquad the solutions obtained by 
   \begin{equation*}
    \begin{aligned}
        x_k^{\text{md}} &= (1 - \alpha_k) x_{k-1}^{\text{ag}} + \alpha_k x_{k-1},\\
        x_k &= \prox_{\lambda_k g}(x_{k-1} - \lambda_k \nabla \psi(x_k^{\text{md}}))\\
        \tilde x_k^{\text{ag}} &= (1 - \alpha_k) x_{k-1}^{\text{ag}} + \alpha_k x_k,
    \end{aligned}
   \end{equation*}
   \qquad satisfy the condition \[
    \psi(\tilde x_k^{\text{ag}}) \leq \psi(x_k^{\text{md}}) + \alpha_k \langle \nabla \psi(x_k^{\text{md}}), x_k - x_{k-1} \rangle + \frac{\alpha_k}{2 \lambda_k} \Vert x_k - x_{k-1} \Vert^2 + \delta \alpha_k.
   \]
   \State Choose initial stepsize \(\hat \beta_k > 0\) and find the smallest integer \(\tau_{2, k} \geq 0\) such that with \[
    \beta_k = \hat \beta_k \gamma_2^{\tau_{2, k}} \quad \text{and} \quad \bar x^{\text{ag}} = \prox_{\beta_k g}(x_{k-1}^{\text{ag}} - \beta_k \nabla \psi(x_{k-1}^{\text{ag}}))
   \]
   \qquad satisfy the condition \[
    \Psi(\bar x_k^{\text{ag}}) \leq \Psi(x_{k-1}^{\text{ag}}) - \frac{\gamma_3}{2 \beta_k} \Vert \bar x_k^{\text{ag}} - x_{k-1}^{\text{ag}}\Vert^2 + \frac{1}{k}.
   \]
   \State Choose \(x_k^{\text{ag}}\) such that \[
    \Psi(x_k^{\text{ag}}) = \min \left\{ \Psi(x_{k-1}^{\text{ag}}), \Psi(\bar x_k^{\text{ag}}), \Psi(\tilde x_k^{\text{ag}}) \right\}.
   \]
\EndFor
\end{algorithmic}
\end{algorithm}

Second, we consider the proximal point updates in \Cref{alg:algorithm_ippm} \cref{alg:inner_inner_step}, which are of the form
\begin{equation} \tag{\ref{eq:fgm-problem} rev.}
    \minimize_{x \in \cX} F(x).
\end{equation}
Here, \(F : \bR^n \to \bR \) is \((H_f, \nu)\)-H\"older-smooth and \(\rho\)-strongly convex with \(H_f \geq 0, \rho > 0\) and \(\nu \in (0, 1]\).
The set \(\cX \subseteq \bR^n\) is closed and convex.
\Cref{alg:fgm} specializes \cite[Algorithm 3]{devolder_first-order_2014} to problems of the form \eqref{eq:fgm-problem} for the standard Euclidean \textit{prox-function}.

\begin{algorithm}[ht]
\caption{Fast gradient method for strongly convex, H\"older smooth functions \citep{devolder_first-order_2014}}
\label{alg:fgm}
\begin{algorithmic}[1]
\Require $x_0 \in \cX$, accuracy parameter \(\delta > 0\), constants \(\rho > 0, H_f \geq 0, \nu \in (0, 1]\).
\State Set \(\mu = \rho\) and \(L = H_f \left( \frac{H_f}{2 \delta} \frac{1-\nu}{1+\nu} \right)^{\frac{1-\nu}{1+\nu}}\), and let \(\{\alpha_k\}_{k \geq 0}\) be such that \(\alpha_0 = 1\) and \(L+\mu \sum_{i = 0}^k \alpha_i = \frac{L \alpha_{k+1}^2}{\sum_{i = 0}^{k+1}\alpha_i}\).
\For{$k=0, 1, \dots$}
   \State \(y_k = \argmin_{x \in \cX} \left\{ \langle \nabla F(x_k), x - x_k \rangle + \frac{L}{2} \Vert x - x_k \Vert^2 \right\}\).
   \State \(z_k = \argmin_{x \in \cX} \left\{ \frac{L}{2} \Vert x - x_0 \Vert^2 + \sum_{k = 0}^k \alpha_i \left[ \langle \nabla F(x_i), x - x_i \rangle + \frac{\mu}{2} \Vert x - x_i \Vert^2 \right]\right\}\).
   \State \(x_{k+1} = \tau_k z_k + (1-\tau_k) y_k\) where \(\tau_k = \frac{\alpha_{k+1}}{\sum_{i = 0}^{k+1} \alpha_i}\).
\EndFor
\end{algorithmic}
\end{algorithm}
Note that the first step is a standard projected gradient update \(
    y_k = \proj_{\cX}(x_k - \frac{1}{L} \nabla F(x_k)),
\)
whereas the second step can be written as \[
    z_k = \proj_{\cX} \left( \frac{L x_0 + \mu \sum_{i = 0}^k \alpha_i x_i - \sum_{i = 0}^k \alpha_i \nabla F(x_i)}{L+\mu \sum_{i = 0}^k \alpha_i}\right).
\]
We highlight that per iteration \cref{alg:fgm} only requires a single gradient evaluation of $F$ (at \(x_k\)), and two projections onto \(\cX\).

\end{document}